\pdfoutput=1
\RequirePackage{ifpdf}
\ifpdf % We~are running pdfTeX in pdf mode
\documentclass[pdftex]{sigma}
\else
\documentclass{sigma}
\fi

\usepackage[all, 2cell]{xy}
\UseTwocells

\usepackage{mathtools}
\usepackage{tikz-cd}
\usepackage{eucal}

\numberwithin{equation}{section}

\newtheorem{Theorem}{Theorem}[section]
\newtheorem*{Theorem*}{Theorem}
\newtheorem{Corollary}[Theorem]{Corollary}
\newtheorem{Lemma}[Theorem]{Lemma}
\newtheorem{Proposition}[Theorem]{Proposition}
{ \theoremstyle{definition}
	\newtheorem{Definition}[Theorem]{Definition}
	\newtheorem{Notation}[Theorem]{Notation}
	
	\newtheorem{Example}[Theorem]{Example}
	\newtheorem{Remark}[Theorem]{Remark} }

\newtheorem{thmx}{Theorem}

\newcommand{\sslash}{/\!\!/}

\newcommand{\hamp}{/\!\!/_{\psi}\,}
\newcommand{\hamn}{/\!\!/_{\psi^*}\,}

\newcommand{\lhamp}{N^L_+\, \backslash\!\!\backslash_{-\psi}\,}
\newcommand{\lhamn}{N^L_-\, \backslash\!\!\backslash_{-\psi^*}\,}

\newcommand{\C}{\mathbf{C}}
\newcommand{\Z}{\mathbf{Z}}

\newcommand{\bP}{\mathbb{P}}

\newcommand{\dD}{\mathrm{D}_{\hbar}}
\newcommand{\GL}{\mathrm{GL}}

\newcommand{\Mat}{\mathrm{Mat}}
\newcommand{\End}{\mathrm{End}}
\newcommand{\Id}{\mathrm{Id}}
\newcommand{\U}{\mathrm{U}_{\hbar}}
\newcommand{\qdet}{\mathrm{qdet}}
\newcommand{\Rep}{\mathrm{Rep}}
\newcommand{\HC}{\mathrm{HC}_{\hbar}}
\newcommand{\ad}{\mathrm{ad}}
\newcommand{\gen}{\mathrm{gen}}
\newcommand{\free}{\mathrm{free}}
\newcommand{\res}{\mathrm{res}}
\newcommand{\Hom}{\mathrm{Hom}}
\newcommand{\QCoh}{\mathrm{QCoh}}
\newcommand{\BiMod}[2]{{}_{#1}\mathrm{BMod}_{#2}}
\newcommand{\rY}{\mathrm{Y}}
\newcommand{\miura}{\mathrm{Miura}}
\newcommand{\gklo}{\mathrm{GKLO}}
\newcommand{\Sym}{\mathrm{Sym}}
\newcommand{\aSym}{\mathrm{Sym}^-}
\newcommand{\rat}{\mathrm{rat}}
\newcommand{\Hecke}{\mathrm{Hecke}}
\newcommand{\Mod}{\mathrm{Mod}}
\newcommand{\Gr}{\mathrm{Gr}}
\newcommand{\IC}{\mathrm{IC}}

\newcommand{\cO}{\mathcal{O}}
\newcommand{\cF}{\mathcal{F}}
\newcommand{\cA}{{}^N\!\tilde{\mathcal{A}}_{\hbar}}
\newcommand{\cM}{\mathcal{M}}
\newcommand{\cS}{\mathcal{S}}
\newcommand{\cQ}{\mathcal{Q}}

\newcommand{\cV}{\mathcal{V}}

\newcommand{\g}{\mathfrak{g}}
\renewcommand{\sl}{\mathfrak{sl}}
\newcommand{\gl}{\mathfrak{gl}}
\newcommand{\n}{\mathfrak{n}}
\newcommand{\h}{\mathfrak{h}_N}

\newcommand{\m}{\mathfrak{m}}
\newcommand{\fH}{\mathfrak{H}}

\newcommand{\bu}{\mathbf{u}}

\newcommand{\pP}{P_{\psi}}

\newcommand{\npP}{P_{\psi}}

\newcommand{\adj}[2]{
\xymatrix{
#1 \ar@<.5ex>[r] & #2 \ar@<.5ex>[l]
}
}

\begin{document}
	
\allowdisplaybreaks

\newcommand{\arXivNumber}{2505.09520}

\renewcommand{\PaperNumber}{017}

\FirstPageHeading

\ShortArticleName{Shuffle Products, Degenerate Affine Hecke Algebras, and Quantum Toda Lattice}

\ArticleName{Shuffle Products, Degenerate Affine Hecke Algebras,\\ and Quantum Toda Lattice}

\Author{Artem KALMYKOV}

\AuthorNameForHeading{A.~Kalmykov}

\Address{Department of Mathematics and Statistics, McGill University, \\
805 Sherbrooke Street West, Montreal, Quebec H3A 0B9, Canada}
\Email{\mail{artem.o.kalmykov@gmail.com}}

\ArticleDates{Received August 14, 2025, in final form February 09, 2026; Published online February 23, 2026}

\Abstract{We revisit an identification of the quantum Toda lattice for $\mathrm{GL}_N$ and the truncated shifted Yangian of $\mathfrak{sl}_2$, as well as related constructions, from a purely algebraic point of view, bypassing the topological medium of the homology of the affine Grassmannian. For instance, we interpret the Gerasimov--Kharchev--Lebedev--Oblezin homomorphism into the algebra of difference operators via a finite analog of the Miura transform. This algebraic identification is deduced by relating degenerate affine Hecke algebras to the simplest example of a rational Feigin--Odesskii shuffle product. As a~bonus, we obtain a presentation of the latter via a mirabolic version of the Kostant--Whittaker reduction.}


\Classification{17B37; 17B38}

\section{Introduction}

Shuffle algebras are ubiquitous in representation theory. To name a few applications, they are related to infinite-dimensional quantum groups (including quantum toroidal algebras), Hall algebras, canonical bases, and so on. They were introduced by Feigin--Odesskii in a series of papers~\cite{FeiginOdesskii, FeiginOdesskiiSklyanin,OdesskiiFeiginEllipticDeformation}. In the simplest case, a shuffle algebra is, as a vector space, the direct sum \smash{$\bigoplus_k \cO(Y_1,\dots,Y_k)^{S_k}$} of suitably defined symmetric functions in $k$ variables, with the product
\begin{gather}
\label{eq::intro_shuffle_prod}
	\cO(Y_1,\dots,Y_k)^{S_k} \times \cO(Y_1,\dots,Y_l)^{S_l} \rightarrow \cO(Y_1,\dots,Y_{k+l})^{S_{k+l}}, \\
	(f*g)(Y_1,\dots,Y_{k+l}) = \frac{1}{(k+l)!} \Sym_{k+l} \biggl( \prod_{i<j}\! \frac{\theta(Y_i- Y_j + 1)}{\theta(Y_i\!-\!Y_j)} f(Y_1,\dots,Y_k) g(Y_{k+1}, \dots,Y_{k+l})\!\biggr),\nonumber
\end{gather}
where $\Sym_{k+l}$ is the symmetrization operator in $k+l$ variables, and $\theta(u)$ is one of the functions%
\begin{equation}
\label{eq::intro_trichotomy}
	\theta^{\mathrm{ell}} (u)= \vartheta_1(u|\tau), \qquad \theta^{\mathrm{trig}}(u) = \sin(\pi u/\tau), \qquad \theta^{\rat}(u) = u,
\end{equation}
corresponding to the elliptic, trigonometric, and rational cases, for some parameter $\tau$. Here $\vartheta_1(u|\tau)$ is the first Jacobi theta-function. For instance, in the original construction of loc.\ cit., the function space is
\[
	\cO(Y_1,\dots,Y_k)^{S_k} = {}^N \!\cS^{\mathrm{ell}}_k := S^k H^0 (\Theta), \qquad {}^N \!\cS^{\mathrm{ell}} = \bigoplus_k {}^N \!\cS^{\mathrm{ell}} _k,
\]
the sum of symmetric powers of global sections of a degree $N$ line bundle on an elliptic curve (put differently, degree $N$ theta functions), and the product is determined by the factor $\theta^{\mathrm{ell}}(u)$, up to shift of variables.

In this paper, we study the rational case. The elliptic one above admits a degeneration to it~\cite{EndelmanHodgesTwistedSU} so that
\begin{equation*}
%\label{eq::intro_polynomial_N}
		{}^N\!\cS^{\rat} = \bigoplus_k {}^N\!\cS_k^{\rat}, \qquad {}^N\!\cS^{\rat}_k = \C[Y_1,\dots,Y_k]^{S_k}_{<N},
\end{equation*}
where \smash{$\C[Y_1,\dots,Y_k]^{S_k}_{<N}$} is the space of symmetric polynomials of degree less than $N$ in each variable, and the factor is $\theta^{\rat}(u)$. More generally, one can drop the restriction on the degree and consider \smash{$\cS^{\rat} = \bigoplus_k \cS^{\rat}_k$} with \smash{$\cS^{\rat}_k = \C[Y_1,\dots,Y_k]^{S_k}$}.

The first goal of the paper is to relate these shuffle algebras to \emph{degenerate affine Hecke algebras} whose definition we recall in Section~\ref{subsect:deg_affine_hecke}. The idea is as follows: observe that without the theta factor, the shuffle product is simply the multiplication on the symmetric algebra of the vector space $\C[Y_1]$ or its $N$-th truncation. Instead of deforming it via inserting a factor, we can deform the symmetrization itself, that is, the action of the symmetric group. One such deformation is provided by the \emph{Demazure operators}, so that a simple transposition acts on $\C[Y_1,\dots,Y_k]$ by
\[
	\sigma_i = P_{i,i+1} + \frac{1-P_{i,i+1}}{Y_i - Y_{i+1}},
\]
where $P_{i,i+1}$ permutes variables $Y_i$ and $Y_{i+1}$. Denoting by $e_k \in \C[S_k]$ the symmetrizer idempotent in the group algebra of the symmetric group, this is the action of $S_k$ on the spherical module $\fH_k e_k \cong \C[Y_1,\dots,Y_k]$ of the degenerate affine Hecke algebra $\fH_k$. Its spherical subalgebra $e_k \fH_k e_k$ can be identified with the symmetric polynomials $\cS^{\rat}_k$.

Based on this observation, one can introduce an a priori different shuffle-like product on the direct sum \smash{$\bigoplus_k \cS^{\rat}_k$} via, as we mentioned, symmetrization with respect to this deformed action:
\[
	(f,g) \mapsto f *^{\Hecke} g:= e_{k+l} f(Y_1,\dots,Y_k) g(Y_{k+1},\dots,Y_{k+l})e_{k+l} \in e_{k+l} \fH_{k+l} e_{k+l} \cong \cS^{\rat}_{k+l}.
\]

The following result of Theorem~\ref{thm:hecke_vs_fo_shuffle} seems to be folklore, but, unfortunately, we were not able to find a reference for its algebraic proof. A geometric proof can be found in an unpublished preprint by Grojnowski \cite{Grojnowski}.
\begin{thmx}
	The rational Feigin--Odesskii shuffle product coincides with the one induced from degenerate affine Hecke algebras.
\end{thmx}

In principle, the theorem follows directly from the explicit formula for the symmetrizer in the spherical representation of a degenerate affine Hecke algebra and, strictly speaking, is not new: for instance, in the trigonometric case, it reflects a well-known relation between Hall--Littlewood polynomials and affine Hecke algebras \cite{Macdonald}. Our main contribution is an alternative proof based on the following simple observation: a canonical \emph{affinization} procedure \cite{KSQuantumGroups}, which, in general, produces a solution of a \emph{quantum Yang--Baxter equation}
\[
	R_{12}(z_1-z_2) R_{13}(z_1-z_3) R_{23}(z_2-z_3) = R_{23}(z_2-z_3) R_{13}(z_1-z_3) R_{12}(z_1-z_2),
\]
out of a representation of symmetric group (or, in the trigonometric case, of a finite Hecke algebra), applied to the spherical module gives a particular example of a \emph{Shibukawa--Ueno $R$\nobreakdash-operator} introduced in \cite{ShibukawaUeno} whose definition we recall in Section~\ref{subsect:shibukawa_ueno}. In general, the latter acts on the space of suitably defined functions and takes as an input any solution to the three-term relation~\eqref{eq::three_term_relation}; in particular, the case of the spherical module of a degenerate affine Hecke algebra corresponds to the rational factor in~\eqref{eq::intro_trichotomy}, while that of a non-degenerate affine Hecke algebra~-- to the trigonometric one. Then the theorem is deduced by, on the one hand, representing the symmetrizer though a certain combination of $R$-matrices, on the other hand, expressing this combination as the usual symmetrization map, but with a factor akin to~\eqref{eq::intro_shuffle_prod}, see Proposition~\ref{prop:su_symmetrizer}.

The main advantage of this proof is that the second part works uniformly for any factor that satisfies the three-term relation. In principle, there are (formal) analogs of theta functions for arbitrary generalized cohomology theory with orientation \cite{Segal}. Unfortunately, we do not know how the respective Shibukawa--Ueno $R$-operators are related to the corresponding \emph{formal (affine) Hecke algebras} \cite{HMSZ}, however, we believe that similar formulas should hold in more general cases.

In particular, it would be interesting to study the elliptic version of the Shibukawa--Ueno operator and its potential relation to \emph{elliptic affine Hecke algebras} \cite{GinzburgKapranovVasserot}. Curiously, the elliptic case admits a dual description on the $\GL_N$-side, in the sense of \emph{Fourier--Mukai transform} on an elliptic curve. Namely, by \cite{FelderPasquier}, the elliptic Shibukawa--Ueno $R$-operator is gauge-equivalent to \emph{Belavin's $R$-matrix} \cite{Belavin,CherednikRmatrix} whose quasi-periodic properties can be interpreted via an essentially unique stable vector bundle of rank $N$ and degree 1 on an elliptic curve. A natural question is: what is the meaning of the shuffle product on the $\GL_N$-side?

We approach this question for the rational case in the second part of the paper where we apply the construction of the first part to the \emph{quantum Toda lattice} for $\GL_N$. The idea to use degenerate affine Hecke algebras is inspired by the results of \cite{CautisKamnitzer} where the authors constructed and proved a~K-theoretic version of the geometric Satake equivalence for a certain subcategory of the $q$-analog of Harish-Chandra bimodules (we recall the rational versions below). Following the Kazhdan--Kostant formulation of this integrable system \cite{Kostant}, let us consider a non-degenerate character $\psi$ of the Lie algebra $\n_+$ of the subgroup of upper unitriangular matrices $N_+$. The action of $\n_+$ via right- and left-invariant vector fields gives rise to a homomorphism ${\mathrm{U}(\n_+) \otimes \mathrm{U}(\n_+) \rightarrow \mathrm{D}(\GL_N)}$ into the algebra of differential operators. Then the Toda Hamiltonians are realized as the image of bi-invariant differential operators, which can be identified with the center~\smash{$\Z_{\gl_N}$ of~$\mathrm{U}(\gl_N)$}, inside the quantum Hamiltonian reduction \smash{$\lhamp \mathrm{D}(\GL_N) \hamp N_+^R$} with respect to the latter embedding and the character $(-\psi,\psi)$. This is an example of the \emph{Kostant--Whittaker reduction} for Harish-Chandra bimodules studied in \cite{BezrukavnikovFinkelberg}.

According to \cite[Appendix B]{BravermanFinkelbergNakajima}, the Toda lattice can be identified with a certain quotient \smash{$\rY^0_{-2N}(\sl_2)$} of the \emph{shifted Yangian} \smash{$\rY_{-2N}(\sl_2)$} whose definition we recall in Section~\ref{subsect:shifted_yangian}. For introduction, we will only use the fact that it is generated by the coefficients of the series
\[
	x_+(u) = \sum_{i=1}^{\infty} x_+^{(i)} u^{-i}, \qquad x_-(u) = \sum_{i=1}^{\infty} x_-^{(i)} u^{-i}, \qquad h(u) = u^{-2N} + \sum_{i=2N+1}^{\infty} h^{(i)} u^{-i},
\]
and that the coefficients of $x_+(u)$ generate the subalgebra isomorphic to the shuffle algebra $\cS^{\rat}$~\cite{TsymbaliukPBWD}. Therefore, it gives a $\GL_N$-presentation of $\cS^{\rat}$. However, citing \cite{FKPRW}, ``we must admit that the identification of the quantum Toda [lattice] for~$\GL_N$ and the [truncated] shifted Yangian for~$\sl_2$ (purely algebraic objects) goes through a topological medium: equivariant homology of the affine Grassmannian of $\GL_N$''; this additional step is the \emph{derived Satake equivalence} of \cite{BezrukavnikovFinkelberg}. For instance, while $\rY_{-2N}(\sl_2)$ is defined via generators and relations, this identification does not specify explicitly the images of the former.

The goal of the second part of the paper is to apply the results of the first part to construct and prove an isomorphism
\[
\rY^{0}_{-2N}(\sl_2) \xrightarrow{\sim} \lhamp \mathrm{D}(\GL_N) \hamp N_+^R,
\]
as well as related constructions, in a purely algebraic and explicit way, bypassing the topological medium, and ``refine'' it to a \emph{mirabolic} version of the Kostant--Whittaker reduction. It is based on another classical fact. For any $\GL_N$-representation $V$, the free left $\mathrm{U}(\gl_N)$-module $\mathrm{U}(\gl_N) \otimes V$ has a right $\mathrm{U}(\gl_N)$ action defined from the diagonal group action of $\GL_N$ on the tensor product. This is an example of a \emph{(free) Harish-Chandra bimodule} whose definition and properties we~review in~Section~\ref{sect:hc_bimodules}. For the tensor product of the vector representation \smash{$V = \bigl(\C^N\bigr)^{\otimes k}$}, there is an~action of the degenerate affine Hecke algebra $\fH_k$ on \smash{$\mathrm{U}(\gl_N) \otimes \bigl(\C^N\bigr)^{\otimes k}$} respecting the $\mathrm{U}(\gl_N)$-bimodule structure. In particular, its spherical subalgebra $\cS_k \cong e_k \fH_k e_k$ acts on the $S_k$-invariants that can be identified with the free Harish-Chandra bimodule $\mathrm{U}(\gl_N) \otimes S^k \C^N$ associated to the $k$-th symmetric power. By restricting the action on the highest-weight vector for all $k$, one obtains a~map of vector spaces
\[
	\cS^{\rat} =\bigoplus_k \cS^{\rat}_k \rightarrow \mathrm{U}(\gl_N) \ltimes S^{\bullet} \C^N \hamp N_+,
\]
where the target is the quantum Hamiltonian reduction of the semi-direct product with $\gl_N$ acting on the symmetric algebra $S^{\bullet} \C^N$ in the standard way. In particular, it is an algebra as well. If we denote by $\m_N$ the \emph{mirabolic subalgebra} of the subgroup of $\GL_N$ preserving the dual of the last basis vector in $\C^N$ (see Section~\ref{subsect:mirabolic} for precise definition), then it turns out that this map restricts to a similarly defined quantum Hamiltonian reduction
\begin{equation}
\label{eq::intro_shuffle_vs_mirabolic_kw}
	{}^N\!\cS^{\rat} \rightarrow \mathrm{U}(\m_N) \ltimes S^{\bullet} \C^N \hamp N_+.
\end{equation}
The following is Theorems~\ref{thm:kw_to_shuffle_mirabolic} and~\ref{thm:kw_to_yangian}; here \smash{$\rY_{-2N}(\sl_2)^{\geq 0}$} is the non-negative part of the shifted Yangian. In fact, the explicit formulas for the second homomorphism are close to another presentation of the quantization of Zastava spaces for $\sl_2$ from \cite[Section 3.22]{FinkelbergRybnikov}.
\begin{thmx}
	The map \smash{${}^N \!\cS^{\rat} \xrightarrow{\sim} \mathrm{U}(\m_N) \ltimes S^{\bullet} \C^N \hamp N_+$} is an algebra isomorphism. It can be extended to a surjective algebra homomorphism \smash{$\rY_{-2N}(\sl_2)^{\geq 0} \rightarrow \mathrm{U}(\gl_N) \ltimes S^{\bullet} \C^N \hamp N_+$}.
\end{thmx}

Unfortunately, we do not know if and how this result can be related to the aforementioned stable bundles on an elliptic curve. We expect (or, rather, hope) that it could viewed as a~quantum rational version of a certain algebra of Hecke modifications of semistable bundles of degree~0, see Remark~\ref{remark:mirabolic_vs_stable_vector_bundle}. In general, it would be interesting to interpret this result at least on the classical level in terms of elliptic analog of the Kostant slice \cite{FriedmanMorgan}.

Besides degenerate affine Hecke algebras, the other main tool to prove the theorem is the \emph{Kirillov projector} $\pP$ introduced in \cite{KalmykovYangians}. This is an element in a certain completion of $\mathrm{U}(\m_N)$ that acts on Whittaker modules and projects onto the space of Whittaker vectors (that is, on which $\n_+$ acts with character $\psi$); as such, it is an analog of the classical \emph{extremal projector} \cite{AsherovaSmirnovTolstoy} in the case of trivial character. For instance, the Kirillov projector gives a canonical isomorphism $V \xrightarrow{\sim} \mathrm{U}(\m_N) \otimes V \hamp N_+$ of any $\GL_N$-representation $V$ with the mirabolic version of the Kostant--Whittaker reduction. It turns out that this isomorphism has rather peculiar properties with respect to the degenerate affine Hecke algebra action mentioned above, see Section~\ref{subsect:mirabolic}. For example, under the isomorphism \eqref{eq::intro_shuffle_vs_mirabolic_kw}, the standard monomial basis of $S^{\bullet} \C^N$ is mapped to a~rational version of the Hall--Littlewood polynomials, see Example~\ref{example:mirabolic_symmetric_power}.

We extend this homomorphism to the Toda lattice. Naturally, $\GL_N$ is embedded into the matrix space; denote by $X = (x_{ji})$ the (transposed) matrix with $x_{ij} \in \cO(\GL_N)$ the $(i,j)$-matrix entry function. Its inverse $X^{-1}$ is well-defined in $\cO(\GL_N)$. Also, denote by \smash{$E^R = \bigl(E_{ij}^R\bigr)$} the matrix composed of left-invariant vector fields associated to the standard generators of the Lie algebra $\gl_N$. The following is Theorem~\ref{thm:toda_lattice_yangian}. It would interesting to compare this explicit form with the approach of \cite{FKPRW} via coproducts on both sides.
\begin{thmx}
	There is an isomorphism $\rY^0_{-2N}(\sl_2) \xrightarrow{\sim} \lhamp \mathrm{D}(\GL_N) \hamp N_+^R$ defined on generating currents by
	\begin{align*}
		&x_+(u)\mapsto -\sum_{i=0}^{\infty} \bigl[\bigl(E^R\bigr)^i X\bigr]_{1N} u^{-i-1}, \qquad x_-(u) \mapsto \sum_{i=0}^{\infty} \bigl[X^{-1} \bigl(E^R\bigr)^i \bigr]_{1N} u^{-i-1}, \\
		&h(u)\mapsto \frac{1}{A(u)A(u-1)},
	\end{align*}
	where the subscript $(1N)$ means the corresponding entry in the matrix product and $A(u)$ is a~\emph{quantum determinant} whose coefficients generate $\Z_{\gl_N}$.
\end{thmx}

One of the main ideas of the proof is to interpret the \emph{GKLO homomorphism} from $\rY_{-2N}(\sl_2)$ into a certain localization ${}^N\!\tilde{\mathcal{A}}$ of the algebra of difference operators, introduced in~\cite{GKLO} for non-shifted Yangians and later studied in \cite{BravermanFinkelbergNakajima, KWWY} in general case, via a finite analog of the \emph{Miura transform} \cite{GinzburgKazhdan}, see Section~\ref{sect:miura}. In general, for any Harish-Chandra bimodule $X$, it relates the Kostant--Whittaker reduction $X \hamp N_+$ to its version with the trivial character, the \emph{parabolic restriction} $X^{\gen} \sslash N_-$, which we review in Section~\ref{sect:parabolic_reduction} (here, the superscript ``gen'' corresponds to a certain localization), akin to the Harish-Chandra homomorphism $\Z_{\gl_N} \rightarrow \mathrm{U}(\h) \cong (\n_- \backslash \mathrm{U}(\gl_N))^{N_-}$. In~particular, it gives an injective homomorphism
\[
	\lhamp \mathrm{D}(\GL_N) \hamp N^R_+ \rightarrow \lhamp \mathrm{D}^{\gen}(\GL_N) \sslash N_-^R.
\]
The following is Theorem~\ref{thm:diff_ops_vs_difference} and the second part of Theorem~\ref{thm:toda_lattice_yangian}.
\begin{thmx}
There is an isomorphism ${}^N\!\tilde{\mathcal{A}} \cong \lhamp \mathrm{D}^{\gen}(\GL_N) \sslash N_-^R$ such that the composition
	\[
		\rY^0_{-2N} (\sl_2) \xrightarrow{\sim} \lhamp \mathrm{D}(\GL_N) \hamp N^R_+ \rightarrow \lhamp \mathrm{D}^{\gen}(\GL_N) \sslash N_-^R \xrightarrow{\sim} {}^N\!\tilde{\mathcal{A}}
	\]
	coincides with the GKLO homomorphism.
\end{thmx}

In Section~\ref{subsect:monopole}, we apply this theorem to get simplest cases of the \emph{monopole operators}, which are rational versions of certain degenerations of Macdonald difference operators, as the images of particular minor functions on $\GL_N$. We express the latter in terms of the shuffle algebra similarly to the geometric and $q$-setting of \cite{FinkelbergTsymbaliuk} and \cite{TsymbaliukGKLO}. Finally, in Section~\ref{subsect:comparison}, we show that this isomorphism coincides with the one in
\cite[Appendix B]{BravermanFinkelbergNakajima} under the derived Satake isomorphism of \cite{BezrukavnikovFinkelberg} between \smash{$\lhamp \mathrm{D}(\GL_N) \hamp N_+^R$} and the equivariant Borel--Moore homology~\smash{$H^{\GL_N[[t]] \rtimes \C^*}_{\bullet}(\Gr)$} of the affine Grassmannian of $\GL_N$.

We would like to mention that we actively use the formulas (adapted to setting of the paper) from the paper \cite{MolevTransvectorAlgebras} that studied Olshanski centralizer construction from the point of view of \emph{Mickelsson algebras} whose particular case we recall in Section~\ref{sect:parabolic_reduction}. It would be interesting to find a direct conceptual link between two subjects.

\subsection*{Future directions} The results of the paper should be more or less directly generalizable to the trigonometric, or $q$-deformed, case: for instance, see \cite{FinkelbergTsymbaliuk} for geometric approach, \cite{CautisKamnitzer} for a certain version of $K$-theoretic Satake equivalence, and \cite{KalmykovSafronov} and references therein for quantum Harish-Chandra bimodules and the corresponding parabolic restriction. In fact, the only missing ingredient is a~$q$-analog of the Kostant--Whittaker reduction which, to the best of our knowledge, was not yet defined; however, see \cite{SevostyanovWhittakerModules} for the $q$-analog of the quantum Toda lattice.

Our original motivation, however, was the elliptic version, where various representation-theoretic constructions should acquire much clearer algebro-geometric meaning, but, unfortunately, much less is known. For instance, to the best of our knowledge, there is no definition of ``elliptic'' Harish-Chandra bimodules, although there is a strong evidence that such a category should exist, see the discussion in the introduction of \cite{KalmykovSafronov}. We believe that the results of this paper can give a hint how to approach this problem at least in type A. For example, as the results of Section~\ref{sect:kw_reduction} combined with calculations with Shibukawa--Ueno $R$-operators of Section~\ref{sect:deg_aha_yangians} suggest, the action of degenerate affine Hecke algebras plays a crucial role to study the category of Harish-Chandra bimodules, therefore, it is reasonable to search for their elliptic generalization in terms of representations of elliptic affine Hecke algebras \cite{ZhaoZhong}. In general, as Remark~\ref{remark:mirabolic_vs_stable_vector_bundle} indicates, it would be interesting to study further this interplay between representation theory and, in a certain sense, ``quantization'' of various algebro-geometric constructions.

\subsection*{Structure of the paper}
In Section~\ref{sect:prelim}, we introduce the necessary notations concerning $\GL_N$ and differential operators $\mathrm{D}(\GL_N)$ as well quantum minors. In Section~\ref{sect:deg_aha_yangians}, we recall basic properties of shifted Yangians, degenerate affine Hecke algebras, and relate the latter to Shibukawa--Ueno $R$-operators and shuffle algebras. In Section~\ref{sect:hc_bimodules}, we recall the definition of Harish-Chandra bimodules and degenerate affine Hecke algebras action on the free ones associated to tensor powers of the vector representation. In Section~\ref{sect:parabolic_reduction}, we study quantum Hamiltonian reduction with a trivial character, the parabolic restriction. In Section~\ref{sect:kw_reduction}, we study the Kostant--Whittaker reduction and its mirabolic version relating them to degenerate affine Hecke algebras, shuffle algebras, and (the positive part of) shifted Yangians. In Section~\ref{sect:miura}, we relate two reductions via finite Miura transform and interpret it as a GKLO map. In~Section~\ref{sect:toda}, we prove an isomorphism between truncated shifted Yangians and Toda lattice, compare it with the geometric identification via homology of the affine Grassmannian, and interpret monopole operators in terms certain minor functions. In~Appendix~\ref{appendix:technical_lemmas}, we prove some technical lemmas from the main sections. While in the main body of the paper we deal with the vector representation, in Appendix~\ref{subsect:dual_rep}, we present the corresponding facts for the dual that are needed for the full Toda lattice.

For simplicity of exposition, we set $\hbar=1$ in the introduction.

\section*{List of notations}

\begin{tabular}{@{}c | c@{}}
	$v_i$ & natural basis of $\C^N$ \\
	$\phi_i$ & dual basis of $\C^N$ \\
	$X = (x_{ji})$ & the \textit{transposed} matrix of matrix entries functions on $\GL_N$ \\
	$\Z_{\gl_N}$ & center of $\U(\gl_N)$ \\
	$A(u)$ & generating polynomial of $Z_{\gl_N}$, certain quantum determinant \\
	$\fH^{\kappa}_k$ & degenerate affine Hecke algebra with quantization parameter $\kappa = \pm \hbar$ \\
	$Y_i$ & polynomial generators of $\fH^{\kappa}_k$ \\
	$\sigma_i$ & simple transposition in $\C[S_k] \subset \fH_k^{\kappa}$ \\
	$e_k$ & symmetrizer in $\C[S_k]$ \\
	$e_k^-$ & antisymmetrizer in $\C[S_k]$ \\
	$ \sigma^{\kappa}_i$ & action of $\sigma_i$ on spherical module $\C[Y_1, \dots, Y_k]$ of $\fH^{\kappa}_k$ \\
	$\C[\hbar][Y_1, \dots, Y_k]_{<N}$ & polynomials of degree less than $N$ in each variable \\
	$\Omega_i$ & action of $Y_i\in \fH^{\hbar}_k$ on Harish-Chandra bimodule $\U(\gl_N) \otimes \bigl(\C^N\bigr)^{\otimes k}$ \\
	$w_i$ & shifted coordinates $E_{ii} + N\hbar - \hbar$ on $\h^* \times \mathbb{A}^1$ \\
	$\U(\h)^{\gen}$ & localized ring of functions on $\h^* \times \mathbb{A}^1$ \\
	$P$ & extremal projector of $\gl_N$ \\
	$\cA$ & algebra of localized difference operators $\U(\h)^{\gen}\bigl[\bu_1^{\pm},\dots,\bu_N^{\pm}\bigr]$ \\
	$\cA^-$ & negative part of difference operators $\U(\h)^{\gen}\bigl[\bu_1^{- 1},\dots,\bu_N^{- 1}\bigr]$ \\
	$\psi$ & non-zero number in $\C$ \\
	$\psi$ & character of $\n_+$ such that $\psi(E_{ij}) = \delta_{i+1,j} \psi$ \\
	$\m_N$ & mirabolic subalgebra \\
	$()^R, ()^L$ & anything related to right (resp.\ left) action of $\GL_N$ on itself \\
	${}^R \iota_{\alpha}$ and analogs & map $\U\bigl(\gl^R_N\bigr) \ltimes S^{\bullet}\C^N \rightarrow \mathrm{D}_{\hbar} (\GL_N)$ or their reductions \\
	$\hat{T}(u)$ & quantum comatrix of $T(u) = u + E$ \\
	$v_{\lambda}, v_{\lambda}^*$ & highest-weight vector in irreducible representation $V_{\lambda}$ (resp.\ its dual $V_{\lambda}^*$) \\
	$\Gr$ & affine Grassmannian $\GL_N((t))/ \GL_N[[t]]$ \\
	$\Gr^{\lambda}$ & $\GL_N[[t]]$-orbit of $t^{\lambda} \in \Gr$ \\
	$\Omega_i^*$ & action of $Y_i \in \fH_k^{-\hbar}$ on Harish-Chandra bimodule $\U(\gl_N) \otimes \bigl(\bigl(\C^N\bigr)^*\bigr)^{\otimes k}$
\end{tabular}
	
\section{Preliminaries}
\label{sect:prelim}

\subsection[Notations for GL(N)]{Notations for $\boldsymbol{\GL_N}$}
\label{subsect:notations}

In this subsection, we fix some notations regarding the group $\GL_N$.

Let $\C^N$ be the vector space with a basis $\{v_1,\dots, v_N\}$. Denote by $\{\phi_1,\dots,\phi_N\}$ the dual basis of~\smash{$\bigl(\C^N\bigr)^*$} such that $\langle v_i, \phi_j \rangle = \delta_{ij}$. Let $\Mat_{N\times N}$ be the affine space of $N\times N$ matrices. We identify its space of functions with the polynomial algebra via
\[
	\cO(\Mat_{N\times N}) \cong \C[x_{ij} \mid i,j=1,\dots, N], \qquad x_{ij} (A) = \langle A\cdot v_j,\phi_i\rangle,\quad A\in \Mat_{N\times N},
\]
and we combine the generators into a matrix
\begin{equation*}
	X = (x_{ji})_{i=1, \dots, N}^{j = 1, \dots, N} \in \C[x_{ij}] \otimes \End\bigl(\C^N\bigr)
\end{equation*}
(observe the transposition).

Let $\GL_N \subset \Mat_{N\times N}$ be the group of invertible matrices. Its ring of functions is the localization of $\cO(\Mat_{N\times N})$ by the determinant function, in particular, we have an embedding $\C[x_{ij}] \hookrightarrow \cO(\GL_N)$. We call the image the \emph{polynomial functions}. In addition, there are entries of the matrix inverse $X^{-1} \in \cO(\GL_N) \otimes \End\bigl(\C^N\bigr)$ such that
\begin{equation*}
	\bigl(X^{-1}\bigr)_{ij} (A) = \langle v_i, A \phi_j\rangle = \bigl\langle A^{-1} v_i, \phi_j \bigr\rangle.
\end{equation*}

We consider $\cO(\GL_N)$ as a left $\GL_N \times \GL_N$-module via the corresponding actions
\[
	\bigl(g^R \cdot f\bigr)(A) = f(Ag), \qquad \bigl(g^L \cdot f\bigr)(A) = f\bigl(g^{-1}A\bigr).
\]
In particular, for every $\alpha=1,\dots, N$,
we have equivariant embeddings
\begin{align}
\label{eq::right_embedding}
	\begin{split}
		&{}^R\iota_{\alpha} \colon\ S^{\bullet} \C^N\rightarrow \cO(\GL_N), \qquad v_i \mapsto X_{i\alpha}, \\
		&{}^R \iota_{\alpha}^* \colon\ S^{\bullet} \bigl(\C^N\bigr)^*\rightarrow \cO(\GL_N), \qquad \phi_i \mapsto \bigl(X^{-1}\bigr)_{\alpha i}
	\end{split}
\end{align}
and
\begin{align}
\label{eq::left_embedding}
	\begin{split}
	&{}^L\iota_{\alpha} \colon\ S^{\bullet} \C^N\rightarrow \cO(\GL_N), \qquad v_i \mapsto \bigl(X^{-1}\bigr)_{i \alpha}, \\
	&{}^L \iota_{\alpha}^* \colon\ S^{\bullet} \bigl(\C^N\bigr)^*\rightarrow \cO(\GL_N), \qquad \phi_i \mapsto X_{\alpha i},
	\end{split}
\end{align}
with respect to the action corresponding to the superscript.

Denote by $H_N\subset \GL_N$ the subgroup of diagonal matrices. Its set of weights $\Lambda$ is given by
\begin{equation*}
	\mathbb{Z}^N = \mathrm{span}_{\mathbb{Z}} (\epsilon_1,\dots,\epsilon_N) \xrightarrow{\sim} \Lambda, \qquad \epsilon_i (\mathrm{diag}(t_1,\dots,t_N)) = t_i.
\end{equation*}

We will always identify the weight and the coweight lattices, we hope it will not lead to any confusion. Denote by $\Lambda^+ \subset \Lambda$ the monoid of dominant weights:
\[
	\Lambda^+ = \bigl\{ (\mu_1,\dots,\mu_N) \in \mathbb{Z}^N \mid \mu_{i} \geq \mu_{i+1},\, i=1,\dots N-1, \bigr\},
\]
which is spanned by the fundamental weights
\[
	\omega_i := (\underbrace{1,\dots,1}_i , 0,\dots,0).
\]
Let
\begin{equation*}
\label{eq::roots}
	\Delta = \{\epsilon_i - \epsilon_j \mid i\neq j\} \subset \Lambda, \qquad \Delta_+ = \{\epsilon_i - \epsilon_j \mid i< j\} \subset \Lambda^+
\end{equation*}
be the set of roots (respectively positive roots). We denote by
\begin{equation}
\label{eq::rho}
	\rho_N = \frac{1}{2} \sum_{\alpha\in\Delta_+} \alpha = \sum_i \omega_i = (N-1,\dots,1,0)
\end{equation}
the half-sum of positive roots.

Let $N_+$ (resp.\ $N_-$) be the subgroups of upper unitriangular (resp.\ lower unitriangular) matrices. For $\lambda\in\Lambda$, denote by $V_{\lambda}$ the irreducible representation of $\GL_N$ whose $N_+$-invariant highest weight vector $v_{\lambda}$ has weight $\lambda$. In particular,
\begin{alignat*}{3}
	&V_{(1,0,\dots,0)}= \C^N, &&\qquad v_{(1,0,\dots,0)}= v_1,& \\
	&V_{(0,\dots,0,-1)}= \bigl(\C^N\bigr)^*, &&\qquad v_{(0,\dots,0,-1)}= \phi_N,&
\end{alignat*}

We identify its Lie algebra $\gl_N$ with $\End\bigl(\C^N\bigr)$ via matrix units such that the action on $\C^N$ and~\smash{$\bigl(\C^N\bigr)^*$} is given by
\[
	\gl_N = \mathrm{span}(E_{ij}\mid i,j=1,\dots, N), \qquad E_{ij} \cdot v_k = \delta_{jk} v_i, \qquad E_{ij} \cdot \phi_k = - \delta_{ik} \phi_j.
\]
We denote by the same letters the corresponding generators of the (asymptotic) universal enveloping algebra $\U(\gl_N)$ in the sense of Definition~\ref{def:auea}; in particular,
\[
	E_{ij} E_{kl} - E_{kl} E_{ij} = \hbar (\delta_{jk} E_{il} - \delta_{li} E_{kj}).
\]
We combine the generators into a matrix $E \in \U(\gl_N) \otimes \End\bigl(\C^N\bigr)$:
\begin{equation}
\label{eq::e_matrix}
	E = \begin{pmatrix}
		E_{11} & \dots & E_{1,N-1} & E_{1N} \\
		\vdots & \ddots & \vdots & \vdots \\
		E_{N-1,1} & \dots & E_{N-1,N-1} & E_{N-1,N} \\
		E_{N1} & \dots & E_{N,N-1} & E_{NN}
	\end{pmatrix},
\end{equation}	

The Lie algebras of $N_+$, $N_-$, $H_N$ are respectively
\[
	\n_+ = \mathrm{span}(E_{ij} \mid i<j), \qquad \n_- = \mathrm{span}(E_{ij} \mid i>j), \qquad \h = \mathrm{span}(E_{ii}).
\]

We will often use an antiautomorphism
\begin{equation}
\label{eq::transposition}
	\varpi_N \colon\ \gl_N \rightarrow \gl^{\mathrm{opp}}_N, \qquad E_{ij} \mapsto E_{ji},
\end{equation}
as well as its version, the automorphism
\begin{equation}
\label{eq::transposition_minus}
	\tilde{\varpi}_N \colon\ \gl_N \rightarrow \gl_N, \qquad E_{ij} \mapsto - E_{ji}.
\end{equation}
Obviously, in terms of the matrix $E$, we have $\varpi(E) = E^{\mathsf T}$.

\subsection{Quantum minors}
\label{subsect:quantum_minors}

Denote by
$T(u) = u + E \in \U(\gl_N)[u] \otimes \End\bigl(\C^N\bigr)$
the \emph{T-operator}.

Let $(a_1,\dots, a_m)$, $(b_1, \dots, b_m)$ be two sequences of indices. Define a {\it quantum minor} as the sum over all permutations of $1,2,\dots,m$, see \cite[formulas~(1.54) and (1.55)]{MolevYangians}:
\begin{align}
\nonumber
T_{b_1\dots b_m}^{a_1\dots a_m}(u) := & \sum\limits_{\sigma} \mathrm{sgn}(\sigma) T_{a_{\sigma(1)}b_1} (u) T_{a_{\sigma(2)},b_{2}} (u-\hbar) \cdots T_{a_{\sigma(m)} b_m} (u-m\hbar+\hbar) \\
		=& \sum\limits_{\sigma} \mathrm{sgn}(\sigma) T_{a_1 b_{\sigma(1)}}(u-m\hbar+\hbar) \cdots T_{a_{m-1} b_{\sigma(m-1)}} (u-\hbar)T_{a_m b_{\sigma(m)}}(u).\label{eq::quantum_minor}
\end{align}

It follows from loc.\ cit.\ that
\begin{equation*}
	\label{eq::t_l_minor_relation}
	u^{-1} \cdots (u-m\hbar + \hbar)^{-1} T_{b_1\dots b_m}^{a_1\dots a_m}(u) \in \U(\gl_N)\bigl[\hspace{-0.2mm}\bigl[u^{-1}\bigr]\hspace{-0.2mm}\bigr] .
\end{equation*}

We will need the following properties of quantum minors.
\begin{Proposition}\quad
	\label{prop::minor_basic}
	\begin{enumerate}\itemsep=0pt
		\item[$(1)$]%\label{prop::minor_basic_perm}
 For any permutation $\sigma$, we have
		\[
		T_{b_1\dots b_m}^{a_{\sigma(1)}\dots a_{\sigma(m)}}(u)
 = \mathrm{sgn}(\sigma) T_{b_1\dots b_m}^{a_1 \dots a_m} = T_{b_{\sigma(1)} \dots b_{\sigma(m)}}^{a_1\dots a_m}(u).
		\]
		In particular, if $a_i = a_j$ or $b_i=b_j$ for some $i\neq j$, then \smash{$T_{b_1\dots b_m}^{a_1 \dots a_m} (u) = 0$}.
		\item[$(2)$]%\label{prop::minor_basic_sum}
 A quantum minor can be decomposed as {\rm \cite[{\it Proposition} 1.6.8]{MolevYangians}}
		\begin{align*}
			T_{b_1 \dots b_m}^{a_1 \dots a_m}(u)&=\sum\limits_{l=1}^m (-1)^{m-l} T_{b_1 \dots b_{m-1}}^{a_1 \dots \hat{a}_l \dots a_m}(u) T_{a_l b_m}(u-\hbar m + \hbar) \\
			&=\sum\limits_{l=1}^m (-1)^{m-l} T_{b_1\dots \hat{b}_l \dots b_m}^{a_1\dots a_{m-1}}(u-\hbar) T_{a_m b_l}(u) \\
			&=\sum\limits_{l=1}^m (-1)^{l-1} T_{a_l b_1} (u) T_{b_2\dots b_m}^{a_1 \dots \hat{a}_l \dots a_m} (u-\hbar) \\
			&=\sum\limits_{l=1}^m (-1)^{l-1} T_{a_1 b_l}(u-\hbar m + \hbar) T_{b_1 \dots \hat{b}_l \dots b_m}^{a_2 \dots a_m}(u).
		\end{align*}
		\item[$(3)$] The following relation holds {\rm \cite[{\it Proposition} 1.7.1]{MolevYangians}}:
		\[
			\bigl[E_{kl}, T^{a_1\dots a_m}_{b_1\dots b_m}(u)\bigr] = \hbar\sum_{i=1}^m \delta_{a_i l} T^{a_1\dots k \dots a_m}_{b_1\dots b_m}(u) - \hbar\sum_{i=1}^m \delta_{b_i k} T^{a_1\dots a_m}_{b_1\dots l \dots b_m}(u),
		\]
		where the indices $k$ and $l$ in the quantum minors replace correspondingly $a_i$ and $b_i$.
	\end{enumerate}
	\begin{proof}
		For the last part, it is enough to consider the coefficient of $u^0$ on both sides of \cite[Proposition 1.7.1]{MolevYangians}.
	\end{proof}
\end{Proposition}

Denote by
\begin{equation}
	\label{eq::quantum_determinant}
		A(u) := (-1)^N T^{1\dots N}_{1\dots N} (-u+N\hbar-\hbar) =: \sum\limits_{i=0}^N A_i u^{N-i}.
\end{equation}
By \cite{MolevNazarovOlshanskii}, the coefficients $\{A_i\}$ generate the center of $\U(\gl_N)$.

We will need an explicit form of the inverse $T(u)^{-1}$ in terms of \emph{quantum comatrix}, see \cite[Proposition 1.9.2]{MolevYangians}:
\begin{equation}
\label{eq::quantum_comatrix}
	\hat{T}_{ij}(u) := (-1)^{i+j} T_{1\dots \hat{i} \dots N}^{1\dots \hat{j} \dots N}(u),
\end{equation}
where the hats indicate the omitted indices. The following equation is an immediate corollary of \cite[Definition 1.9.1 and Proposition 1.9.2]{MolevYangians}:
\begin{equation}
\label{eq::inverse_vs_comatrix}
	T(u)^{-1} = (-1)^N \hat{T}(u+N\hbar-\hbar) A(-u)^{-1}.
\end{equation}

Recall the antiautomorphism \eqref{eq::transposition_minus}. The next lemma follows immediately from the definition~\eqref{eq::quantum_minor} and $\tilde{\varpi}_N[T(u)] = -T(-u)^{\mathsf T}$, where the superscript means matrix transposition.
\begin{Lemma}
\label{lm:quantum_determinant_transposition}
	Let $\tilde{\varpi}_N$ be the automorphism of $\U(\gl_N)$ given by $E_{ij} \mapsto - E_{ji}$. Then
	\[
		\tilde{\varpi}_N (A(u)) = (-1)^N A(-u + N\hbar - \hbar).
	\]
\end{Lemma}

\subsection{Differential operators}
\label{subsect:diff_ops}

Let $\dD(\GL_N)$ be the Rees algebra of the differential operators on $\GL_N$ with respect to the degree filtration. For $\xi \in \gl_N$, consider the corresponding asymptotic vectors fields generating the right (resp.\ the left) translations along $\xi$:
\begin{equation*}
	\xi^R(f) (g):= \hbar \frac{\mathrm{d}}{\mathrm{d} t} f\bigl(g{\rm e}^{t\xi}\bigr) \Big|_{t=0}, \qquad \xi^L(f) (g):= \hbar \frac{\mathrm{d}}{\mathrm{d} t} f\bigl({\rm e}^{-t\xi} g\bigr) \Big|_{t=0}.
\end{equation*}
In terms of either trivialization, we have $\dD(\GL_N) \cong \U(\gl_N) \otimes \cO(\GL_N)$.

Denote by $p_{ij} := \hbar \partial/\partial x_{ij}$ and by $D = (p_{ij})$ the corresponding matrix. We have
\begin{equation*}
	E_{ij}^R = \sum_{\alpha=1}^{N} x_{\alpha i} p_{\alpha j}, \qquad E_{ij}^L = -\sum_{\alpha=1}^N x_{j \alpha} p_{i\alpha}.
\end{equation*}
In particular, in terms of matrices from Section~\ref{subsect:notations}, we can write
\[
	E^R = \bigl(E_{ij}^R\bigr) = XD, \qquad E^L = \bigl(E_{ij}^L\bigr) = - DX + \hbar N \cdot \mathrm{Id}_N,
\]
so that
\begin{equation}
\label{eq::right_vs_left_matrix}
	E^L = X^{-1}\bigl(-E^R + \hbar N \cdot \mathrm{Id}_N\bigr) X.
\end{equation}

There is an embedding $\Z_{\gl_N} \hookrightarrow \dD(\GL_N)$ as bi-invariant differential operators. Denote by~$A^R(u)$ \big(resp.\ $A^L(u)$\big) the quantum determinant \eqref{eq::quantum_determinant} constructed out of the matrix $E^R$ \big(resp.~$E^L$\big).
\begin{Proposition}
\label{prop:center_left_vs_right}
	We have $A^L(u) = (-1)^N A^R(-u+N\hbar-\hbar)$.
	\begin{proof}
		It is enough to show that their actions on $\cO(\Mat_{N\times N})$ coincide. By the Howe duality \cite[Theorem~5.16]{ChengWang}, under the $(\gl_N,\gl_N)$-action by
		\[
			E^{(1)}_{ij} \mapsto \sum_{\alpha=1}^N x_{\alpha i} p_{\alpha j}, \qquad E^{(2)}_{ij} \mapsto \sum_{\alpha=1}^N x_{i \alpha} p_{j \alpha},
		\]
		there is an isomorphism
		\[
			\cO(\Mat_{N\times N}) \cong \bigoplus_{\lambda} V_{\lambda} \otimes V_{\lambda}
		\]
		of $(\gl_N,\gl_N)$-modules, where $\lambda$ runs over certain dominant weights. In particular, we have \smash{$A^{(1)}(u) = A^{(2)}(u)$}, where the quantum determinants are constructed out of corresponding matrices $E^{(i)}$ for $i=1,2$. Denote by
		\[
			\tilde{\varpi}_N \colon \ \U(\gl_N) \rightarrow \U(\gl_N), \qquad \tilde{\varpi}_M(E_{ij}) = -E_{ji}
		\]
		an algebra automorphism of $\U(\gl_N)$, a variant of \eqref{eq::transposition}. Then \smash{$E^{(1)}_{ij} = E^R_{ij}$} and \smash{$E^{(2)}_{ij} = \tilde{\varpi}_N \bigl(E^{L}_{ij}\bigr)$}, therefore,
		\[
			A^L(u) = \tilde{\varpi}_N\bigl(A^R(u)\bigr).
		\]
		The result follows from Lemma~\ref{lm:quantum_determinant_transposition}.
	\end{proof}
\end{Proposition}

Denote by $\hat{T}^R(u)$ \big(resp.\ $\hat{T}^L(u)$\big) the quantum comatrix \eqref{eq::quantum_comatrix} associated to $E^R$ \big(resp.\ $E^L$\big).

\begin{Proposition}
\label{prop:comatrix_left_right}
	We have $X^{-1}\cdot \hat{T}^R(-u+N\hbar-\hbar)\cdot X = (-1)^{N+1} \hat{T}^L(u)$.
	\begin{proof}
		It follows from Propositions~\ref{prop:canonical_matrix_desciption} and~\ref{prop:conjugation_center} that
		\[
			A^R(u) X = \frac{u-E^R}{u-E^R+\hbar} X A^R(u),
		\]
		where the equality is understood as in Remark~\ref{remark:pwoer_series} and multiplication by $A^R(u)$ is entry-wise. Therefore, by \eqref{eq::inverse_vs_comatrix} and $T^R(u) = u+E^R$, we have
		\begin{align*}
			X^{-1}\cdot \hat{T}^R(-u+N\hbar-\hbar)\cdot X &=(-1)^N X^{-1} \frac{1}{-u+E^R} A^R(u) X \\
			&= (-1)^{N+1} X^{-1} \frac{1}{u-E^R+\hbar} X A^R(u).
		\end{align*}
		It follows from \eqref{eq::right_vs_left_matrix} and Proposition~\ref{prop:center_left_vs_right} that
		\begin{align*}
			(-1)^{N+1} X^{-1} \frac{1}{u-E^R+\hbar} X A^R(u) &= - \frac{1}{u+E^L - \hbar N + \hbar} A^L(-u+N\hbar-\hbar) \\
			&=-T^L(u-N\hbar+\hbar)^{-1} A^L(-(u-N\hbar+\hbar)).
		\end{align*}
		Applying \eqref{eq::inverse_vs_comatrix} for $E^L$, we get
		\[
			X^{-1}\cdot \hat{T}^R(-u+N\hbar-\hbar)\cdot X = (-1)^{N+1} \hat{T}^L(u),
		\]
		as required.
	\end{proof}
\end{Proposition}

\section{Degenerate affine Hecke algebras and Yangians}
\label{sect:deg_aha_yangians}

The goal of this section is to establish a relation between rational Feigin--Odesskii shuffle algebras and degenerate affine Hecke algebras summarized in Theorem~\ref{thm:hecke_vs_fo_shuffle}. Using, on the one hand, the action of latter on certain Harish-Chandra bimodules and their reductions as in Section~\ref{sect:kw_reduction}, on~the~other hand, connection between the former and (truncated shifted) Yangians of $\sl_2$ which we recall in Section~\ref{subsect:shifted_yangian}, we apply the results of this section to prove the main result of the paper, Theorem~\ref{thm:toda_lattice_yangian}. This relation is deduced through Shibukawa--Ueno $R$-operators whose definition and properties we recall in Section~\ref{subsect:shibukawa_ueno}.

\subsection[Shibukawa-Ueno R-operator]{Shibukawa--Ueno $\boldsymbol{R}$-operator}
\label{subsect:shibukawa_ueno}
In this subsection, we recall the definition of Shibukawa--Ueno $R$-operators \cite{ShibukawaUeno} and study the induced symmetrizers. By \cite{FelderPasquier}, its elliptic version is gauge-equivalent to Belavin's $R$-matrix that is used to define the original elliptic shuffle algebras \cite{FeiginOdesskii, FeiginOdesskiiSklyanin,OdesskiiFeiginEllipticDeformation}. At the same time, we will show in the next subsection that its rational version is related to degenerate affine Hecke algebras.

Let $\theta(z)$ be a skew-symmetric (i.e., $\theta(-z)=-\theta(z)$) holomorphic function satisfying the three-term relation
\begin{align}
		&\theta(x+y) \theta(x-y) \theta(z+w) \theta(z-w) + \theta(x+z) \theta(x-z) \theta(w+y) \theta(w-y) \nonumber \\
		&\qquad{} +\theta(x+w) \theta(x-w) \theta(y+z) \theta(y-z) = 0.\label{eq::three_term_relation}
\end{align}
Its analytic solution is essentially one of the following three functions \cite{WhittakerWatson}:
\begin{equation}
\label{eq::generalized_theta_functions}
	\theta^{\mathrm{ell}}(z) := \vartheta_1(z,\tau), \qquad \theta^{\mathrm{trig}}(z) := \sin(\pi z/\tau_1), \qquad \theta^{\rat} (z):= z,
\end{equation}
where $\tau \in \mathbb{H}$ is in the upper half-plane, $\tau_1 \in \C\backslash \{0\}$, and
\[
	\vartheta_1(z,\tau) = - \sum_{m\in \mathbb{Z}} \exp \Biggl[ \pi {\rm i} \tau \biggl(m+\frac{1}{2} \biggr)^2 + 2\pi {\rm i}\biggl( m+\frac{1}{2} \biggr) \biggl( z+ \frac{1}{2} \biggr) \Biggr]
\]
is the Jacobi's first theta function. Denote by
\begin{equation*}
	G^{\theta}(z,w) := \frac{\theta(z+w) \theta'(0)}{\theta(z)\theta(w)}.
\end{equation*}

Let $\cM_k$ the space of meromorphic functions in variables $z_1,\dots,z_k$. Denote by
\begin{equation}
\label{eq::permutation}
	(P_{ij} f)(z_1,\dots,z_i,\dots,z_j,\dots,z_k) := f(z_1,\dots,z_j,\dots,z_i,\dots,z_k)
\end{equation}
the permutation operator in $\cM_k$. Let $\kappa \in \C$ be a number such that $\theta(k\kappa) \neq 0$ for any $k\in \mathbb{Z}$.
\begin{Definition}[{\cite{ShibukawaUeno}}]
\label{def:shibukawa_ueno}
	The {\it Shibukawa--Ueno $R$-operator} associated to $\theta$ is the operator on $\cM_k$ defined by
	\[
		R_{ij}^{\theta}(\xi|\kappa):= G^{\theta}(z_i-z_j,\kappa) \mathrm{Id}_{\cM_k} - G^{\theta}(z_i-z_j,\xi) P_{ij} .
	\]
	The Shibukawa--Ueno {\it braiding element} is
	\[
		\check{R}^{\theta}_i (\xi|\kappa) := P_{i,i+1} \circ R^{\theta}_{i,i+1}(\xi |\kappa) = G^{\theta} (z_{i+1}-z_i,\kappa) P_{i,i+1} - G^{\theta}(z_{i+1}-z_i,\xi) \mathrm{Id}_{\cM_k}.
	\]
	The {\it elliptic} (resp.\ {\it trigonometric, rational}) Shibukawa--Ueno $R$-operator is
	\[
		R^{\mathrm{ell}}_{ij}(\xi|\kappa) := R_{ij}^{\theta^{\mathrm{ell}}}(\xi|\kappa), \qquad R^{\mathrm{trig}}_{ij}(\xi|\kappa) := R_{ij}^{\theta^{\mathrm{trig}}}(\xi|\kappa), \qquad R^{\mathrm{rat}}_{ij}(\xi|\kappa) := R_{ij}^{\theta^{\mathrm{rat}}}(\xi|\kappa)
	\]
	for the functions of \eqref{eq::generalized_theta_functions}. If dependence on $\kappa$ is clear from the context, we denote $R^{\theta}_{ij}(\xi) := R^{\theta}_{ij}(\xi|\kappa)$, and similarly for other operators.
\end{Definition}

\begin{Proposition}[{\cite{ShibukawaUeno}}]
	The Shibukawa--Ueno $R$-operator satisfies the quantum Yang--Baxter equation
	\[
		R^{\theta}_{12}(\xi_1-\xi_2|\kappa) R^{\theta}_{13}(\xi_1-\xi_3|\kappa) R^{\theta}_{23}(\xi_2-\xi_3|\kappa) = R^{\theta}_{23}(\xi_2-\xi_3|\kappa) R^{\theta}_{13}(\xi_1-\xi_3|\kappa)R^{\theta}_{12}(\xi_1-\xi_2|\kappa)
	\]
	for any $\theta(z)$ satisfying \eqref{eq::three_term_relation}.
\end{Proposition}

Consider the $\kappa$-symmetrization operator \cite{CherednikRmatrix}
\begin{equation*}
	\check{H}_k^{\theta} (\kappa):= \bigl(\check{R}^{\theta}_{k-1}(-\kappa|\kappa) \cdots \check{R}^{\theta}_{1}(-k\kappa+\kappa|\kappa)\bigr)\cdots \bigl(\check{R}^{\theta}_{k-1}(-\kappa|\kappa)\check{R}^{\theta}_{k-2} (-2\kappa|\kappa)\bigr)\bigl(\check{R}^{\theta}_{k-1}(-\kappa|\kappa)\bigr).
\end{equation*}
Denote by
\begin{equation*}
	(\Sym_k f)(z_1,\dots,z_k) := \sum_{\sigma\in S_k} f(z_{\sigma(1)},\dots,z_{\sigma(k)})
\end{equation*}
the standard symmetrizer.

\begin{Proposition}
\label{prop:su_symmetrizer}
	We have
	\[
		\bigl(\check{H}_k^{\theta} (\kappa) f\bigr)(z_1,\dots,z_k) = \Sym_k \biggl(\prod_{i<j} G^{\theta}(z_i - z_j,\kappa) f(z_1,\dots,z_k) \biggr).
	\]
	\begin{proof}
		For brevity, we fix $\kappa$ and will not indicate dependence on it explicitly. Observe that \smash{$\check{H}^{\theta}_k := \prod_{i<j} P_{ij} \circ H^{\theta}_k$}, where the product is taken in lexicographic order on pairs $(i,j)$ and
		\[
			H_k^{\theta} = \bigl(R_{12}^{\theta}(-\kappa) \cdots R_{1k}(-k\kappa+\kappa)\bigr) \cdots \bigl(R^{\theta}_{k-2,k-1}(-\kappa) R^{\theta}_{k-2,k}(-2\kappa)\bigr) \bigl(R^{\theta}_{k-1,k}(-\kappa)\bigr).
		\]
		Therefore, it would be enough to show that
		\[
			\bigl(H_k^{\theta} f\bigr)(z_1,\dots,z_k) = \Sym_k \biggl(\prod_{i<j} G^{\theta}(z_i - z_j,\kappa) f(z_1,\dots,z_k) \biggr).
		\]
		We can prove by induction on $k$. The base $k=1$ is tautological. For the step, let
		\begin{align*}
			F(z_0,z_1,z_2,\dots,z_{k}) :={}& \bigl(H_k^{\theta} f\bigr) (z_0,z_1,\dots,z_k) \\
 ={}& \Sym_k \biggl(\prod_{1\leq i < j \leq k} G^{\theta}(z_i-z_j) f(z_0,z_1,\dots,z_k) \biggr),
		\end{align*}
		where symmetrization is over the variables $z_1,\dots,z_k$; in particular, $F$ is symmetric with respect to them.
		Then it would be enough to prove that
		\[
			R_{01}(-\kappa) \cdots R_{0k}(-k\kappa) F = (1+ P_{01} + \dots + P_{0k}) \prod_{i=1}^{k} G^{\theta}(z_0-z_i,\kappa) F(z_0,\dots,z_k).
		\]
		Denote by \smash{$\cM_{k+1}^{S_l}$} the space of meromorphic functions in $z_0,z_1,\dots,z_k$ symmetric with respect to the variables $z_1,\dots,z_l$. For any $1\leq l \leq k$, define
		\[
			X_l \colon\ \cM_{k+1}^{S_l} \rightarrow \cM_{k+1}^{S_l}, \qquad X_l = R^{\theta}_{01}(-\kappa) \cdots R_{0l}(-l\kappa).
		\]
		Let us show by induction on $l$ that, as operators,
		\[
			X_l = (1+P_{01} + \dots + P_{0l}) \prod_{i=1}^{l} G^{\theta}(z_0-z_i).
		\]
		For $l=1$, it follows from $\theta(-z)=-\theta(z)$ for \eqref{eq::generalized_theta_functions}. For the step, we have
		\begin{align*}
			X_{l+1}={}& X_l R_{0,l+1}(-l\kappa-\kappa) \\
 ={}&(1+P_{01}+\dots+P_{0l})\\
 &{\times}\,\prod_{i=1}^{l} G^{\theta}(z_0-z_i,\kappa) \cdot \bigl(G^{\theta}(z_0-z_{l+1},\kappa) + G^{\theta}(z_{l+1}-z_0,l\kappa+\kappa) P_{0,l+1}\bigr).
		\end{align*}
		Therefore, it is enough to show that
		\begin{gather*}
			(1+P_{01} + \dots + P_{0l}) \prod_{i=1}^{l} G^{\theta}(z_0-z_i,\kappa) \cdot G^{\theta}(z_{l+1} - z_0,l\kappa+\kappa) P_{0,l+1}\\
 \qquad = P_{0,l+1} \prod_{i=1}^{l+1} G^{\theta}(z_0-z_i,\kappa).
		\end{gather*}
		The right-hand side is \smash{\raisebox{0.2pt}{$\prod_{i=0}^{l} G^{\theta}(z_{l+1}-z_i,\kappa) P_{0,l+1}$}}. Observe that we have $P_{0i} P_{0,l+1} = P_{0,l+1}$ for any $1\leq i \leq l$ as operators on \smash{$\cM\raisebox{-0.5pt}{${}_{k+1}^{S_{l+1}}$}$}.
Therefore, the left-hand side is
		\[
			\sum_{i=0}^l \prod_{\substack{j=0,\dots, l, \\ j \neq i}} G^{\theta}(z_i - z_j,\kappa)\cdot G^{\theta}(z_{l+1} - z_i,l\kappa+\kappa) P_{0,l+1}.
		\]
		Hence, the equality would follow from
		\[
			\prod_{i=0}^{l} G^{\theta}(z_{l+1}-z_i,\kappa) = \sum_{i=0}^l \prod_{\substack{j=0,\dots, l, \\ j \neq i}} G^{\theta}(z_i - z_j,\kappa)\cdot G^{\theta}(z_{l+1}- z_i,l\kappa+\kappa).
		\]
		This is the content of Lemma~\ref{lm:generalized_three_term} for $u=z_{l+1}$.
	\end{proof}
\end{Proposition}

Consider the $\kappa$-antisymmetrizer
\[
	\check{A}^{\theta}_k(\kappa) :=\bigl(\check{R}^{\theta}_{k-1}(\kappa|\kappa)\bigr) \bigl(\check{R}^{\theta}_{k-2}(2\kappa|\kappa) \check{R}^{\theta}_{k-1}(\kappa|\kappa)\bigr) \cdots \bigl(\check{R}^{\theta}_{1}(k\kappa-\kappa|\kappa) \cdots \check{R}^{\theta}_{k-1}(\kappa|\kappa)\bigr).
\]

Denote by
\[
	\bigl(\aSym_k f\bigr)(z_1,\dots,z_k) := \sum_{\sigma\in S_k} (-1)^{\mathrm{sgn}(\sigma)} f(z_{\sigma(1)},\dots,z_{\sigma(k)})
\]
the standard antisymmetrization map.
\begin{Corollary}
\label{cor:su_anti_symm}
	We have
	\[
		\bigl(\check{A}^{\theta}_k f\bigr) (z_1,\dots,z_k) = \prod_{i<j} G^{\theta}(z_i-z_j,-\kappa) \cdot \aSym_k f(z_1,\dots,z_k).
	\]
	\begin{proof}
		Let us consider $\Sym_k$ and $\Sym_k^-$ as elements of $\C[S_k]$ embedded into semi-direct product~$\C[S_k] \ltimes \cM_k$. Then Proposition~\ref{prop:su_symmetrizer} for the quantization parameter $-\kappa$ is equivalent to~equality
		\[
			\check{H}^{\theta}_k (-\kappa) = \Sym_k\cdot \prod_{i<j} G^{\theta}(z_i-z_j,-\kappa)
		\]
		in $\C[S_k] \ltimes \cM_k$. Consider the antiautomorphism
		\[
			\varpi_k \colon\ \C[S_k] \ltimes \cM_k \rightarrow (\C[S_k] \ltimes \cM_k)^{\mathrm{opp}}, \qquad \varpi_k |_{\cM_k} = \mathrm{Id}, \qquad \varpi_k (P_{ij}) = - P_{ij}.
		\]
		Then $\varpi_k \bigl(\check{R}^{\theta}_i (\xi|-\kappa)\bigr) = \check{R}^{\theta}_i (\xi|\kappa)$, in particular, $\varpi_k\bigl(\check{H}^{\theta}_k(-\kappa)\bigr) = \check{A}^{\theta}_k(\kappa)$ by definition. Combining with the equality above and using $\varpi_k(\Sym_k) = \Sym_k^-$, we obtain
		\[
			\check{A}^{\theta}_k(\kappa) = \varpi_k\bigl(\check{H}^{\theta}_k(-\kappa)\bigr) = \prod_{i<j} G^{\theta}(z_i-z_j,-\kappa) \Sym_k^-,
		\]
		as required.
	\end{proof}
\end{Corollary}

\subsection{Degenerate affine Hecke algebras and shuffle products}
\label{subsect:deg_affine_hecke}

In this subsection, we recall some basic facts about degenerate affine Hecke algebras and relate them to Feigin--Odesskii shuffle algebras \cite{FeiginOdesskii}. Unless otherwise stated, we refer the reader to~\cite{KirillovHecke} for all the statements and proofs (upon replacing $h$ in loc.\ cit.\ with $-\kappa$).

\begin{Definition}	\label{def:deg_affine_hecke}
	The {\it degenerate affine Hecke algebra} $\fH^{\kappa}_k$ is a $\C[\hbar]$-algebra generated by~\smash{$\{\sigma_i\}_{i=1}^{k-1}$} with relations
	\[
	(\sigma_i)^2 = 1, \qquad \sigma_i \sigma_j = \sigma_j \sigma_i, \qquad \sigma_{i} \sigma_{i+1} \sigma_i = \sigma_{i+1} \sigma_{i} \sigma_{i+1}
	\]
	for all $|i-j| > 1$, and $\{Y_i\}_{i=1}^k$ satisfying
	\[
	Y_i Y_j = Y_j Y_i, \qquad Y_{i+1} = \sigma_i Y_i \sigma_i - \kappa \sigma_i
	\]
	for all $i$, $j$. Here, $\kappa = \pm \hbar$.
\end{Definition}

It contains the group algebra $\C[S_k]$ of the symmetric group and the polynomial algebra $\C[\hbar][Y_1,\dots,Y_k]$. The latter can be endowed with a $\fH^{\kappa}_k$-module structure via induction from the trivial $\C[S_k]$-representation
\begin{equation}
\label{eq::hecke_polynomial_representation}
	\C[\hbar][Y_1,\dots,Y_k] \cong \fH^{\kappa}_k \otimes_{\C[S_k]} \C
\end{equation}
such that generators $\{Y_i\}$ act by multiplication. In what follows, for any element $\sigma\in S_k$, we denote by $\sigma^{\kappa}$ its action on the module above, in particular,
\begin{equation}
\label{eq::hecke_permutation_action}
	\sigma^{\kappa}_i = P_{i,i+1} + \kappa \frac{1 - P_{i,i+1}}{Y_i - Y_{i+1}},
\end{equation}
where the operator $P_{i,i+1}$ is defined by \eqref{eq::permutation}. It follows that the usual symmetric polynomials $\C[\hbar][Y_1,\dots,Y_k]^{S_k}$ lie in the center of $\fH^{\kappa}_k$ (and, in fact, coincide with it).

By affinization procedure \cite[Section 8.7.1]{KSQuantumGroups}, the collection of operators
\begin{equation*}
	\check{R}^{\Hecke}_i(\xi) := \sigma_i+ \frac{1}{\xi} \in \C[S_k]
\end{equation*}
satisfies the braid relations
\begin{align*}
	&\check{R}^{\Hecke}_i(\xi_{i+1} - \xi_{i+2}) \check{R}^{\Hecke}_{i+1}(\xi_i - \xi_{i+2}) \check{R}^{\Hecke}_i(\xi_{i} - \xi_{i+1}) \\
	&\qquad=\check{R}^{\Hecke}_{i+1}(\xi_{i} - \xi_{i+1}) \check{R}^{\Hecke}_i(\xi_{i} - \xi_{i+2}) \check{R}^{\Hecke}_{i+1}(\xi_{i+1} - \xi_{i+2}).
\end{align*}

Consider the (anti-)symmetrizer
\begin{equation}
\label{eq::symmetrizer}
	e_k := \frac{1}{k!} \sum_{\sigma\in S_k} \sigma, \qquad {\rm e}^-_k := \frac{1}{k!} \sum_{\sigma\in S_k} (-1)^{\sigma} \sigma
\end{equation}
in $\C[S_k]$. The following is an analog of \cite{Jucis}.

\begin{Proposition}
\label{prop:r_matrix_idempotents}
	We have
	\begin{align*}
		&e_k= \frac{1}{k!} \bigl[\check{R}^{\Hecke}_{k-1}(1) \cdots \check{R}^{\Hecke}_1(k-1)\bigr] \cdots \bigl[\check{R}^{\Hecke}_{k-1}(1) \check{R}^{\Hecke}_{k-2}(2)\bigr] \bigl[\check{R}^{\Hecke}_{k-1}(1)\bigr], \\
		&e_k^-= \frac{(-1)^{\frac{k(k+1)}{2}}}{k!} \bigl[\check{R}^{\Hecke}_{k-1}(-1)\bigr] \bigl[\check{R}^{\Hecke}_{k-2}(-2) \check{R}^{\Hecke}_{k-1}(-1)\bigr] \cdots \bigl[\check{R}^{\Hecke}_{1}(-k+1) \cdots \check{R}^{\Hecke}_{k-1}(-1)\bigr].
	\end{align*}
	\begin{proof}
		It is enough to check the relation in some faithful representation. As such, we can take~\smash{$\bigl(\C^N\bigr)^{\otimes k}$} with the usual permutation action for $N\gg 0$. Then it is the content of \cite[Proposition 1.6.2]{MolevYangians} and \cite[Proposition 1.6.3]{MolevYangians}.
	\end{proof}
\end{Proposition}

Consider the action of Hecke $R$-matrices in the polynomial representation \eqref{eq::hecke_polynomial_representation}. Using explicit formula \eqref{eq::hecke_permutation_action} and Definition~\ref{def:shibukawa_ueno}, we see that it is essentially given by the rational Shibukawa--Ueno braiding element
\begin{equation*}
	\check{R}^{\Hecke}_i(\xi) = \kappa \check{R}^{\rat}_i (-\kappa \xi).
\end{equation*}
The following proposition is a rational analog of the relation between Hall--Littlewood polynomials \cite[Section~III]{Macdonald} and affine Hecke algebras (see \cite{LusztigCharacter}), and follows immediately from Proposition~\ref{prop:r_matrix_idempotents} and~\ref{prop:su_symmetrizer}, and Corollary~\ref{cor:su_anti_symm}.
\begin{Proposition}
\label{prop:symmetrizer_hecke_vs_shuffle}
	The action of $($anti-$)$symmetrizer in the representation \eqref{eq::hecke_polynomial_representation} is given by
	\begin{align*}
		&({\rm e}^{\kappa}_k f)(Y_1,\dots,Y_k)= \frac{1}{k!} \Sym_k \left(\prod_{i<j} \frac{Y_i-Y_j+\kappa}{Y_i - Y_j} f(Y_1,\dots,Y_k) \right), \\
		&(({\rm e}^-_k)^{\kappa} f)(Y_1,\dots,Y_k)= \frac{1}{k!} \prod_{i<j} \frac{Y_i-Y_j -\kappa}{Y_i-Y_j} \Sym_k^- f(Y_1,\dots,Y_k).
	\end{align*}
\end{Proposition}

The formula above implies that there is an isomorphism
\begin{equation}
\label{eq::hecke_spherical}
	\C[\hbar] [Y_1,\dots,Y_k]^{S_k} \xrightarrow{\sim} {\rm e}^{\kappa}_k \fH^{\kappa}_k {\rm e}^{\kappa}_k, \qquad f(Y_1,\dots,Y_k) \mapsto {\rm e}^{\kappa}_k f(Y_1,\dots,Y_k) {\rm e}^{\kappa}_k.
\end{equation}

\begin{Definition}
\label{def:shuffle_algebra}
	The algebra $\cF^{\kappa}$ is the vector space
	\begin{equation*}
		\label{eq::shuffle_algebra}
		\cF^{\kappa} = \bigoplus_{k=0}^{\infty} \cF^{\kappa}_k, \qquad \cF^{\kappa}_k= {\rm e}^{\kappa}_k \fH^{\kappa}_k {\rm e}^{\kappa}_k = \C[\hbar] [Y_1,\dots,Y_k]^{S_k},
	\end{equation*}
	with the $\C[\hbar]$-linear product $*^{\kappa} \colon \cF^{\kappa}_k \otimes_{\C[\hbar]} \cF^{\kappa}_l \rightarrow \cF^{\kappa}_{k+l}$ defined by
	\begin{equation*}
		\label{eq::shuffle_product}
		(f*^{\kappa}g) (Y_1,\dots,Y_{k+l}) = {\rm e}^{\kappa}_{k+l} f(Y_1,\dots,Y_k) g(Y_{k+1},\dots,Y_{k+l}) {\rm e}^{\kappa}_{k+l}.
	\end{equation*}
	The truncated algebra ${}^N \!\cF^{\kappa}$ is the vector space
	\begin{equation*}
		{}^N \!\cF^{\kappa} = \bigoplus_{k=0}^{\infty} {}^N \!\cF^{\kappa}_k, \qquad {}^N \!\cF^{\kappa}_k = \C[\hbar] [Y_1,\dots,Y_k]^{S_k}_{<N}
	\end{equation*}
	of symmetric polynomials of degree less than $N$ in each variable with the same product.
\end{Definition}

By \eqref{eq::hecke_permutation_action}, the product does restrict to ${}^N \!\cF^{\kappa}$. Notice that it has the same Poincar\'{e} polynomial as $S^{\bullet}\C^N$. Also, observe that when $\hbar=0$, the multiplication map becomes the standard shuffle product
\[
	(f * g)(Y_1,\dots,Y_{k+l}) = \frac{1}{(k+l)!} \Sym_{k+l} (f(Y_1,\dots,Y_k) g(Y_{k+1}, \dots ,Y_{k+l})).
\]
There is another deformation of the latter which is a simplest case of the \emph{Feigin--Odesskii shuffle algebras} \cite{FeiginOdesskii} (note that we use a different normalization constant).

\begin{Definition}
\label{def::fo_shuffle_algebra}
	The {\it rational Feigin--Odesskii shuffle algebra} $\cS^{\kappa}$ is the vector space
	\[
		\cS^{\kappa} = \bigoplus \cS^{\kappa}_k, \qquad \cS^{\kappa}_k = \C[Y_1,\dots,Y_k]^{S_k}
	\]
	with the $\C[\hbar]$-linear product \smash{$*^{\kappa} \colon \cS^{\kappa}_k \otimes_{\C[\hbar]} \cS^{\kappa}_l \rightarrow \cS^{\kappa}_{k+l}$} defined by
	\begin{gather*}
		(f *^{\kappa} g)(Y_1,\dots,Y_{k+l})\\
 \qquad{} = \frac{1}{(k+l)!} \Sym_{k+l} \Biggl( \prod_{\substack{i = 1,\dots, k \\ j = k+1,\dots, k+l}} \frac{Y_i - Y_j + \kappa}{Y_i - Y_j} f(Y_1,\dots,Y_k)g(Y_{k+1},\dots,Y_{k+l}) \Biggr).
	\end{gather*}
	The {\it $N$-truncated} rational Feigin--Odesskii shuffle algebra ${}^N \!\cS^{\kappa}$ is the vector space
	\begin{equation*}
		{}^N \!\cS^{\kappa} = \bigoplus_{k=0}^{\infty} {}^N \!\cS^{\kappa}_k, \qquad {}^N \!\cS^{\kappa}_k = \C[\hbar] [Y_1,\dots,Y_k]^{S_k}_{<N}
	\end{equation*}
	of symmetric polynomials of degree less than $N$ in each variable with the same product.
\end{Definition}

\begin{Remark}
	The original construction of elliptic algebras \cite{FeiginOdesskii, FeiginOdesskiiSklyanin,OdesskiiFeiginEllipticDeformation} defines a shuffle product with factor $\vartheta_1(z)$ (up to shifts) on the symmetric algebra $S^{\bullet} H^0 (\Theta_N)$ of global sections of a~degree~$N$ line bundle on an elliptic curve. By \cite{EndelmanHodgesTwistedSU}, the latter admits a degeneration such that~${\vartheta_1(z) \mapsto z}$ and the symmetric algebra becomes ${}^N \!\cS^{\kappa}$.
\end{Remark}

The following theorem is a corollary of Proposition~\ref{prop:symmetrizer_hecke_vs_shuffle}.

\begin{Theorem}
\label{thm:hecke_vs_fo_shuffle}
	The identity map $\mathrm{Id}\colon \cF_k \rightarrow \cS_k$ is an algebra isomorphism.
\end{Theorem}

\begin{Remark}
	One can similarly relate the trigonometric Shibukawa--Ueno operator to the Laurent polynomial representation of the non-degenerate affine Hecke algebra \cite[Section 12.3]{ChariPressley} via the affinization procedure of \cite{KSQuantumGroups}. It would be interesting to relate the elliptic Shibukawa--Ueno operators to elliptic affine Hecke algebras \cite{GinzburgKapranovVasserot}.
\end{Remark}

\subsection[Shifted Yangian of sl(2)]{Shifted Yangian of $\boldsymbol{\smash{\sl_2}}$}
\label{subsect:shifted_yangian}

In this subsection, we recall the definition and some properties of the shifted Yangian of $\sl_2$. For convenience, we will use its $\gl_2$-version $\rY_{m_1,m_2}(\gl_2)$.

The following definition is \cite[Section 2.1]{FrassekPestunTsymbaliuk}; we use its equivalent form in terms of currents. In particular, we assume that the weight $(m_1,m_2)$ is antidominant, i.e., $m_1 \leq m_2$.

\begin{Definition}
\label{def:shifted_yangian}
	The {\it shifted Yangian} $\rY_{m_1,m_2}(\gl_2)$ is an algebra over $\C[\hbar]$ generated by \linebreak \smash{$\bigl\{x_+^{(i)}, x_-^{(i)}\bigr\}_{i\geq 1}$} and \smash{$\bigl\{d_j^{(i)}\bigr\}_{i > m_j}$} for $j=1,2$ that we combine into the series
	\[
		x_+(u) = \sum_{i\geq 1} x_+^{(i)} u^{-i}, \qquad x_-(u) = \sum_{i\geq 1} x_-^{(i)} u^{-i}, \qquad d_j(u) = u^{-m_j} + \sum_{i > m_j} d_j^{(i)} u^{-i}
	\]
	with values in \smash{$\rY_{m_1,m_2}(\gl_2) \bigl(\bigl(u^{-1}\bigr)\bigr)$} satisfying
	\begin{align}
		\label{eq::yangian_d_d_comm}
		&[d_1(u),d_2(v)]= 0, \\
		\label{eq::yangian_e_f_comm}
		&(u-v)[x_+(u),x_-(v)]= \hbar \bigl(d_1(v)^{-1} d_2(v) - d_1(u)^{-1}d_2(u)\bigr), \\
		\label{eq::yangian_d_e_comm}
		&(u-v) [d_j(u),x_+(v)]= \hbar(\delta_{j1} - \delta_{j2})d_j(u) (x_+(u) - x_+(v)), \\
		\label{eq::yangian_d_f_comm}
		&(u-v) [d_j(u),x_-(v)]= \hbar(\delta_{j2} - \delta_{j1}) (x_-(u)-x_-(v)) d_j(u), \\
		\label{eq::yangian_e_e_comm}
		&(u-v) [x_{\pm}(u),x_{\pm}(v)]= \mp \hbar(x_{\pm}(u)-x_{\pm}(v))^2.
	\end{align}
	The {\it positive $($negative$)$ subalgebra} $\rY(\gl_2)^{\pm}$ is the subalgebra of $\rY_{m_1,m_2}(\gl_2)$ generated by the modes of $x_{\pm}(u)$. The {\it non-negative subalgebra}~\smash{$\rY_{m_1,m_2}(\gl_2)^{\geq 0}$} \big(resp.\ {\it non-positive subalgebra} \smash{$\rY_{m_1,m_2}(\gl_2)^{\leq 0}$}\big) is the subalgebra generated by the modes of $x_+(u)$, $d_j(u)$ (resp.\ $x_-(u)$, $d_j(u)$) for~${j=1,2}$.
\end{Definition}

\begin{Remark}
	By definition, positive and negative subalgebras do not depend on the shift. In~fact, it is also the case for non-negative and non-positive ones: we have
	\[
		\rY_{m_1,m_2} (\gl_2)^{\geq 0} \xrightarrow{\sim} \rY_{0,0} (\gl_2)^{\geq 0}, \qquad x_+(u) \mapsto x_+(u),\qquad d_j(u) \mapsto u^{m_j} d_j(u),
	\]
	and similarly for the non-positive one. Nonetheless, we prefer to keep track of the shifts in what follows.
\end{Remark}

Let $m=m_1-m_2$. The definition of the shifted Yangian $\rY_{m}(\sl_2)$ can be found in \cite[Sec\-tion~B.1]{BravermanFinkelbergNakajima}, see also its form in terms of currents $e(u)$, $f(u)$, $h(u)$ in \cite[Section 6.1]{FKPRW}. Its relation with $\rY_{m_1,m_2}(\gl_2)$ is as follows.
\begin{Proposition}[{\cite[Proposition 2.19]{FrassekPestunTsymbaliuk}}]
	There is a $\C[\hbar]$-algebra embedding
	\[
		\rY_{m}(\sl_2) \rightarrow \rY_{m_1,m_2}(\gl_2),
	\]
	uniquely determined via
	\[
		e(u) \mapsto x_+(u), \qquad f(u) \mapsto x_-(u), \qquad h(u) \mapsto d_1(u)^{-1} d_2(u).
	\]
\end{Proposition}

Let $(m_1,m_2) = (-N,N)$. Let
\[
	 \U(\h)^{\gen} = \C[\hbar][w_1,\dots,w_N]
\biggl[ \frac{1}{w_i - w_j + k\hbar} \: \bigg| \: i,j=1,\dots, N,\, k\in \mathbb{Z} \biggr]
\]
be a localization of the algebra of functions on $\mathbb{A}^N \times \mathbb{A}^1$ (see also \eqref{eq::generic_algebra}). Denote by
\begin{equation}
\label{eq::diff_ops}
	\cA := \U(\h)^{\gen} \bigl[\bu_1^{\pm}, \dots,\bu_N^{\pm}\bigr]
\end{equation}
the algebra of difference operators such that
\[
	 \bu_i^{\pm} f(w) = f(w \pm \hbar \epsilon_i) \bu_i^{\pm},\qquad f(w) \in \U(\h)^{\gen},
\]
where $f(w \pm \hbar \epsilon_i) := f(w_1,\dots,w_i \pm \hbar,\dots,w_N,\hbar)$. The homomorphism of the following theorem is called the \emph{GKLO homomorphism}, see \cite{GKLO} for the non-shifted version, \cite{KWWY} for dominant shifts, and \cite{BravermanFinkelbergNakajima} for general definition. We use the version adapted to $\gl_2$.

\begin{Theorem}[{\cite[Theorem~2.35]{FrassekPestunTsymbaliuk}}]
\label{thm:gklo}
	There is a $\C[\hbar]$-algebra homomorphism
	\[
		\gklo_N \colon\ \rY_{-N,N}(\gl_2) \rightarrow \cA,
	\]
	defined by
	\begin{align*}
			&x_+(u)\mapsto -\sum_{i=1}^N \frac{1}{u-w_i} \prod_{j\neq i} \frac{1}{(w_i - w_j)} \bu_i^{-1}, \\
			&x_-(u)\mapsto \sum_{i=1}^N \frac{1}{u-w_i-\hbar} \prod_{j\neq i} \frac{1}{(w_i - w_j)} \bu_i, \\
			&d_1(u)\mapsto \prod_{i=1}^N (u-w_i), \qquad d_2(u) \mapsto \prod_{i=1}^N (u-w_i - \hbar)^{-1}.
	\end{align*}
\end{Theorem}

This map is not injective; for conjectural description of the kernel, see \cite[Remark B.21]{BravermanFinkelbergNakajima}. The following is a special case of truncated shifted Yangians, see \cite[Appendix B]{BravermanFinkelbergNakajima}.
\begin{Definition}
\label{def:truncated_shifted_yangian}
	The algebra $\rY^0_{-2N}(\sl_2)$ is the image of $\rY_{-N,N}(\gl_2)$ in $\cA$ under $\gklo_N$.
\end{Definition}

Recall the shuffle algebra of Definition~\ref{def:shuffle_algebra}. The following is \cite[Theorem~6.20]{TsymbaliukPBWD} combined with Theorem~\ref{thm:hecke_vs_fo_shuffle}, see also \cite{Enriquez}.

\begin{Theorem}
	\label{thm:yangian_hecke_shuffle}
	There is an algebra isomorphism \smash{$\rY(\gl_2)^{\pm} \xrightarrow{\sim} \cF^{\pm \hbar}$} defined on generators by \smash{$x_{\pm}^{(i)} \mapsto Y_1^{i-1} \in \cF^{\pm \hbar}_1$}.
\end{Theorem}

\section{Harish-Chandra bimodules}
\label{sect:hc_bimodules}

In this section, we recall the definition of the category of Harish-Chandra bimodules and the action of degenerate affine Hecke algebras on its certain objects. The descent of this action to the Kostant--Whittaker reduction of Section~\ref{sect:kw_reduction} will define a homomorphism from the shuffle algebra (and, more generally, truncated shifted Yangian of $\sl_2$) to the quantum Toda lattice.

Let $G$ be an affine algebraic group over $\C$ with the Lie algebra $\g$. Denote by $\Rep(G)$ the category of $G$-representations.

\begin{Definition}
\label{def:auea}
	The {\it asymptotic universal enveloping algebra} $\U(\g)$ is a tensor algebra over $\C[\hbar]$ generated by $\g$ with the relations
	\[
		xy -yx = \hbar [x,y],\qquad x,y\in \g.
	\]
\end{Definition}

Observe that $\U(\g)$ is naturally an object in $\Rep(G)$.

\begin{Definition}
%\label{def::hc_bimodules}
	A {\it Harish-Chandra bimodule} is a left $\U(\g)$-module in $\Rep(G)$. In other words, it is a $G$-representation $X$ equipped with a left action of $\U(\g)$ such that the action map~${\U(\g) \otimes X \rightarrow X}$ is $G$-equivariant. The category of Harish-Chandra bimodules is denoted by~$\HC(G)$.
\end{Definition}

For $X \in \HC(G)$ and $\xi\in\g$, we denote by $\ad_{\xi}\colon X \rightarrow X$ the derivative of the $G$-action on $X$. There is a natural \emph{right} action of $\U(\g)$ defined on generators by
\[
	x \xi := \xi x- \hbar \ad_{\xi} (x),\qquad x\in X,\ \xi \in \g.
\]

For any $V\in \Rep(G)$, we denote by $\free(V) := \U(\g) \otimes V$ the left free Harish-Chandra bimodule whose left $\U(\g)$-action is by multiplication on the left component and the $G$-action is diagonal. In fact, this assignment defines a functor $\free \colon \Rep(G) \rightarrow \HC(G)$. In what follows, we will often drop the tensor product sign, for instance, for each $v\in V$ and $x,y\in \U(\g)$, we denote $xvy:= (x \otimes v)y$. For any $V,W\in \Rep(G)$, we have a natural identification
\begin{equation}
	\label{eq::hom_hc_rep_g}
	\Hom_{\HC(\GL)} (\U(\g) \otimes V, \U(\g) \otimes W ) \cong \Hom_{\Rep(G)} (V, \U(\g) \otimes W),
\end{equation}
since $\U(\g) \otimes V$ is a left free $\U(\g)$-module.

\begin{Remark}
\label{remark:hc_hbar_zero}
	The subcategory of $\HC(G)$ of modules where $\hbar$ acts by zero can be identified with the category of modules over the symmetric algebra $S^{\bullet}(\g)$ in $\Rep(G)$, in other words, with the category \smash{$\QCoh_{[\g^*/G]}$} of $G$-equivariant quasi-coherent sheaves on $\g^*$. In particular, for $V\in \Rep(G)$, the quotient $\U(\g) \otimes V/(\hbar) \cong S^{\bullet}(\g) \otimes V$ corresponds to the trivial vector bundle $\g^* \times V \rightarrow \g^*$ with the diagonal $G$-equivariant structure.
\end{Remark}

The functor $\free$ is monoidal, i.e., for every $V,W\in\Rep(G)$, we have natural isomorphisms
\[
(\U(\g) \otimes V) \otimes_{\U(\g)} (\U(\g) \otimes W) \cong \U(\g) \otimes V \otimes W.
\]
In particular, we have the following construction.
\begin{Notation}
\label{def:hc_algebra}
	Let $A$ be an algebra in $\Rep(G)$. We denote by $\U(\g) \ltimes A$ the semi-direct product algebra where the multiplication is defined by
	\begin{gather*}
		(\U(\g) \otimes A) \otimes (\U(\g) \otimes A) \rightarrow (\U(\g) \otimes A) \otimes_{\U(\g)} (\U(\g) \otimes A)\\
 \qquad{} \cong \U(\g) \otimes A \otimes A \rightarrow \U(\g) \otimes A,
	\end{gather*}
	where the first arrow is the projection map and the last one is the multiplication on $A$. In other words, it satisfies
	\[
		a \xi = \xi a - \hbar\ad_{\xi}(a),\qquad a\in A,\ \xi \in \g.
	\]
\end{Notation}

\begin{Example}\label{ex:embedding_diff_ops}
	Recall the right embeddings \eqref{eq::right_embedding}. By identifying $\U(\gl_N)$ with \emph{left}-invariant vector fields (that is, generated by vector fields generating right translations), we can extend them to
	\begin{align*}
		&{}^R \iota_{\alpha} \colon\ \U(\gl_N) \ltimes S^{\bullet} \C^N \rightarrow \U(\gl_N) \ltimes \cO(\GL_N) \cong \dD(\GL_N), \\
		&{}^R \iota_{\alpha}^* \colon\ \U(\gl_N) \ltimes S^{\bullet} \bigl(\C^N\bigr)^* \rightarrow \U(\gl_N) \ltimes \cO(\GL_N) \cong \dD(\GL_N),
	\end{align*}
	that we denote by the same letter. Here $\dD(\GL_N)$ is the Rees algebra of differential operators on $\GL_N$, see Section~\ref{subsect:diff_ops}. Similarly, one can extend their left versions \eqref{eq::left_embedding} by identifying $\U(\gl_N)$ with \emph{right}-invariant vector fields.
\end{Example}

\begin{Example}
\label{ex:trick}
	In what follows, we will use the following trick to relate the results of \cite{MolevTransvectorAlgebras} to the setting of the paper. Denote by $\alpha=N+1$. Then there is an algebra map
	\[
		\U(\gl_N) \ltimes S^{\bullet} \C^N \rightarrow \U(\gl_{N+1}), \qquad E_{ij} \mapsto E_{ij}, \qquad v_i \mapsto E_{i\alpha}
	\]
	for all $i,j=1,\dots, N$. Likewise, there is an algebra map
	\[
		\U(\gl_N) \otimes S^{\bullet} \bigl(\C^N\bigr)^* \rightarrow \U(\gl_{N+1}), \qquad E_{ij} \mapsto E_{ij}, \qquad \phi_i \mapsto E_{\alpha i}.
	\]
	Moreover, the anti-automorphism \eqref{eq::transposition} for $\gl_{N+1}$ gives
	\[
		\varpi_N \colon\ \U(\gl_N) \ltimes S^{\bullet} \C^N \xrightarrow{\sim} \bigl[\U(\gl_N) \ltimes S^{\bullet} \bigl(\C^N\bigr)^*\bigr]^{\mathrm{opp}}, \qquad E_{ij} \mapsto E_{ji}, \qquad v_i \mapsto \phi_i,
	\]
	extending transposition $\varpi_N$ for $\gl_N$.
\end{Example}

\subsection{Action of the degenerate affine Hecke algebra}
\label{subsect:deg_affine_hecke_action}

For any $V \in \Rep(G)$, there is a natural map
\begin{equation*}
\label{eq::center_action}
	\Z_{\g} \rightarrow \End_{\HC(G)} (\U(\g) \otimes V), \qquad z \mapsto (x\otimes v \mapsto zx \otimes v).
\end{equation*}

Let $G = \GL_N$. One can easily check that the map
\begin{equation}
\label{eq::canonical_homomorphism}
	\Omega \colon \ \U(\gl_N) \otimes \C^N \rightarrow \U(\gl_N) \otimes \C^N, \qquad x v_i \mapsto \sum_j xE_{ij} v_j,
\end{equation}
is an endomorphism of Harish-Chandra bimodules. It will be convenient to present it in a matrix form
\begin{equation}
\label{eq::omega_matrix_form}
	(\Omega(v_1),\dots,\Omega(v_N))^{\mathsf T} = E \cdot (v_1,\dots,v_N)^{\mathsf T},
\end{equation}
where $E$ is the matrix \eqref{eq::e_matrix}.

Recall the quantum determinant $A(u)$ of \eqref{eq::quantum_determinant} which we can consider as a polynomial with coefficients in \smash{$\End_{\HC(\GL_N)}\bigl(\U(\gl_N) \otimes \C^N\bigr)$}. Since $\Omega$ commutes with the action of the center, substitution of the variable $u$ with $\Omega$ is a well-defined operation.

\begin{Proposition}
\label{prop:cayley_hamilton}
	We have $A(\Omega) = 0$.
	\begin{proof}
According to matrix form of \eqref{eq::omega_matrix_form}, we need to show that $A(E)=0$. This is \cite[Theorem~7.2.1]{MolevYangians}.
	\end{proof}
\end{Proposition}

Recall the degenerate affine Hecke algebra $\fH^{\kappa}_k$ of Definition~\ref{def:deg_affine_hecke}. According to \cite[Example~2.1]{Cherednik}, there is an action of $\fH^{\hbar}_k$ on the Harish-Chandra bimodule \smash{$\U(\gl_N) \otimes \bigl(\C^N\bigr)^{\otimes k}$}. Recall that there is an isomorphism
\[
\U(\gl_N) \otimes \bigl(\C^N\bigr)^{\otimes k} \cong \bigl(\U(\gl_N) \otimes \C^N\bigr) \otimes_{\U(\gl_N)} \dots \otimes_{\U(\gl_N)} \bigl(\U(\gl_N) \otimes \C^N\bigr).
\]
Define
\begin{align}
		&\Omega_i := 1 \otimes \dots \otimes \Omega \otimes \dots \otimes 1 \curvearrowright \bigl(\U(\gl_N) \otimes \C^N\bigr) \otimes_{\U(\gl_N)} \dots \otimes_{\U(\gl_N)} \bigl(\U(\gl_N) \otimes \C^N\bigr), \nonumber \\
		&\sigma_i (x \otimes v_{a_1} \otimes \dots \otimes v_{a_{i}} \otimes v_{a_{i+1}} \otimes \dots \otimes v_{a_k}) := x \otimes v_{a_1} \otimes \dots \otimes v_{a_{i+1}} \otimes v_{a_{i}} \otimes \dots \otimes v_{a_k},\label{eq::affine_hecke_action_vector}
\end{align}
where $\Omega$ acts on the $i$-th factor and $\sigma_i$ is just a permutation of $(i,i+1)$-th factor on \smash{$\bigl(\C^N\bigr)^{\otimes k}$} extended to \smash{$\U(\gl_N) \otimes \bigl(\C^N\bigr)^{\otimes k}$}.

Recall the degenerate affine Hecke algebra $\fH^{\hbar}_k$ from Definition~\ref{def:deg_affine_hecke}. The following can be checked directly.
\begin{Proposition}
\label{prop:affine_hecke_action}
	The assignment
	\[
		\sigma_i \mapsto \sigma_i, \qquad Y_i \mapsto \Omega_i
	\]
	defines an action of $\fH^{\hbar}_k$ on the Harish-Chandra bimodule \smash{$\U(\gl_N) \otimes \bigl(\C^N\bigr)^{\otimes k}$}.
\end{Proposition}

For the following corollary, we use the embedding of $\U(\gl_N) \otimes S^k \C^N$ into $\U(\gl_N) \ltimes S^{\bullet} \C^N$ as in Definition~\ref{def:hc_algebra}.
\begin{Corollary}
\label{cor:action_spherical_on_symmetric}
	The action of Proposition~$\ref{prop:affine_hecke_action}$ induces an action of $\C[\hbar][Y_1,\dots,Y_k]^{S_k}$ on the Harish-Chandra bimodule $\U(\gl_N) \otimes S^k \C^N$ such that
	\begin{gather*}
		f\bigl(\vec{Y}\bigr) \cdot (v_{i_1} \cdots v_{i_k}) = \sum_J c_J \Omega_1^{j_1}(v_{i_1}) \cdots \Omega_k^{i_k}(v_{i_k}), \qquad f\bigl(\vec{Y}\bigr) = \sum_{J = (j_1,\dots,j_k)} c_J Y_1^{j_1} \cdots Y_k^{j_k}.
	\end{gather*}
	\begin{proof}
		Observe that the action is well-defined since $f\bigl(\vec{Y}\bigr)$ is symmetric. The spherical subalgebra \smash{${\rm e}^{\hbar}_k \fH^{\hbar}_k {\rm e}^{\hbar}_k$} acts on the invariants \smash{$\bigl[\U(\gl_N) \otimes \bigl(\C^N\bigr)^{\otimes k}\bigr]^{S_k}$} that, by definition, can be identified with $\U(\gl_N) \otimes S^k \C^N$. The explicit formula follows from Proposition~\ref{prop:affine_hecke_action}.
	\end{proof}
\end{Corollary}

\section{Parabolic restriction and difference operators}
\label{sect:parabolic_reduction}

In this section, we recall the notion of the parabolic restriction for Harish-Chandra bimodules essentially following \cite{KalmykovSafronov}, as well as the classical construction of Mickelsson algebras \cite{ZhelobenkoSalgebras}. In~particular, we identify the one associated to the symmetric algebra $S^{\bullet} \C^N$ with a negative half of the algebra of difference operators, see Theorem~\ref{thm:par_res_to_difference_ops}. It will be used in Section~\ref{sect:miura} to define a~quantum Toda lattice counterpart of the GKLO homomorphism of Theorem~\ref{thm:gklo}, one of the key steps for the main result of the paper Theorem~\ref{thm:toda_lattice_yangian}.

\subsection{Harish-Chandra bimodules for torus}
%\label{subsect:hc_torus}
Let $H_N \subset \GL_N$ be the subgroup of diagonal matrices as in Section~\ref{subsect:notations}. In this subsection, we explicitly describe the category $\HC(H_N)$ and its generic version $\HC(H_N)^{\gen}$ that we will use for the parabolic restriction.

Since $H_N$ is commutative, the algebra $\U(\h)$ can be identified with the functions \smash{$\cO\bigl(\h^* \!\times\! \mathbb{A}^1\bigr)$}, where $\mathbb{A}^1$ corresponds to $\hbar$. We will mostly use $\rho_N$-shifted coordinates \eqref{eq::rho} on $\h^* \times \mathbb{A}^1$
\begin{equation*}
	w_i := E_{ii} +N\hbar - i \hbar \in \U(\h) \subset \U(\h)^{\gen}.
\end{equation*}
Somewhat abusing notations, we will denote a typical element of $\U(\h)$ by
\[
	f(w) := f(w_1,\dots,w_N,\hbar) \in \U(\h)^{\gen}.
\]
We hope it will not lead to any confusion.

Since irreducible representations of $H_N$ are classified by characters in $\Lambda$, we have
\begin{equation}
\label{eq::hc_torus_decomposition}
	\HC(H_N) \cong \bigoplus_{\lambda \in \Lambda} \Mod_{\U(\h)} [\lambda],
\end{equation}
where $\Mod_{\U(\h)} [\lambda]$ is the category of $\U(\h)$-modules to which we assign degree $\lambda$. Under this identification, if $X\in \Mod_{\U(\h)} [\lambda]$, then the left action is that of $\U(\h)$ and the right one is
\begin{equation}
\label{eq::hc_torus_right_action}
	x \cdot f(w) := f(w - \hbar \lambda) \cdot x, \qquad x\in X,\ f(w)\in \U(\h).
\end{equation}
The tensor structure is given component-wise by
\begin{equation}
\label{eq::hc_torus_tensor_structure}
	\Mod_{\U(\h)}[\lambda] \times \Mod_{\U(\h)}[\mu] \rightarrow \Mod_{\U(\h)}[\lambda+\mu], \qquad (X,Y) \mapsto X \otimes_{\U(\h)} Y,
\end{equation}
where the right action on $X$ is used.

We will also need to extend the scalars. Consider the localized algebra
\begin{equation}
\label{eq::generic_algebra}
	\U(\h)^{\gen} := \C[\hbar][w_1,\dots,w_N]
\biggl[ \frac{1}{w_i - w_j + k\hbar} \:\bigg|\: i,j=1,\dots, N,\, k\in \mathbb{Z} \biggr],
\end{equation}
and define the category of \emph{generic} Harish-Chandra bimodules $\HC(H_N)^{\gen}$ as
\[
	\HC(H_N)^{\gen} \cong \bigoplus_{\lambda \in \Lambda} \Mod_{\U(\h)^{\gen}} 	[\lambda],
\]
with similar definitions of the right action \eqref{eq::hc_torus_right_action} and tensor structure \eqref{eq::hc_torus_tensor_structure}.

\subsection{Parabolic restriction}
In this subsection, we recall the notion of parabolic restriction as in \cite{KalmykovSafronov} and explicitly compute in some concrete examples.

For any weight element $x\in \U(\gl_N)$ and $f(w) \in \U(\h)$, observe that
\[
f(w) x = x f(w + \hbar \cdot \mathrm{wt}(x)).
\]
In particular, the tensor product $\U(\gl_N)^{\gen} :=\U(\h)^{\gen} \otimes_{\U(\h)} \U(\g) $ has a well-defined algebra structure. Denote by $\HC(\GL_N)^{\gen}$ the corresponding generic version of Harish-Chandra bimodules for $\GL_N$ and by
\[
X^{\gen} := \U(\gl_N)^{\gen} \otimes_{\U(\gl_N)} X \in \HC(\GL_N)^{\gen}
\]
for any $X\in \HC(\GL_N)$.

\begin{Definition}
	The {\it parabolic restriction functor} is a functor $\res_- \colon \HC(\GL_N) \rightarrow \HC(H_N)^{\gen}$ defined by
	\[
	X \mapsto X^{\gen} \sslash {N_-} := (\n_- \backslash X^{\gen})^{N_-} = \U(\h)^{\gen} \otimes_{\U(\h)} (\n_- \backslash X)^{N_-} , \qquad X \in \HC(\GL_N),
	\]
	where the quotient is taken by the \emph{left} action of $\n_-$ and the invariants are with respect to the diagonal action of $N_-$.
\end{Definition}

In fact, one can explicitly compute the parabolic restriction as follows. By \cite{AsherovaSmirnovTolstoy}, there exists an element $P$ of $\h$-weight zero in a certain completion of \smash{$\U(\gl_N)^{\gen}\bigl[\hbar^{-1}\bigr]$} whose action is well-defined on any right $\U(\gl_N)^{\gen}$-module $F$ with a locally-finite action of $\n_-$ satisfying $F \n_- \subset \hbar F$; this element is called the \emph{extremal projector}. It satisfies the following properties:
\begin{align}
	&P\in 1+ \n_- \U(\gl_N)^{\gen} \cap \U(\gl_N)^{\gen} \n_+, \nonumber \\
	&E_{ij} P= P E_{ji} = 0,\qquad 1 \leq i < j \leq N, \nonumber \\
	&P^2 = P.\label{eq::extremal_proj}
\end{align}

\begin{Notation}
\label{notation:par_res_quotient_extremal_proj}
	For any Harish-Chandra bimodule $X\in \HC(\GL_N)$, we denote by
	\[
		X^{\gen} \mapsto \n_- \backslash X^{\gen}, \qquad x \mapsto Px
	\]
	 the image of $x\in X^{\gen}$ under the projection map. This is a well-defined notation since for any~${\xi \in \n_-}$, we have $P\xi x= 0$.
\end{Notation}
The extremal projector defines an isomorphism
\begin{equation}
	\label{eq::extremal_proj_quotient_inv_iso}
	P \colon \n_- \backslash X^{\gen} /\n_+ \xrightarrow{\sim} (\n_-\backslash X^{\gen})^{N_-}, \qquad [x] \mapsto PxP.
\end{equation}

In particular, for any $V \in \Rep(\GL_N)$, there is a canonical isomorphism
\begin{equation}
	\label{eq::res_trivialization}
	\U(\h)^{\gen} \otimes V \xrightarrow{\sim} \U(\gl_N)^{\gen} \otimes V \sslash N_-, \qquad f(w) \otimes v \mapsto f(w) PvP,
\end{equation}
see \cite{ZhelobenkoSalgebras}.

\begin{Remark}
	\label{remark:transposition}
	Usually, the cited results use the positive version of the parabolic restriction \smash{$X \mapsto (X^{\gen}/\n_+)^{N_+}$}. For Miura transform of Section~\ref{sect:miura}, the negative one will be more convenient as the former is a map of \emph{left} $\U(\h)^{\gen}$-modules; however, both versions are equivalent thanks to the property $\varpi_N(P) = P$. For instance, assume $\varpi_N$ can be extended to an isomorphism of a Harish-Chandra bimodule $X$. Then there is an anti-isomorphism
	\[
	\varpi_N \colon\ (\n_- \backslash X^{\gen})^{N_-} \xrightarrow{\sim} (X^{\gen}/\n_+)^{N_+}, \qquad PxP \mapsto \varpi_N (PxP) = P\varpi_N(x) P.
	\]
\end{Remark}

\begin{Example}
	\label{ex:parabolic_res_center}
	If $V=\C$ is the trivial representation, then $\U(\g)^{\gen}\sslash N_- \cong \U(\h)^{\gen}$. In~particular, let us determine the image of the center $\Z_{\gl_N}$. Recall that by \eqref{eq::quantum_determinant}, its generators are given by the quantum minor $A(u)$. Using explicit formula \eqref{eq::quantum_minor} and the properties of the extremal projector, one can conclude that
	\begin{equation}
		\label{eq::parabolic_res_center}
		PA(u) = A(u)P = \prod_{i=1}^N (u-w_i),
	\end{equation}
	where, as usual, the first part is the image of $A(u)$ in the quotient $\n_- \backslash \U(\gl_N)$ and the second part is its image in the quotient $\U(\gl_N)/\n_+$ by the right action of $\n_+$, see \cite[Theorem~7.1.1]{MolevYangians}. In~particular, this defines an isomorphism $\Z_{\gl_N} \cong \U(\h)^W$, where $W = S_N$ acts by permutations of coordinates $\{w_i\}$. In terms of coordinates $\{E_{ii}\}$, it corresponds to the $\rho_N$-shifted action
	\[
		\sigma \bullet h := \sigma(h+\rho) - \rho, \qquad h\in \h,\ \sigma \in W.
	\]
	This is the standard Harish-Chandra isomorphism $\Z_{\gl_N} \cong \U(\h)^{W}$, see \cite[Theorem~1.10]{HumphreysBGG}.
\end{Example}

\subsection{Vector representation}
In this subsection, we describe explicitly the parabolic restriction of $\U(\gl_N) \otimes \C^N$. We will use a~slightly renormalized invariant vectors. Define
\begin{equation}
	\label{eq::parabolic_res_dual_difference_gens}
	\bar{v}_i = \sum\limits_{N \geq i_1 > \dots > i_s > i} P(w_i - w_{j_1}) \cdots (w_i - w_{j_r}) v_{i_1} E_{i_{2} i_1} \cdots E_{i i_s} \cdot \in \n_- \backslash \U(\gl_N) \otimes \C^N,
\end{equation}
where $\{j_1,\dots, j_r\}$
is the complementary subset to $\{i_1,\dots,i_s\}$ inside $\{i+1,\dots,N\}$. According to \cite[formulas~(3.2) and (3.3)]{MolevTransvectorAlgebras}, the trick of Example~\ref{ex:trick}, and Remark~\ref{remark:transposition}, those vectors are indeed invariant under $N_-$ and satisfy
\begin{equation}
	\label{eq::parabolic_res_vector_generators}
	\bar{v}_i = (w_i - w_{i+1}) \cdots (w_i - w_N) \cdot Pv_i P
\end{equation}
(more precisely, in the notations of loc.\ cit., they correspond to $\varpi_N(s_{\alpha i})$ for $\alpha=N+1$).

Using the parabolic restriction, one can deduce the following explicit description of the Harish-Chandra bimodule $\U(\gl_N) \otimes \C^N $ as bimodule over the center $\Z_{\gl_N}$. Before the statement, we need the following lemma which is essentially a part of \cite[Theorem~2]{KhoroshkinNazarovVinberg}; for reader's convenience, we repeat its proof in our notations.
\begin{Lemma}
	\label{lm:parabolic_res_injective}
	For any $\GL_N$-representations $V$ and $W$, the map
	\begin{gather*}
		\Hom_{\HC(\GL_N)} (\U(\gl_N) \otimes V, \U(\gl_N) \otimes W)\\
 \qquad \rightarrow \Hom_{\HC(H_N)^{\gen}} (\U(\gl_N) \otimes V\sslash N_-, \U(\gl_N) \otimes W \sslash N_- )
	\end{gather*}
	is injective.
	\begin{proof}
		Due to \eqref{eq::hom_hc_rep_g}, we have
		\[
		\Hom_{\HC(\GL_N)} (\U(\gl_N) \otimes V, \U(\gl_N) \otimes W) \cong (\U(\gl_N) \otimes W \otimes V^*)^{\GL_N}.
		\]
		For brevity, let us assume that $V = \C$. Let $s_{\hbar} \in (\U(\gl_N) \otimes W)^{\GL_N}$. Consider the induced homomorphism
		\[
		s \in (\U(\gl_N) \otimes W)^{\GL_N} / (\hbar) \cong (\cO(\gl_N^*) \otimes W)^{\GL_N}.
		\]
		In other words, $s$ is a $\GL_N$-equivariant section of the trivial vector bundle $\gl_N^* \times W \rightarrow \gl_N^*$. Since $\U(\gl_N) \otimes W$ is a free $\C[\hbar]$-module, we can assume that $s \neq 0$ by diving by powers of $\hbar$ if~necessary.
		
		Recall that parabolic restriction is isomorphic to the functor $X \mapsto \n_- \backslash X^{\gen} /\n_+$ as in \eqref{eq::extremal_proj_quotient_inv_iso}. For the category $\QCoh_{[\g^*/G]}$ as a subcategory of Harish-Chandra bimodules on which $\hbar$ acts by zero as in Remark~\ref{remark:hc_hbar_zero}, this functor is simply the restriction along the embedding $(\h^*)^{\mathrm{reg}} \hookrightarrow \g^*$. Since its $\GL_N$-orbit is dense, the restriction $\res(s)$ is non-zero as well. Since the parabolic restriction functor commutes with $\C[\hbar]$-action, we obtain that $\res(s_{\hbar})$ is non-zero too.
	\end{proof}
\end{Lemma}

The next lemma concerns the canonical homomorphism $\Omega$ of \eqref{eq::canonical_homomorphism} under the parabolic restriction.

\begin{Lemma}
	\label{lm:parabolic_res_canonical_homomorphism}
	For any $i$, we have $P\Omega(v_i)P = w_i Pv_iP$.
	\begin{proof}
		We have
		\[
		P \Omega(v_i) P = \sum_{b=1}^N PE_{ib} v_b P.
		\]
		Using the properties of the extremal projector \eqref{eq::extremal_proj}, we obtain
		\[
		\sum_{b>i} P E_{ib} v_b P = \sum_{b>i} P v_b E_{ib} P + \hbar(N-i) P v_i P = \hbar(N-i) P v_i P,
		\]
		as well as
		\[
		\sum_{b<i} P E_{ib} v_b P = 0, \qquad P E_{ii} v_i P = E_{ii} Pv_i P.
		\]
		Therefore, we have
		\[
			P\Omega(v_i) P = (E_{ii} + \hbar(N-i)) Pv_i P = w_i Pv_i P,
		\]
		as required.
	\end{proof}
\end{Lemma}

The equality in the following proposition is understood as in Remark~\ref{remark:pwoer_series}.

\begin{Proposition}
	\label{prop:conjugation_center}
	For any $v\in \C^N \subset \U(\gl_N) \otimes \C^N$, we have
	\[
	A(u)^{-1} v A(u) = \frac{u-\Omega+\hbar}{u-\Omega} v
	\]
	as elements of $\U(\gl_N) \otimes \C^N\bigl[\hspace{-0.2mm}\bigl[u^{-1}\bigr]\hspace{-0.2mm}\bigr]$.
	\begin{proof}
		Both sides are actions of certain elements of \smash{$\End_{\HC(\GL_N)} \bigl(\U(\gl_N) \otimes \C^N\bigr)\bigl[\hspace{-0.2mm}\bigl[u^{-1}\bigr]\hspace{-0.2mm}\bigr]$}. By Lemma~\ref{lm:parabolic_res_injective}, it is enough to show the equality under their parabolic restriction functor. It follows from Lemma~\ref{lm:parabolic_res_canonical_homomorphism} that
		\[
		\frac{u-\Omega+\hbar}{u-\Omega} Pv_i P = \frac{u-w_i+\hbar}{u-w_i} Pv_i P.
		\]
		At the same time, recall that by Example~\ref{ex:parabolic_res_center}, we have \smash{$P A(u) = \prod_{j=1}^{N} (u-w_j)$}, therefore,
		\[
		P v_i A(u) P = P\prod_{j=1}^{N} v_i (u-w_j) P = \prod_{j=1}^{N}(u - w_i + \delta_{ij} \hbar) \cdot P v_i P.
		\]
		Hence
		\[
		PA(u)^{-1} v_i A(u) P = \frac{u - w_i + \hbar}{u- w_i} Pv_i P,
		\]
		and the statement follows.
	\end{proof}
\end{Proposition}

\begin{Remark}
\label{remark:pwoer_series}
	Strictly speaking, we understand the equality of Proposition~\ref{prop:conjugation_center} in the sense of power series expansion in $u^{-1}$. The right-hand side can be presented as
	\[
	\frac{u-\Omega+\hbar}{u-\Omega} = 1 + \frac{\hbar}{u-\Omega} = 1 + \hbar \sum_{i=0}^{\infty} \Omega^i u^{-i-1} \in \End_{\HC(\GL_N)}\bigl(\U(\gl_N) \otimes \C^N\bigr)\bigl[\hspace{-0.2mm}\bigl[u^{-1}\bigr]\hspace{-0.2mm}\bigr],
	\]
	while the left-hand side can be equivalently presented as conjugation by $\qdet(T(-u+N\hbar-\hbar))$ in the notations of Section~\ref{subsect:quantum_minors} which also defines an element of \smash{$\End_{\HC(\GL_N)}\!\bigl(\U(\gl_N) \otimes \C^N\bigr)\!\bigl[\hspace{-0.2mm}\bigl[u^{-1}\bigr]\hspace{-0.2mm}\bigr]$}. Another way is via the coefficients \eqref{eq::quantum_determinant}: for any $v\in \C^N$, we have
	\[
	A_i v - A_{i-1}\Omega v = v A_i - (\Omega-\hbar) v A_{i-1}.
	\]
\end{Remark}

\begin{Remark}
	In fact, using Lemma~\ref{lm:parabolic_res_canonical_homomorphism} and \eqref{eq::parabolic_res_center}, one can re-deduce the Cayley--Hamilton theorem of Proposition~\ref{prop:cayley_hamilton}: observe that for every $i$, we have
	\[
	PA(\Omega) v_i P = \prod_{j=1}^N P (\Omega-w_j) v_i P = \prod_{j=1}^N (w_i -w_j) P v_i P = 0.
	\]
	Then we apply Lemma~\ref{lm:parabolic_res_injective}.
\end{Remark}

\subsection{Mickelsson algebras and difference operators}
\label{subsect:mickelsson}
By \cite[Corollary 4.18]{KalmykovSafronov}, the natural lax monoidal structure on $\res$ given by tensor multiplication
\begin{align}
\label{eq::par_res_tensor_structure}
\begin{split}
	X^{\gen} \sslash N_- \otimes_{\U(\h)^{\gen}} Y^{\gen} \sslash N_- & \rightarrow (X \otimes_{\U(\gl_N)} Y)^{\gen} \sslash N_-, \\
	[x] \otimes [y] &\mapsto \bigl[x\otimes_{\U(\gl_N)^{\gen}} y\bigr]
\end{split}
\end{align}
is an isomorphism. In particular, it gives the following classical construction.

\begin{Definition}[{\cite{ZhelobenkoSalgebras}}]
\label{def:mickelsson_algebra}
	Let $A\in \HC(\GL_N)$ be an algebra object. The {\it Mickelsson algebra} associated to $A$ is the $\U(\h)^{\gen}$-bimodule $A^{\gen} \sslash N_-$ with the multiplication given by
	\[
		A^{\gen} \sslash N_- \otimes_{\U(\h)^{\gen}} A^{\gen} \sslash N_- \xrightarrow{\sim} \bigl(A^{\gen} \otimes_{\U(\gl_N)^{\gen}} A^{\gen}\bigr) \sslash N_- \rightarrow A^{\gen} \sslash N_-,
	\]
	where the first arrow is \eqref{eq::par_res_tensor_structure} and the second arrow is multiplication on $A$. In the notations of~\eqref{eq::res_trivialization}, the product is denoted by
	\[
		PaP \otimes PbP \mapsto PaPbP,\qquad a,b \in A^{\gen}.
	\]
\end{Definition}

We will compute it explicitly for $\U(\gl_N) \ltimes S^{\bullet} \C^N$. Consider the algebra $\cA$ of \eqref{eq::diff_ops}. Denote by \smash{$\cA^{\pm}\subset \cA$} the $\U(\h)^{\gen}$-subalgebras generated by respectively positive or negative powers of $\{\bu_i\}$.

\begin{Theorem}
\label{thm:par_res_to_difference_ops}
	The left $\U(\h)^{\gen}$-module map
	\[
	\U(\gl_N)^{\gen} \ltimes S^{\bullet} \C^N \sslash N_- \xrightarrow{\sim} \cA^-, \qquad \bar{v}_i \mapsto \bu_i^{-1}
	\]
%	and
%	\[
%	\U(\gl_N)^{\gen} \ltimes S^{\bullet} \bigl(\C^N\bigr)^* \sslash N_- \xrightarrow{\sim} \cA^+, \qquad \bar{\phi}_i \mapsto \bu_i
%	\]
	is an algebra isomorphism.
	\begin{proof}
		By \eqref{eq::res_trivialization}, the monomials $Pv_1^{k_1} \cdots v_N^{k_N}P$ for all possible powers define a basis of the reduction $\U(\gl_N) \ltimes S^{\bullet} \C^N \sslash N_-$ over $\U(\h)^{\gen}$. By invertibility of the monoidal isomorphisms~\eqref{eq::par_res_tensor_structure}, the monomials
		\[
		(Pv_1 P)^{k_1} \cdots (Pv_N P)^{k_N}
		\]
		for all possible powers also constitute a basis. Therefore, by \eqref{eq::parabolic_res_vector_generators}, we conclude that, as a left $\U(\h)^{\gen}$-module, the reduction $\U(\gl_N)^{\gen} \ltimes S^{\bullet} \C^N \sslash N_-$ is a free symmetric algebra in $\{\bar{v}_i\}$.
		
		Since the extremal projector preserves weights and $v_i$ has weight $\epsilon_i$, we have
		\begin{equation}
			\label{prop_eq::function_commutation}
			f(w) \bar{v}_i = \bar{v}_i f(w + \hbar \epsilon_i)
		\end{equation}
		by \eqref{eq::hc_torus_right_action}. Let us show that $\bar{v}_i \bar{v}_j = \bar{v}_j \bar{v}_i$. Indeed, using \eqref{eq::parabolic_res_dual_difference_gens} and \eqref{eq::parabolic_res_vector_generators}, we have
		\[
		P v_i P = \sum \prod_{j\in I} (w_i - w_j)^{-1} v_{i_1} E_{I},
		\]
		where $E_{I} = E_{i_2 i_1} \cdots E_{i i_s} \in \U(\n_+)$ for $I = \{i_1 > \dots > i_s\}$. Consider the product $Pv_i P v_j P$. By pushing terms of $\n_+$ to the right using
		\[
		E_{ab} v_j P = v_j E_{ab} P + \hbar \delta_{bj} v_a P = \hbar \delta_{bj} v_a P
		\]
		for all $a<b$, we conclude that the only terms contributing non-trivially are those with $I=\varnothing$ and $|I|=1$, so that
		\[
		Pv_i Pv_j P = Pv_i v_j P + \delta_{j > i} \hbar (w_i - w_j)^{-1} Pv_j v_i P.
		\]
		Assume without loss of generality that $j>i$. Since in the symmetric algebra we have $v_i v_j = v_j v_i$, we obtain
		\[
		Pv_i P v_j P = \frac{w_i - w_j + \hbar}{w_i - w_j} Pv_j P v_i P.
		\]
		Recall that \smash{$Pv_i P = \prod_{j=i+1}^N (w_i - w_j)^{-1} \bar{v}_i$}. Substituting it in the equality above using the property $P^2 = P$, we conclude that $\bar{v}_i \bar{v}_j = \bar{v}_j \bar{v}_i$, in particular,
		\[
		\U(\gl_N)^{\gen} \ltimes S^{\bullet} \C^N \sslash N_- \cong \U(\h)^{\gen}[\bar{v}_1,\dots,\bar{v}_N]
		\]
		with relation \eqref{prop_eq::function_commutation}. Also, it is clear that \smash{$\cA^{-} \cong \U(\h)^{\gen{}}\bigl[\bu_1^{-1},\dots,\bu_N^{-1}\bigr]$} as left $\U(\h)^{\gen}$-modules. Therefore, the assignment $\bar{v}_i\mapsto \bu^{-1}_i$ extends to an algebra isomorphism
		\[
			\U(\gl_N)^{\gen} \ltimes S^{\bullet} \C^N \sslash N_- \xrightarrow{\sim} \cA^-,
		\]
		as claimed.
	\end{proof}
\end{Theorem}

Observe that the isomorphism is actually graded:
\[
	\U(\gl_N)^{\gen} \otimes S^k \C^N \sslash N_- \xrightarrow{\sim} \U(\h)^{\gen} \bigl[\bu_1^{-1},\dots,\bu_N^{-1}\bigr]_k,
\]
where the target is the space of polynomials of total degree $k$ in $\bu^{-1}$ variables. By Corollary~\ref{cor:action_spherical_on_symmetric}, the source has an action of the spherical subalgebra ${\rm e}^{\hbar}_k \fH^{\hbar}_k {\rm e}^{\hbar}_k$ of the degenerate affine Hecke algebra.

\begin{Proposition}
\label{prop:spherical_action_par_res}
	The spherical subalgebra $\C[\hbar][Y_1,\dots,Y_k]^{S_k}$ acts on localized difference operators $\U(\h)^{\gen} \bigl[\bu_1^{-1},\dots,\bu_N^{-1}\bigr]_k$ by
	\[
		f\bigl(\vec{Y}\bigr) \cdot \bu_{i_1}^{-1} \cdots \bu_{i_k}^{-1} = \sum_J c_J \bigl(w_{i_1}^{j_1} \bu^{-1}_{i_1}\bigr) \cdots (w_{i_k}^{j_k} \bu^{-1}_{i_k}), \qquad f\bigl(\vec{Y}\bigr) = \sum_{J=(j_1,\dots,j_k)} c_J Y_1^{j_1} \cdots Y_k^{j_k}
	\]
	and commutes with the left $\U(\h)^{\gen}$-action.
	\begin{proof}
		Recall that $\bu_{i_1}^{-1} \cdots \bu_{i_k}^{-1}$ is proportional to $Pv_{i_1} P \cdots P v_{i_k}P$. By Corollary~\ref{cor:action_spherical_on_symmetric}, the action of a monomial on the latter is
		\[
			Y_1^{j_1} \cdots Y_k^{j_k}\cdot Pv_{i_1} P \cdots P v_{i_k}P = P\Omega^{j_1}(v_{i_1}) P \cdots P \Omega^{j_k}(v_{i_k}) P
		\]
		(it is well-defined once the function is symmetric). The result follows by Lemma~\ref{lm:parabolic_res_canonical_homomorphism}.
	\end{proof}
\end{Proposition}

\section{Kostant--Whittaker reduction and Yangians}
\label{sect:kw_reduction}

In this section, we will study the Whittaker variant of the constructions from Section~\ref{sect:parabolic_reduction}. We use the action of degenerate affine Hecke algebras from Section~\ref{sect:hc_bimodules} to find a presentation of the Kostant--Whittaker reduction of $\U(\gl_N) \otimes S^k \C^N$ as a certain quotient of symmetric polynomial in $k$ variables over the center $\Z_{\gl_N}$, see Corollary~\ref{cor:kw_symmetric_power}. It is used to construct a homomorphism from the non-negative part $\rY_{-N,N}(\sl_2)^{\geq 0}$ of the shifted Yangian to the Kostant--Whittaker reduction of $\U(\gl_N) \ltimes S^{\bullet} \C^N$ by using the results of Section~\ref{sect:deg_aha_yangians}. The main construction of the paper, Theorem~\ref{thm:toda_lattice_yangian}, is essentially a combination of the one from this section with the corresponding dual counterpart of Appendix~\ref{subsect:dual_rep}. In Section~\ref{subsect:mirabolic}, we discuss a mirabolic analog of the results from previous subsections.

Let $\psi \in \C^*$ be a non-zero number. Somewhat abusing notations, we denote by
\begin{equation}
\label{eq::nilp_character}
	\psi\colon\ \n_+ \rightarrow \C, \qquad \psi(E_{ij}) = \delta_{i+1,j} \psi
\end{equation}
a non-degenerate character of $\n_+$. For any $x\in \n_+$, denote by $x^{\psi} := x - \psi(x)$. Consider the shift
\[
\n_+^{\psi}:= \{ x - \psi(x) \mid x\in \n_+\} \subset \U(\gl_N).
\]
For any $X\in \HC(\GL_N)$, denote by $X/\n_+^{\psi}$ the quotient by the right action of $\n_+^{\psi}$. Recall that by Kostant's construction \cite{Kostant}, the space of invariants \smash{$\U(\gl_N) \hamp N_+ := \bigl(\U(\gl_N)/\n_+^{\psi}\bigr)^{N_+}$} is isomorphic to the center $\Z_{\gl_N}$ of $\U(\gl_N)$. It admits the following categorification.
\begin{Definition}[{\cite{BezrukavnikovFinkelberg}}]
\label{def:kw_reduction}
	The {\it Kostant--Whittaker reduction} is a functor
	\[
	\res^{\psi}_+ \colon\ \HC(\GL_N) \rightarrow \BiMod{\Z_{\gl_N}}{\Z_{\gl_N}}, \qquad X \mapsto X \hamp N_+ := \bigl(X/\n_+^{\psi}\bigr)^{N_+}.
	\]
\end{Definition}

There is a Whittaker analog of the extremal projector \eqref{eq::extremal_proj} introduced in \cite{KalmykovYangians} under the name of the Kirillov projector.
\begin{Theorem}
\label{eq::kirillov_proj_properties}
	There is an element $\pP$ in a certain completion of $\U(\gl_N)$ acting on any locally nilpotent $\n_+^{\psi}$-module $F$ such that $\n_+^{\psi} F \subset \hbar F$ and satisfying the properties
	\begin{align*}
		\begin{split}
			&(x - \psi(x)) \pP= 0,\qquad x\in \n_+, \\
			&\pP (E_{ij} + \delta_{ij} \hbar (N-i))= 0,\qquad 1 \leq j \leq i \leq N-1, \\
			&(\pP)^2 = \pP,\qquad \pP|_{(X/\n_+^{\psi})^{N_+}} = \mathrm{id}.
		\end{split}
	\end{align*}
	\begin{proof}
		In the notations of \cite[Definition 5.2]{KalmykovYangians}, we can apply \eqref{eq::transposition} to $P_{\m_N}^{\psi}(N\hbar-\hbar,\dots,\hbar)$; by \cite[Theorem~5.3]{KalmykovYangians}, it satisfies the properties of the theorem for $\psi(E_{ij}) = \delta_{i+1,j}$. For general $\psi\in \C^*$, we can apply the automorphism of $\U(\gl_N)$ defined by $x \mapsto x/\psi$ for $x\in\gl_N$ and $\hbar \mapsto \hbar/\psi$.
	\end{proof}
\end{Theorem}

As in Section~\ref{sect:parabolic_reduction}, we use the following.
\begin{Notation}
	For any $X \in \HC(\GL_N)$, we denote by
	\[
		X \rightarrow X/\n_+^{\psi}, \qquad x \mapsto x \pP
	\]
	the image of $x\in X$ under the projection map. This is well-defined, since $x(\xi - \psi(\xi)) \pP=0$ for any $\xi \in \n_+$.
\end{Notation}

The Kirillov projector defines an isomorphism
\begin{equation}
\label{eq::kw_kirillov_proj}
	\mathfrak{b}^{\rho}_{\mathrm{mir}} \backslash X / \n_+^{\psi} \xrightarrow{\sim} X \hamp N_+, \qquad [x] \mapsto \pP x \pP,
\end{equation}
where
\[
	\mathfrak{b}^{\rho}_{\mathrm{mir}} := \mathrm{span} (E_{ij} + \delta_{ij} \cdot \hbar(N-i) \mid 1 \leq j \leq i \leq N-1 ).
\]
For every $V\in \Rep(\GL_N)$, we also have an analog of the trivialization \eqref{eq::res_trivialization}, a left $\Z_{\gl_N}$-module isomorphism
\begin{equation}
\label{eq::kw_trivialization}
	\Z_{\gl_N} \otimes V \xrightarrow{\sim} \U(\gl_N) \otimes V \hamp N_{\pm}, \qquad v \mapsto \npP v \npP.
\end{equation}
We can compute explicitly this trivialization for the vector representation $\C^N$ using the canonical homomorphism \eqref{eq::canonical_homomorphism}. The following can be proved analogously to \cite[Proposition 6.18]{KalmykovYangians}.
\begin{Lemma}
	\label{lm:kirillov_projector_canonical_homomorphism}
	For any $1 \leq k \leq N$, we have
	\[
	\pP v_i \pP = \psi^{1-i} \pP \Omega^{i-1} (v_1)\pP.
	\]
\end{Lemma}

It implies another description of $\U(\gl_N) \otimes \C^N \hamp N_+$.
\begin{Proposition}\quad
\label{prop:kw_vector}
	\begin{enumerate}\itemsep=0pt
		\item[$1.$] Action on $v_1$ gives a left $\Z_{\gl_N}$-module isomorphism
		\[
		\Z_{\gl_N} [\Omega]/(A(\Omega)) \xrightarrow{\sim} \U(\gl_N) \otimes \C^N \hamp N_+, \qquad f \mapsto f(v_1).
		\]
		
		\item[$2.$] There is an algebra isomorphism
		\[
		\Z_{\gl_N}[\Omega]/(A(\Omega)) \xrightarrow{\sim} \End_{\BiMod{\Z_{\gl_N}}{\Z_{\gl_N}}} \bigl(\U(\gl_N) \otimes \C^N \hamp N_+\bigr).
		\]
	\end{enumerate}
	\begin{proof}
		Combining \eqref{eq::kw_trivialization}, Lemma~\ref{lm:kirillov_projector_canonical_homomorphism} and Proposition~\ref{prop:cayley_hamilton}, we conclude that the map
		\[
		\Z_{\gl_N}[\Omega]/(A(\Omega)) \rightarrow \U(\gl_N) \otimes \C^N \hamp N_+, \qquad f \mapsto f(v_1)
		\]
		is an isomorphism of modules over $\Z_{\gl_N}[\Omega]$. By Proposition~\ref{prop:conjugation_center}, we have
		\[
		vA(u) = A(u) \frac{u-\Omega+\hbar}{u-\Omega} v.
		\]
		By expanding in powers of $u$ as in Remark~\ref{remark:pwoer_series}, we can conclude that an endomorphism of $\U(\gl_N) \otimes \C^N \hamp N_+$ commuting with the left $\Z_{\gl_N}$-action commutes with the right $\Z_{\gl_N}$-action if and only if it commutes with $\Omega$. Therefore, such an endomorphism commutes with the \emph{algebra} action of $\Z_{\gl_N}[\Omega]/(A(\Omega))$, giving a canonical isomorphism
		\[
		\Z_{\gl_N}[\Omega]/(A(\Omega)) \xrightarrow{\sim} \End_{\BiMod{\Z_{\gl_N}}{\Z_{\gl_N}}} \bigl(\U(\gl_N) \otimes \C^N \hamp N_+\bigr),
		\]
		as claimed.
	\end{proof}
\end{Proposition}

\subsection{Kostant--Mickelsson--Whittaker algebras}

Similarly to \eqref{eq::par_res_tensor_structure}, the Kostant--Whittaker reduction is endowed with a monoidal structure via
\begin{align}
\label{eq::kw_tensor_structure}
\begin{split}
	X \hamp N_+ \otimes_{\Z_{\gl_N}} Y \hamp N_+ &\xrightarrow{\sim} \bigl(X \otimes_{\U(\gl_N)} Y\bigr) \hamp N_+, \\
	[x] \otimes [y] &\mapsto \bigl[x \otimes_{\U(\gl_N)} y\bigr],
\end{split}
\end{align}
see \cite{BezrukavnikovFinkelberg}. In particular, we have an analog of Definition~\ref{def:mickelsson_algebra}.

\begin{Definition}
\label{def:kmw_algebra}
	Let $A$ be an algebra object in $\HC(\GL_N)$. The {\it Kostant--Mickelsson--Whittaker algebra} associated to $A$ is the $\Z_{\gl_N}$-bimodule $A \hamp N_+$ with the product
	\[
		A \hamp N_+ \otimes_{\Z_{\gl_N}} A \hamp N_+ \xrightarrow{\sim} \bigl(A \otimes_{\U(\gl_N)} A\bigr) \hamp N_+ \rightarrow A \hamp N_+,
	\]
	where the first arrow is the map \eqref{eq::kw_tensor_structure} and the second arrow is multiplication on $A$. In terms of identification \eqref{eq::kw_kirillov_proj}, the multiplication is denoted by
	\[
		\npP a \npP \otimes \npP b \npP \mapsto \npP a \npP b \npP,\qquad a,b\in A.
	\]
\end{Definition}

We explicitly compute the Kostant--Mickelsson--Whittaker algebras $\U(\gl_N) \ltimes T\bigl(\C^N\bigr) \hamp N_+$, where \smash{$T\bigl(\C^N\bigr) := \bigoplus_{i=0}^{\infty} \bigl(\C^N\bigr)^{\otimes i}$} is the tensor algebra, and \smash{$\U(\gl_N) \ltimes S^{\bullet} \C^N \hamp N_+$}.

Recall the action of the degenerate affine Hecke algebra $\fH^{\hbar}_k$ on \smash{$\U(\gl_N) \otimes \bigl(\C^N\bigr)^{\otimes k}$} from Section~\ref{subsect:deg_affine_hecke_action}. It descends to an action on \smash{$\U(\gl_N) \otimes \bigl(\C^N\bigr)^{\otimes k} \hamp N_+$}. Recall the central elements \smash{$A(u) = \sum_i A_i u^{N-i}$} of \eqref{eq::quantum_determinant}. Define the elements \smash{$A^{(j)}_i \in \Z_{\gl_N} [\Omega_1,\dots,\Omega_k]$} for $1\leq j \leq k+1$ inductively via \smash{$A^{(1)}_i = A_i$} for all $i$ and
\begin{equation*}
	A_i^{(j+1)} - \Omega_{j} A_{i-1}^{(j+1)} = A_i^{(j)} - (\Omega_j-\hbar) A_{i-1}^{(j)}, \qquad A_{i}^{(j)} = \delta_{i0},\qquad i\leq 0.
\end{equation*}
Using Remark~\ref{remark:pwoer_series}, it can be rewritten as
\begin{equation}
\label{eq::modified_center}
	A^{(j+1)}(u) = \frac{u-\Omega_j+\hbar}{u-\Omega_j} A^{(j)}(u), \qquad A^{(j)}(u) := \sum_{i=0}^N A^{(j)}_i u^{N-i}.
\end{equation}
In particular, we have
\begin{equation}
%\label{eq::tensor_center_gens}
	A^{(j)}(u) = A(u) \prod_{a=1}^{j-1} \frac{u-\Omega_a + \hbar}{u-\Omega_a}.
\end{equation}

Consider the algebra $\bigoplus_{k=0}^{\infty} \Z_{\gl_N} [\Omega_1,\dots,\Omega_k]$ with product
\begin{align}
	\label{eq::pre_shuffle_product_center}
	&f(\Omega_1,\dots,\Omega_k) \otimes A(u)\mapsto A^{(k+1)}(u) f(\Omega_1,\dots,\Omega_k), \\
	\label{eq::pre_shuffle_product_poly}
	&f(\Omega_1,\dots,\Omega_k) \otimes g(\Omega_1,\dots,\Omega_{k+l})\mapsto f(\Omega_1,\dots,\Omega_k) g(\Omega_{k+1},\dots,\Omega_{k+l})
\end{align}
for $f(\Omega_1,\dots,\Omega_k) \in \Z_{\gl_N}[\Omega_1,\dots,\Omega_k]$ and $g(\Omega_1,\dots,\Omega_l) \in \Z_{\gl_N}[\Omega_1,\dots,\Omega_l]$ regarded as elements in the corresponding graded components.

The following is a generalization of Proposition~\ref{prop:kw_vector}.

\begin{Proposition}\qquad
\label{prop:kw_tensor_algebra}
	\begin{enumerate}\itemsep=0pt
		\item[$1.$] The action on $v_1^{\otimes k}$ gives an isomorphism of left $\Z_{\gl_N}$-modules
		\[
		\Z_{\gl_N} [\Omega_1,\dots,\Omega_k]/\bigl(A^{(j)}(\Omega_j) \mid j=1, \dots, k\bigr) \xrightarrow{\sim} \U(\gl_N) \otimes \bigl(\C^N\bigr)^{\otimes k} \hamp N_+,
		\]
		compatible with the action of $\fH^{\hbar}_k$.
		
		\item[$2.$] There is an algebra isomorphism
		\[
			 \bigoplus_{k=0}^{\infty} \Z_{\gl_N} [\Omega_1,\dots,\Omega_k]/\bigl(A^{(j)}(\Omega_j) \mid j=1,\dots, k \bigr) \xrightarrow{\sim} \U(\gl_N) \ltimes \bigoplus_{k=0}^{\infty} \bigl(\C^N\bigr)^{\otimes k} \hamp N_+,
		\]
		where multiplication on the right-hand side is descended from \eqref{eq::pre_shuffle_product_center} and \eqref{eq::pre_shuffle_product_poly}.
		
		\item[$3.$] Denote by
		\[
			\End_{\Z_{\gl_N}}^{\C[\hbar][\Omega_1,\dots,\Omega_k]^{S_k}}\bigl(\Z_{\gl_N} \otimes \bigl(\C^N\bigr)^{\otimes k}\bigr)
		\]
		the space of $\Z_{\gl_N}$-linear endomorphism of the free left $\Z_{\gl_N}$-module \smash{$\Z_{\gl_N} \otimes \bigl(\C^N\bigr)^{\otimes k}$} commuting with the action of the spherical subalgebra \smash{$\C[\hbar][\Omega_1,\dots,\Omega_k]^{S_k}$} of $\fH^{\hbar}_k$. Then we have an isomorphism
		\begin{gather*}
			\End_{\Z_{\gl_N}}^{\C[\hbar][\Omega_1,\dots,\Omega_k]^{S_k}}\bigl(\Z_{\gl_N} \otimes \bigl(\C^N\bigr)^{\otimes k}\bigr) \cong \End_{\BiMod{\Z_{\gl_N}}{\Z_{\gl_N}}} \bigl(\U(\gl_N) \otimes \bigl(\C^N\bigr)^{\otimes k} \hamp N_+\bigr).
		\end{gather*}
	\end{enumerate}
	\begin{proof}
		For the first part, observe that we have an isomorphism
		\begin{align*}
 &\U(\gl_N) \otimes \bigl(\C^N\bigr)^{\otimes k} \hamp N_+ \\
&\qquad \cong \bigl(\U(\gl_N) \otimes \C^N \hamp N_+\bigr) \otimes_{\Z_{\gl_N}} \dots \otimes_{\Z_{\gl_N}} \bigl(\U(\gl_N) \otimes \C^N \hamp N_+\bigr) \\
			&\qquad \cong \Z_{\gl_N} [\Omega_1]/(A(\Omega_1)) \otimes_{\Z_{\gl_N}} \dots \otimes_{\Z_{\gl_N}} \Z_{\gl_N} [\Omega_k]/(A(\Omega_k)).
		\end{align*}
		By Proposition~\ref{prop:conjugation_center}, the right action of the center on each tensor component satisfies
		\[
			v A(u) = A(u) \frac{u-\Omega_i+\hbar}{u-\Omega_i}.
		\]
		Therefore, by pushing the central generators to the left using this identity, we obtain the desired isomorphism
		\[
			\U(\gl_N) \otimes \bigl(\C^N\bigr)^{\otimes k} \hamp N_+
			\cong \Z_{\gl_N} [\Omega_1,\dots,\Omega_k]/\bigl(A^{(j)}(\Omega_j) \mid j=1,\dots,k\bigr)
		\]
		by definition \eqref{eq::modified_center}.
		
		For the second part, observe that the right action on the right-hand side can be expressed as
		\[
			A^{(k+1)}(u) = A(u) \prod_{i=1}^k \frac{u-\Omega_i+\hbar}{u-\Omega_i}.
		\]
		Therefore, the equation \eqref{eq::pre_shuffle_product_center} is satisfied. Likewise, the equation \eqref{eq::pre_shuffle_product_poly} holds by definition.
		
		For the third part, an endomorphism commuting with the left $\Z_{\gl_N}$-action commutes with the right one if and only if it commutes with the coefficients of the product
		\[
			\prod_{i=1}^k \frac{u-\Omega_i+\hbar}{u-\Omega_i}.
		\]
		By Lemma~\ref{lm:symmetric_func}, they generate the ring of symmetric functions in $\{\Omega_i\}$. Therefore, we have a~desired isomorphism.
	\end{proof}
\end{Proposition}

\begin{Corollary}
\label{cor:kw_symmetric_power}
	Denote by
	\[
	I_k := \bigl(A^{(j)}(\Omega_j)\mid j=1,\dots, k\bigr) \cap \Z_{\gl_N}[\Omega_1,\dots,\Omega_k]^{S_k}
	\]
	the ideal inside $\Z_{\gl_N}[\Omega_1,\dots,\Omega_k]^{S_k}$. For any $k$, there is an isomorphism
	\[
	\U(\gl_N) \otimes S^k\C^N \hamp N_+ \cong \Z_{\gl_N}[\Omega_1,\dots,\Omega_k]^{S_k}/I_k.
	\]
	In particular, as vector space,
	\[
	\U(\gl_N) \otimes S^{\bullet} \C^N \hamp N_+ \cong \bigoplus_{k=0}^{\infty} \Z_{\gl_N} [\Omega_1,\dots,\Omega_k]^{S_k} /I_k.
	\]
	\begin{proof}
		By naturality, the subspace
		\[
		\U(\gl_N) \otimes S^k \C^N \hamp N_+ \subset \U(\gl_N) \otimes \bigl(\C^N\bigr)^{\otimes k} \hamp N_+
		\]
		can identified with $S_k$-invariants. Since $S_k$ is a finite group, the image of the map
		\[
		\Z_{\gl_n} [\Omega_1,\dots,\Omega_k]^{S_k} \rightarrow \Z_{\gl_N} [\Omega_1,\dots,\Omega_k]/\bigl(A^{(j)}(\Omega_j) \mid j=1,\dots,k\bigr)
		\]
		is equal to $S_k$-invariants of the target. However, it can be identified with $\Z_{\gl_N}[\Omega_1,\dots,\Omega_k]^{S_k} / I_k$ by definition.
	\end{proof}
\end{Corollary}

The following theorem relates the positive part of the shifted Yangian from Definition~\ref{def:shifted_yangian} to the Kostant--Mickelsson--Whittaker algebra associated to $S^{\bullet} \C^N$. Observe similarity of the explicit formulas below to the quantization of Zastava spaces for $\mathrm{SL}_2$ constructed in \cite{FinkelbergRybnikov}.

\begin{Theorem}
\label{thm:kw_to_yangian}
	There is a surjective homomorphism of $\C[\hbar]$-algebras
	\[
	\rY_{-N,N}(\gl_2)^{\geq 0} \rightarrow \U(\gl_N) \ltimes S^{\bullet} \C^N \hamp N_+,
	\]
	defined by
	\begin{align*}
		d_1(u) &\mapsto A(u), \qquad d_2(u) \mapsto A(u-\hbar)^{-1}, \\
		x_+(u) &\mapsto -\frac{1}{u-\Omega_1} (v_1) =- \sum_{i=0}^{\infty} \Omega_1^i(v_1) u^{-i-1} \in \bigl(\U(\gl_N) \otimes \C^N \hamp N_+ \bigr)\bigl[\hspace{-0.2mm}\bigl[u^{-1}\bigr]\hspace{-0.2mm}\bigr].
	\end{align*}
	\begin{proof}
		Relation \eqref{eq::yangian_d_d_comm} is clear. Let us verify \eqref{eq::yangian_d_e_comm} that reads
		\[
		-A(z) \frac{1}{w-\Omega_1} (v_1) + \frac{1}{w-\Omega_1} (v_1) A(z) = -\frac{\hbar}{z-w} A(z) \biggl(\frac{1}{z-\Omega_1} - \frac{1}{w-\Omega_1}\biggr) (v_1).
		\]
		Indeed, the right-hand side is simply
		\[
		A(z) \frac{\hbar}{(z-\Omega_1)(w-\Omega_1)} (v_1),
		\]
		and, due to Proposition~\ref{prop:conjugation_center}, the left-hand side is
		\[
		-A(z) \frac{1}{w-\Omega_1} (v_1) + A(z) \frac{1}{w-\Omega_1} (v_1) \biggl( 1 + \frac{\hbar}{z-\Omega_1}\biggr) = A(z) \frac{\hbar}{(z-\Omega_1)(w-\Omega_1)}(v_1).
		\]
		The relation for $d_2(z)$ can be verified similarly.
		
		As for the relation \eqref{eq::yangian_e_e_comm}, observe that due to the second part of Proposition~\ref{prop:kw_tensor_algebra} and Corollary~\ref{cor:kw_symmetric_power}, the product
		\[
			\bigl(\U(\gl_N) \otimes S^k \C^N \hamp N_+\bigr) \otimes \bigl(\U(\gl_N) \otimes S^l \C^N \hamp N_+\bigr) \rightarrow \U(\gl_N) \otimes S^{k+l} \C^N \hamp N_+
		\]
		is induced from
		\begin{align*}
			\Z_{\gl_N}[\Omega_1,\dots,\Omega_k]^{S_k} \otimes \Z_{\gl_N}[\Omega_1,\dots,\Omega_l]^{S_l} &\rightarrow \Z_{\gl_N}[\Omega_1,\dots,\Omega_{k+l}]^{S_{k+l}}, \\
			 f(\Omega_1,\dots,\Omega_k) \otimes g(\Omega_1,\dots,\Omega_l) &\mapsto {\rm e}^{\hbar}_{k+l} \cdot f(\Omega_1,\dots,\Omega_k) g(\Omega_{k+1},\dots,\Omega_{k+l}),
		\end{align*}
		where \smash{${\rm e}^{\hbar}_{k+l}$} is the symmetrizer of Section~\ref{subsect:deg_affine_hecke}. However, upon extending the scalars to $\Z_{\gl_N}$, this is exactly the product of the shuffle algebra \smash{$\cF^{\hbar}$} of Definition~\ref{def:shuffle_algebra}, in other words, the ``quotient'' map
		\[
			\cF^{\hbar} = \bigoplus_{k=0}^{\infty} \C[\hbar][\Omega_1,\dots,\Omega_k]^{S_k} \rightarrow \bigoplus_{k=0}^{\infty} \Z_{\gl_N}[\Omega_1,\dots,\Omega_k]^{S_k} /I_k \cong \U(\gl_N) \ltimes S^{\bullet} \C^N \hamp N_+
		\]
		is an algebra map. By Theorem~\ref{thm:yangian_hecke_shuffle}, there is an algebra isomorphism $Y(\gl_2)^{+} \xrightarrow{\sim} \cF^{\hbar}$. Twisting the resulting composition $Y(\gl_2)^+ \rightarrow \U(\gl_N) \ltimes S^{\bullet} \C^N \hamp N_+$ by negative sign (which does not influence the defining relations) sends $x_+(u)$ to $-\frac{1}{u-\Omega_1} (v_1)$.
	\end{proof}
\end{Theorem}

\subsection{Other representations}

In this subsection, we compare the isomorphism of Proposition~\ref{prop:kw_tensor_algebra} with the trivialization \eqref{eq::kw_trivialization} provided by the Kirillov projector.

\begin{Proposition}
\label{prop:kw_triv_tensor_prod}
	The natural isomorphism of \eqref{eq::kw_trivialization}
	\[
		\varphi\colon\ \Z_{\gl_N} \otimes \C^N \otimes \dots \otimes \C^N \rightarrow \Z_{\gl_N} [\Omega_1,\dots,\Omega_k]/\bigl(A^{(j)}(\Omega_j) \mid j=1,\dots,k\bigr)
	\]
	of left $\Z_{\gl_N}$-modules is determined by two conditions: for every multi-index $\{N \geq i_1 \geq \dots \geq i_k \geq 1 \}$, we have
	\[
		\varphi(v_{i_1} \otimes \dots \otimes v_{i_k}) = \psi^{k - i_1-\dots - i_k} \Omega_1^{i_1-1} \cdots \Omega^{i_k-1}_k,
	\]
	and for every $\sigma \in S_k$
	\[
		\varphi(v_{i_{\sigma(1)}} \otimes \dots \otimes v_{i_{\sigma(k)}}) = \sigma^{\hbar} \cdot \varphi(v_{i_1} \otimes \dots \otimes v_{i_k} ),
	\]
	where the action on the right is that of the degenerate affine Hecke algebra.
	\begin{proof}
		In what follows, we will drop the tensor product sign by considering \smash{$\U(\gl_N) \otimes \bigl(\C^N\bigr)^{\otimes k}$} as a submodule of \smash{$\U(\gl_N) \ltimes \bigoplus_{i=0}^{\infty} \bigl(\C^N\bigr)^{\otimes i}$} and using notations of Definition~\ref{def:kmw_algebra}.
		
		The second part follows by construction of the action of $\fH^{\hbar}_k$. As for the first part, we essentially need to show that
		\[
			\pP v_{i_1} \pP \cdots \pP v_{i_k} \pP = \pP v_{i_1} \cdots v_{i_k} \pP
		\]
		for sequences $\{i_1 \geq \dots \geq i_k\}$ since
		\[
			\pP v_{i_1} \pP \cdots \pP v_{i_k} \pP = \psi^{k-i_1-\dots-i_k} \pP \Omega_1^{i_1-1} \cdots \Omega_k^{i_k-1}\bigl(v_1^{\otimes k}\bigr) \pP.
		\]
		We will prove it by induction on each term of the sequence $\{i_1 \geq \dots \geq i_k\}$. The base is
		\[
			\pP v_1 \pP \cdots \pP v_1 \pP = \pP v_1 \cdots v_1 \pP,
		\]
		since $v_1 \pP$ is already a Whittaker vector. Assume we showed
		\[
		\pP v_{i_1} \pP \cdots \pP v_{i_k} \pP = \pP v_{i_1} \cdots v_{i_k} \pP.
		\]
		Take $1\leq \alpha \leq k$ such that $i_{\alpha-1} \geq i_{\alpha} +1 \geq i_{\alpha+1}$ (i.e., the non-increasing condition is still~sat\-is\-fied). By Lemma~\ref{lm:kirillov_projector_canonical_homomorphism}, we have
		\[
			\pP v_{i_1} \pP \cdots \pP v_{i_{\alpha}+1} \pP \cdots \pP v_{i_k} \pP = \psi^{-1} \pP v_{i_1} \pP \cdots \pP \Omega_{\alpha}(v_{i_{\alpha}}) \pP \cdots \pP v_{i_k} \pP.
		\]
		By induction assumption, the right-hand side multiplied by $\psi$ is equal to
		\[
			\pP v_{i_1} \cdots \Omega_{\alpha}(v_{i_{\alpha}}) \cdots v_{i_k} \pP = \pP v_{i_1} \cdots v_{i_{\alpha-1}} \Biggl(\sum_{\beta=1}^N E_{i_{\alpha} \beta} v_{\beta}\Biggr) v_{i_{\alpha+1}} \cdots v_{i_k} \pP.
		\]
		Consider the sum for $\beta > i_{\alpha}$. Observe that due to $\{i_1 \leq \dots \leq i_{\alpha}\}$, the terms $E_{i_{\alpha}\beta}$ commute with the head $v_{i_{\alpha+1}} \cdots v_{i_k}$, hence we have
		\begin{align*}
			&\pP v_{i_1} \cdots v_{i_{\alpha-1}} \biggl(\sum_{\beta > i_{\alpha}} E_{i_{\alpha} \beta} v_{\beta}\biggr) v_{i_{\alpha+1}} \cdots v_{i_k} \pP \\
			& \qquad{} = \pP v_{i_1} \cdots v_{i_{\alpha-1}} \biggl(\sum_{\beta > i_{\alpha}} v_{\beta} v_{i_{\alpha+1}} \cdots v_{i_k} E_{i_{\alpha} \beta}\biggr)\pP+ \hbar(N-i_{\alpha}) \pP v_{i_1} \cdots v_{i_{\alpha}} \cdots v_{i_k} \pP \\
			& \qquad{} = \psi \pP v_{i_1} \cdots v_{i_{\alpha}+1} \cdots v_{i_k} \pP + \hbar(N-i_{\alpha}) \pP v_{i_1} \cdots v_{i_{\alpha}} \cdots v_{i_k} \pP
		\end{align*}
		due to the properties of \eqref{eq::kirillov_proj_properties} of the Kirillov projector. Consider the sum for $\beta \leq i_{\alpha}$. Since we assumed that $i_{1} \geq \dots \geq i_{\alpha-1} > i_{\alpha}$, the elements $E_{i_{\alpha}\beta}$ also commute with the tail $v_{i_1} \cdots v_{i_{\alpha-1}}$. Using again the properties of the Kirillov projector, we have
		\[
			\pP v_{i_1} \cdots v_{i_{\alpha-1}} \biggl(\sum_{\beta \leq i_{\alpha}} E_{i_{\alpha} \beta} v_{\beta}\biggr) v_{i_{\alpha+1}} \cdots v_{i_k} \pP = -\hbar (N-i_{\alpha}) \pP v_{i_1} \cdots v_{i_{\alpha}} \cdots v_{i_k} \pP.
		\]
		Taking the sum, we obtain
		\[
			\pP v_{i_1} \pP \cdots \pP v_{i_{\alpha}+1} \pP \cdots \pP v_{i_k} \pP = \pP v_{i_1} \cdots v_{i_{\alpha}+1} \cdots v_{i_k} \pP,
		\]
		as required.
	\end{proof}
\end{Proposition}

Therefore, for every irreducible representation $V$ with an embedding $V \hookrightarrow \bigl(\C^N\bigr)^{\otimes k}$, there is an explicit realization of the corresponding Kostant--Whittaker reduction
\[
	\Z_{\gl_N} \otimes V \cong \U(\gl_N) \otimes V \hamp N_+ \rightarrow \Z_{\gl_N}[\Omega_1,\dots,\Omega_k]/\bigl(A^{(j)}(\Omega_j) \mid j=1,\dots, k\bigr).
\]
By construction, it restricts to a $\C[\hbar]$-linear map
\[
	V[\hbar] \rightarrow \C[\hbar][\Omega_1,\dots,\Omega_k]_{<N},
\]
where the target is the space of polynomials of degree less than $N$ in each variable. In the next subsection, we will give its canonical interpretation.

\subsection{Mirabolic reduction and Feigin--Odesskii shuffle algebras}
\label{subsect:mirabolic}

Denote by $\m_N \subset \gl_N$ the \emph{mirabolic subalgebra}:
\begin{equation*}
\m_N = \mathrm{span}(E_{ij} \mid i=1,\dots, N-1,\, j=1,\dots, N ).	
\end{equation*}
Observe that $\n_+ \subset \m_N$. In particular, there is a well-defined functor
\begin{equation}
\label{eq::kw_mirabolic}
	\Rep(\GL_N) \rightarrow \Mod_{\C[\hbar]}, \qquad V \mapsto \U(\m_N) \otimes V \hamp N_+ := \bigl(\U(\m_N) \otimes V /\n_+^{\psi}\bigr)^{N_+},
\end{equation}
a mirabolic version of the Kostant--Whittaker reduction from Definition~\ref{def:kw_reduction} for free Harish-Chandra bimodules. It follows from \cite[Theorem~5.3]{KalmykovYangians} that the action of the Kirillov projector of \eqref{eq::kirillov_proj_properties} is well-defined on \smash{$\U(\m_N) \otimes V /\n_+^{\psi}$} and gives an isomorphism
\begin{equation}
\label{eq::kw_mirabolic_triv}
	V[\hbar] \xrightarrow{\sim} \U(\m_N) \otimes V \hamp N_+, \qquad v \mapsto \pP v \pP.
\end{equation}
The following is an immediate corollary of Lemma~\ref{lm:kirillov_projector_canonical_homomorphism} and Proposition~\ref{prop:kw_triv_tensor_prod}.

\begin{Proposition}\quad
\label{prop:mirabolic_polynomial}
	\begin{enumerate}\itemsep=0pt
		\item[$1.$] There is a natural isomorphism
		\[
			\C[\hbar][\Omega_1,\dots,\Omega_k]_{<N} \xrightarrow{\sim} \U(\m_N) \otimes \bigl(\C^N\bigr)^{\otimes k} \hamp N_+, \qquad f\bigl(\vec{\Omega}\bigr) \mapsto f\bigl(\vec{\Omega}\bigr)\bigl(v_1^{\otimes k}\bigr).
		\]
		\item[$2.$] The composition with the trivialization \eqref{eq::kw_mirabolic_triv}
		\[
			\varphi \colon \bigl(\C^N\bigr)^{\otimes k}[\hbar] \xrightarrow{\sim} \C[\hbar][\Omega_1,\dots,\Omega_k]_{<N}
		\]
		is uniquely characterized by the properties
		\begin{align*}
			&\varphi(v_{i_1} \otimes \dots \otimes v_{i_k})= \psi^{k-i_1-\dots-i_k} \Omega_1^{i_1-1} \cdots \Omega_k^{i_k-1}, \\
			&\varphi\bigl(v_{i_{\sigma(1)}} \otimes \dots \otimes v_{i_{\sigma(k)}}\bigr)= \sigma^{\hbar} \cdot \varphi(v_{i_1} \otimes \dots \otimes v_{i_k}),
		\end{align*}
		where the action is of that of the degenerate affine Hecke algebra.
	\end{enumerate}
\end{Proposition}

\begin{Example}
\label{example:mirabolic_symmetric_power}
	Consider $S^k \C^N \subset V^{\otimes k}$. It has a basis
	\[
		v_{i_1} \cdots v_{i_k} := \frac{1}{k!} \sum_{\sigma\in S_k} v_{i_{\sigma(1)}} \otimes \dots \otimes v_{i_{\sigma(k)}}, \qquad N \geq i_1 \geq \dots \geq i_k \geq 1.
	\]
	Therefore, by Propositions~\ref{prop:mirabolic_polynomial} and~\ref{prop:symmetrizer_hecke_vs_shuffle}, we have
	\[
		\pP v_{i_1} \cdots v_{i_k} \pP \mapsto \frac{\psi^{k-i_1-\dots-i_k}}{k!} \Sym_k \biggl( \prod_{i<j} \frac{\Omega_i - \Omega_j + \hbar}{\Omega_i - \Omega_j} \cdot \Omega_1^{i_1-1} \cdots \Omega_k^{i_k -1} \biggr),
	\]
	which is a rational analog of the Hall--Littlewood polynomial associated to the partition $({i_1-1},\allowbreak\dots,{i_k-1})$.
\end{Example}

\begin{Example}
%\label{example:mirabolic_exterior_power}
	Consider \smash{$\Lambda^k \C^N \subset \bigl(\C^N\bigr)^{\otimes k}$}. It has a basis
	\[
		v_{i_1} \wedge \dots \wedge v_{i_k} := \frac{1}{k!} \sum_{\sigma \in S_k} (-1)^{\sigma} v_{i_{\sigma(1)}} \otimes \dots \otimes v_{i_{\sigma(k)}}, \qquad N \geq i_1 > \dots > i_k \geq 1.
	\]
	Therefore, by Proposition~\ref{prop:mirabolic_polynomial} and Proposition~\ref{prop:symmetrizer_hecke_vs_shuffle}, we have
	\[
		\pP v_{i_1} \wedge \dots \wedge v_{i_k} \pP \mapsto \frac{\psi^{k-i_1-\dots-i_k}}{k!} \prod_{i<j} (\Omega_i - \Omega_j - \hbar)\cdot s_{\lambda}\bigl(\vec{\Omega}\bigr), \qquad \lambda = (i_1-k,\dots,i_k-1),
	\]
	where \smash{$s_{\lambda} \bigl(\vec{\Omega}\bigr)$} is the Schur polynomial associated to partition $\lambda$, see \cite[Section I.3]{Macdonald}, and the factor can be viewed as a deformation of the Vandermonde determinant.
\end{Example}

By \cite[Theorem~6.3]{KalmykovYangians}, the functor \eqref{eq::kw_mirabolic} is monoidal, in~particular, if $A$ is an algebra in $\Rep(\GL_N)$, one can define a mirabolic analog $\U(\m_N) \ltimes A$ of Definition~\ref{def:hc_algebra}, and its reduction $\U(\m_N) \ltimes A \hamp N_+$ is an algebra as well, similarly to Definition~\ref{def:kmw_algebra}. Recall the rational Feigin--Odesskii shuffle algebra ${}^N\!\cS^{\hbar}$ of Definition~\ref{def::fo_shuffle_algebra}. The following is a corollary of Theorem~\ref{thm:hecke_vs_fo_shuffle}, Theorem~\ref{thm:kw_to_yangian}, and Example~\ref{example:mirabolic_symmetric_power}.
\begin{Theorem}\label{thm:kw_to_shuffle_mirabolic}
	There is a $\C[\hbar]$-algebra isomorphism ${}^N\!\cS^{\hbar} \xrightarrow{\sim} \U(\m_N) \ltimes S^{\bullet} \C^N \hamp N_+$.
\end{Theorem}

\begin{Remark}
\label{remark:mirabolic_vs_stable_vector_bundle}
	As we mentioned in the introduction, the shuffle algebra ${}^N\!\cS^{\hbar}$ is a rational degeneration of its elliptic version, in view of Section~\ref{sect:deg_aha_yangians} and the results of \cite{EndelmanHodgesTwistedSU}. The original construction of \cite{FeiginOdesskiiSklyanin} of the latter used Belavin's $R$-matrix that can be associated to a stable bundle of rank $N$ and degree 1 (for instance, see ``Geometric interpretation'' after \cite[formula~(1.56)]{LOZ} in loc.\ cit.\ for precise statement). Unfortunately, we do not know an interpretation of Theorem~\ref{thm:kw_to_shuffle_mirabolic} in these terms. However, it seems to be a step in the right direction based on the following observation.
	
	Let $E$ be an elliptic curve with the origin $0\in E$. Denote by $\xi_{n,r}$ a unique stable bundle $E$ of rank $r$ whose determinant is $\cO_E (n\cdot 0)$. As it is argued in \cite{FeiginOdesskii}, the elliptic shuffle algebra quantizes a certain Poisson bracket on the space of extensions $\mathrm{Ext}^1(\xi_{N,1},\xi_{0,1})$. Its ring of functions is the symmetric algebra on
	\[
		\mathrm{Ext}^1(\xi_{N,1},\xi_{0,1})^* \cong \Hom^0 (\xi_{0,1}, \xi_{N,1}) \cong \Hom^0 (\xi_{-N,1}, \xi_{0,1}).
	\]
	Following loc.\ cit., we can apply the Fourier--Mukai transform
	\[
		\Phi(\xi_{-N,1}) = \xi_{1,N}[-1], \qquad \Phi(\xi_{0,1})[-1] = \cO_{0}[-1],
	\]
	where $\cO_0$ is the skyscraper sheaf at the origin. Therefore, it transforms $\Hom^0 (\xi_{-N,1}, \xi_{0,1})$ to the space of morphisms $\xi_{1,N} \rightarrow \cO_0$ which can be interpreted as a mirabolic reduction of the fiber of $\xi_{1,N}$ at the origin. To relate to the Kostant--Whittaker reduction, recall that it quantizes the Kostant slice \cite{Kostant}. The latter can be identified with $\h/W$. Its elliptic version is the coarse moduli space of degree 0 semistable bundles of rank $N$. At the same time, we can complete a~non-zero $\varphi \in \Hom^0 (\xi_{-N,1}, \xi_{0,1})$ to an exact triangle
	\[
		 \xi_{-N,1} \rightarrow \xi_{0,1} \rightarrow \cF \rightarrow \xi_{-N,1}[1],
	\]
	where $\cF$ is some torsion sheaf. Since $\Phi$ is an exact functor, we obtain
	\[
		\xi_{1,N}[-1] \rightarrow \cO_{0}[-1] \rightarrow \mathcal{E} \rightarrow \xi_{1,N},
	\]
	where $\mathcal{E} = \Phi(\cF)$ is a degree zero semistable bundle of rank $N$. It can be interpreted as a \emph{Hecke modification}
	\[
		0 \rightarrow \mathcal{E} \rightarrow \xi_{1,N} \rightarrow \cO_0 \rightarrow 0
	\]
	at $0\in E$. Therefore, it seems that $\U(\m_N) \otimes S^{\bullet} \C^N \hamp N_+$ should be a certain quantum (and rational) version of the algebra of Hecke modifications of special kind.
\end{Remark}

\section{Finite Miura transform}
\label{sect:miura}

In this section, we relate the Kostant--Whittaker reduction of Section~\ref{sect:kw_reduction} to the parabolic restriction of Section~\ref{sect:parabolic_reduction} following the construction of \cite[Section 6]{GinzburgKazhdan}, a finite analog of the Miura transform for W-algebras. We compare it with the GKLO homomorphism for the positive part $\rY_{-N,N}(\sl_2)^{\geq 0}$ from Section~\ref{sect:deg_aha_yangians} in Theorem~\ref{thm:miura_vs_gklo}. It will be extended to the whole Yangian in Section~\ref{sect:toda} and used to prove the main result Theorem~\ref{thm:toda_lattice_yangian}.

In this section, we relate the Kostant--Whittaker reduction of Section~\ref{sect:kw_reduction} to the parabolic restriction of Section~\ref{sect:parabolic_reduction} following the construction of \cite[Section 6]{GinzburgKazhdan}.

\begin{Definition}
	Let $X\in \HC(\GL_N)$. The {\it Miura bimodule} associated to $X$ is the double quotient $\n_- \backslash X^{\gen} /\n^{\psi}_+$.
\end{Definition}

The definition is motivated by the following simple fact: if $X$ is an algebra object in $\HC(\GL_N)$, then there is an action
\[
	X^{\gen} \sslash N_- \curvearrowright \n_- \backslash X^{\gen} / \n^{\psi}_+ \curvearrowleft X \hamp N_+
\]
of the corresponding reductions.

Recall the Harish-Chandra map of Example~\ref{ex:parabolic_res_center}:
\begin{equation}
\label{eq::hc_map_reminder}
	\Z_{\gl_N} \rightarrow \U(\gl_N)^{\gen} \sslash N_- \cong \U(\h)^{\gen}, \qquad A(u) \mapsto PA(u) = \prod_{i=1}^{N} (u-w_i).
\end{equation}
Following \cite[Section 6]{GinzburgKazhdan}, consider the map
\begin{equation}
\label{eq::kw_to_miura}
	\U(\h)^{\gen} \otimes_{\Z_{\gl_N}} (X\hamp N_+) \xrightarrow{1\otimes p_X} \U(\h)^{\gen} \otimes_{\Z_{\gl_N}} \n_- \backslash X / \n^{\psi}_+ \xrightarrow{\mathrm{act}} \n_- \backslash X^{\gen} / \n_+,
\end{equation}
for $X\in \HC(\GL_N)$, where
\[
	p_X \colon\ X \hamp N_+ = \bigl(X/\n_+^{\psi}\bigr)^{N_+} \rightarrow \n_- \backslash X / \n_+^{\psi}
\]
is the projection map and
\[
	\mathrm{act} \colon\ \U(\h)^{\gen} \otimes \n_- \backslash X \rightarrow \n_- \backslash X^{\gen}
\]
is the left action map. Observe that the latter factors through the tensor product over $\Z_{\gl_N}$ by~\eqref{eq::hc_map_reminder}.
Also, consider the quotient map
\begin{equation}
\label{eq::par_res_to_miura}
	X^{\gen} \sslash N_- = (\n_- \backslash X^{\gen}) ^{N_-} \rightarrow \n_- \backslash X^{\gen} /\n^{\psi}_+.
\end{equation}

\begin{Proposition}\quad
	\begin{enumerate}\itemsep=0pt
		\item[$1.$] The map \eqref{eq::kw_to_miura} is an isomorphism {\rm \cite[{\it Lemma} 6.2.1]{GinzburgKazhdan}}.
		\item[$2.$] The map \eqref{eq::par_res_to_miura} is an isomorphism {\rm \cite[{\it Proposition} 7.6]{KalmykovYangians}}.
	\end{enumerate}
\end{Proposition}

\begin{Remark}
	For free Harish-Chandra bimodules, the isomorphism of the second part holds without genericity assumption by \cite[Theorem~7.1.8]{GinzburgKazhdan}.
\end{Remark}

This motivates the following definition.
\begin{Definition}
\label{def:miura_transform}
	For a Harish-Chandra bimodule $X\in \HC(\GL_N)$, the {\it Miura transform} is the composition
	\[
		\miura_X\colon\ X \hamp N_+ \rightarrow \n_- \backslash X^{\gen} / \n_+^{\psi} \rightarrow X^{\gen} \sslash N_-.
	\]
\end{Definition}

\begin{Remark}
\label{remark:miura_w_action}
	In view of the map \eqref{eq::par_res_to_miura}, the Miura transform induces a $W$-action on $X^{\gen} \sslash N_-$ extending the one of Example~\ref{ex:parabolic_res_center} such that $X\hamp N_+ \subset (X^{\gen} \sslash N_-)^W$.
\end{Remark}

\begin{Lemma}
\label{lm:kw_to_par_res_algebra}
	If $X$ is an algebra object in $\HC(\GL_N)$, then the Miura transform is an algebra~map.
	\begin{proof}
		Observe that the map \eqref{eq::kw_to_miura} can be equivalently presented as the action
		\[
			X \hamp N_+ \xrightarrow{[1_X] \otimes \mathrm{id} }\bigl(\n_- \backslash X^{\gen} / \n^{\psi}_+\bigr) \otimes (X \hamp N_+) \rightarrow \bigl(\n_- \backslash X^{\gen} / \n^{\psi}_+\bigr)
		\]
		on the class of identity $1_X\in X$ in $\n_- \backslash X^{\gen} / \n^{\psi}_+$, and similarly for \eqref{eq::par_res_to_miura}. Then the lemma follows from a simple general result: let $A,B$ be two algebras and $M$ be an $(A,B)$-bimodule with an~element~$m\in M$ such that the composition
		\[
			A \xrightarrow{\mathrm{id}_A \otimes m} A \otimes M \xrightarrow{\mathrm{act}} M
		\]
		is an isomorphism. Then the composition with its inverse
		\[
			B \xrightarrow{m \otimes \mathrm{id}_B} M \otimes B \xrightarrow{\mathrm{act}} M \rightarrow A
		\]
		is an algebra map.
	\end{proof}
\end{Lemma}

Consider $\U(\gl_N) \ltimes S^{\bullet} \C^N$. Combining Theorems~\ref{thm:kw_to_yangian} and~\ref{thm:par_res_to_difference_ops}, and Lemma~\ref{lm:kw_to_par_res_algebra}, we obtain an algebra homomorphism
\begin{equation}
\label{eq::miura_yangian_positive}
	\rY_{-N,N}(\gl_2)^{\geq 0} \rightarrow \U(\gl_N) \ltimes S^{\bullet} \C^N \hamp N_+ \rightarrow \U(\gl_N)^{\gen} \ltimes S^{\bullet} \C^N \sslash N_- \cong \cA^-.
\end{equation}

\begin{Theorem}
\label{thm:miura_vs_gklo}
	The map \eqref{eq::miura_yangian_positive} coincides with the GKLO map of Theorem~$\ref{thm:gklo}$.
	\begin{proof}
		It follows from Example~\ref{ex:parabolic_res_center} that
		\[
			d_1(u) \mapsto \prod_{i=1}^{N}(u-w_i), \qquad d_2(u) \mapsto \prod_{i=1}^{N} (u-w_i - \hbar)^{-1}.
		\]
		By arguments of \cite[Corollary 7.12]{KalmykovYangians} adapted to the positive case, we have
		\begin{gather*}
			\U(\gl_N) \ltimes \C^N \hamp N_+ \ni \frac{1}{u-\Omega_1}(v_1)\\
\qquad{} \mapsto \sum_{i=1}^N \frac{1}{u-w_i} \prod_{j<i} (w_i - w_j)^{-1} \cdot Pv_i P\in \U(\gl_N)^{\gen} \otimes \C^N \sslash N_-.
		\end{gather*}
		By \eqref{eq::parabolic_res_dual_difference_gens}, we obtain
		\[
			\frac{1}{u-\Omega_1} (v_1) \mapsto \sum_{i=1}^N \frac{1}{u-w_i} \prod_{j\neq i} \frac{1}{w_i-w_j} \cdot \bar{v}_i,
		\]
		therefore, under the composition \eqref{eq::miura_yangian_positive}, we have
		\[
			x_+(u) \mapsto \sum_{i=1}^N \frac{1}{u-w_i} \prod_{j\neq i} \frac{1}{w_i-w_j} \cdot \bu_i^{-1},
		\]
		as required.
	\end{proof}
\end{Theorem}

\section{Toda lattice}
\label{sect:toda}

In this section, we show an isomorphism between the quantum Toda lattice for $\GL_N$ and the truncated shifted Yangian from Section~\ref{sect:deg_aha_yangians} in Theorem~\ref{thm:toda_lattice_yangian}. The result is based on identification of the Miura transform of the quantum Toda lattice with the algebra of difference operators, see Theorem~\ref{thm:diff_ops_vs_difference}, and comparing their images inside the latter. In Section~\ref{subsect:comparison}, we compare our result with its geometric counterpart, in terms of equivariant homology of the affine Grassmannian. Finally, in Section~\ref{subsect:monopole}, we show that the images of matrix coefficient functions associated to fundamental weights are given by degeneration of (rational) Macdonald operators.

The main object of this section is the algebra $\dD(\GL_N)$ of differential operators from Section~\ref{subsect:diff_ops}. Recall that it has two $\GL_N$ actions.

\begin{Notation}
	In what follows, we use the superscript $R$ (resp.\ $L$) to denote everything related to the right action (resp.\ left action) of $\GL_N$ on itself. For instance, let $\HC\bigl(\GL_N^R\bigr)$ be the Harish-Chandra category $\HC(\GL_N)$ such that $\dD(\GL_N)$ is an object $\HC\bigl(\GL_N^R\bigr)$ via:
	\begin{itemize}\itemsep=0pt
		\item The left $\U(\gl_N)$-action is induced from multiplication by the vector fields generating the right translations;
		\item The $\GL_N$-action is induced from the right action of $\GL_N$ on itself.
	\end{itemize}
	Similarly, we denote by $\HC\bigl(\GL_N^L\bigr)$ the Harish-Chandra category $\HC(\GL_N)$ such that the $\U(\gl_N)$-bimodule structure on $\dD(\GL_N)$ is induced from the left action. For any $V\!\in\! \Rep(\GL_N)$, define
	\[
		\U\bigl(\gl_N^R\bigr) \otimes V \in \HC\bigl(\GL_N^R\bigr), \qquad \U\bigl(\gl_N^L\bigr) \otimes V \in \HC\bigl(\GL_N^L\bigr)
	\]
	the corresponding free Harish-Chandra bimodules of Section~\ref{sect:hc_bimodules}.
	
	For the rest of the section, we fix $\psi=1$ in \eqref{eq::nilp_character}. For any $X\in \HC\bigl(\GL^L_N\bigr)$, denote by
	\[
		\lhamp X := (X/\n_+^{L,-\psi})^{N^L_+}
	\]
	the Kostant--Whittaker reduction of Section~\ref{sect:kw_reduction} for left Harish-Chandra bimodules. In particular,
	\[
		\lhamp \dD(\GL_N) \in \HC\bigl(\GL_N^R\bigr).
	\]
	As in Notation~\ref{notation:par_res_quotient_extremal_proj}, for $x\in X$, we denote by
	\[
		x P_{-\psi}^L \in X /\n_+^{L,-\psi}, \qquad P_{-\psi}^L x P_{-\psi}^L \in \lhamp X
	\]
	correspondingly the image of $x$ under the projection map \smash{$X \rightarrow X /\n_+^{L,-\psi}$} and the action of the Kirillov projector from \eqref{eq::kirillov_proj_properties} on the latter. Recall that if \smash{$x P_{-\psi}^L$} is $N^L_+$-invariant, then we have \smash{$P_{-\psi}^L x P_{-\psi}^L = x P_{-\psi}^L$}.
	
	Likewise, for any $X\in \HC\bigl(\GL^R_N\bigr)$, denote by
	\[
		X \sslash N_-^R := \bigl(\n_-^R \backslash X\bigr)^{N_-^R}, \qquad X \hamp N^R_+ := \bigl(X/\n_+^{R,\psi}\bigr)^{N_+^R}
	\]
	the parabolic restriction of Section~\ref{sect:parabolic_reduction} and the Kostant--Whittaker reduction respectively for right Harish-Chandra bimodules. In particular, both
	\[
		\dD(\GL_N) \sslash N^R_- , \qquad \dD(\GL_N) \hamp N^R_+
	\]
	are objects of $\HC\bigl(\GL_N^L\bigr)$. Similarly to the left action, we use notations
	\[
		x \pP^R \in X/ \n_+^{R,\psi}, \qquad \pP^R x \pP^R \in X \hamp N_+^R,
	\]
	as well their versions
	\[
		P^R x \in \n_-^R \backslash X, \qquad P^R x P^R \in X \sslash N_-^R
	\]
	for the extremal projector as in Notation~\ref{notation:par_res_quotient_extremal_proj}.
\end{Notation}

In these notations, the quantum Toda lattice is $\lhamp \dD(\GL_N) \hamp N^R_+$.
\begin{Remark}
	For simplicity, let $\hbar = 1$. Character $\psi \colon \n_+ \rightarrow \C$ gives a multiplicative character of $N_+$ that we denote by ${\rm e}^{\psi} \colon N_+ \rightarrow \C^*$. Let $U \subset \GL_N$ be the big Bruhat cell $U = N_+ w_0 H_N N_+$, where $w_0$ is the longest element in the Weyl group. Denote by
	\begin{equation}
	\label{eq::function_biinv}
		C^{\infty}(U) \supset C^{\infty}_{(N^L_+ \times N_+^R, -\psi \times \psi)} (U) := \bigl\{f(n_1 g n_2) = {\rm e}^{-\psi}(n_1) {\rm e}^{\psi}(n_2) f(g),\, n_1,n_2\in N_+\bigr\}
	\end{equation}
	the space of smooth functions on $U$ equivariant with respect to $N_+^L \times N_+^R$-action and character \smash{$\bigl({\rm e}^{-\psi}, {\rm e}^{\psi}\bigr)$}. It can be identified with $C^{\infty}(H_N)$, smooth functions on the torus.
	
	There is a well-defined action of $\lhamp \dD(\GL_N) \hamp N^R_+$ on \smash{$C\raisebox{1.5pt}{${}^{\infty}_{(N^L_+ \times N_+^R, -\psi \times \psi)}$} (U) \cong C^{\infty} (H_N)$}. The quantum Toda Hamiltonians can be obtained as the images of the embedding
	\[
		Z(\U(\gl_N)) \rightarrow \lhamp \dD(\GL_N) \hamp N^R_+
	\]
	as bi-invariant differential operators twisted by character ${\rm e}^{\rho_N} \colon H_N \rightarrow \C^*$ \eqref{eq::rho}. For instance, consider the operator
	\[
		C^R = \sum_{i=1}^N \bigl(E^R_{ii}\bigr)^2 + 2\sum_{i<j} E^R_{ji} E^R_{ij} + 2\bigl(\rho_N^{\vee}\bigr)^R = \sum_{i,j=1}^N E^R_{ij} E^R_{ji} + (N-1) \sum_{i=1}^N E^R_{ii},
	\]
	a slight modification of the Casimir operator, where
	\[
		\rho_N^{\vee} = \sum_{i=1}^N (N-i) E_{ii}.
	\]
	Using transformation properties \eqref{eq::function_biinv}, one can check that
	\[
		{\rm e}^{\rho_N} C^R {\rm e}^{-\rho_N} = \sum_{i=1}^N \partial_i^2 - 2\psi^2 \sum_{i=1}^{N-1} {\rm e}^{x_i - x_{i+1}} - (\rho_N,\rho_N),
	\]
	the standard Toda Hamiltonian. Here $\{{\rm e}^{x_i}\}$ are natural coordinates of diagonal entries on $H_N$.
\end{Remark}

The following is a direct corollary of definitions \eqref{eq::canonical_homomorphism} and \eqref{eq::canonical_homomorphism_dual}.
\begin{Proposition}
\label{prop:canonical_matrix_desciption}
	Under the embeddings of Example~$\ref{ex:embedding_diff_ops}$, we have
	\begin{alignat*}{3}
		&{}^R \iota_{\alpha} (\Omega(v_i)) = \bigl(E^R X\bigr)_{i\alpha}, \qquad && {}^R\iota_{\alpha}^*(\Omega^*(\phi_i)) = \bigl(X^{-1} E^R\bigr)_{\alpha i},& \\
		&{}^L \iota_{\alpha}^* (\Omega^*(\phi_i)) = \bigl(X E^L\bigr)_{\alpha i}, \qquad && {}^L\iota_{\alpha}(\Omega(v_i)) = \bigl(E^L X^{-1}\bigr)_{i \alpha}.&
	\end{alignat*}
\end{Proposition}

Since $v_1 \in \C^N$ \big(resp.\ $\phi_N \in \bigl(\C^N\bigr)^*$\big) is a highest-weight vector, the embeddings of Example~\ref{ex:embedding_diff_ops} factor through the maps of $\GL_N^R$ Harish-Chandra bimodules
\begin{align}
	{}^R\iota_N \otimes_{\U(\gl_N^R)} {}^R\iota^*_1 \colon\ \U\bigl(\gl^R_N\bigr) \ltimes \bigl(S^{\bullet} \C^N \otimes S^{\bullet} \bigl(\C^N\bigr)^*\bigr) &\rightarrow \lhamp \dD(\GL_N), \nonumber \\
	 v_i &\mapsto X_{iN}, \qquad \phi_i \mapsto \bigl(X^{-1}\bigr)_{1i}.\label{eq::matrix_units_embedding}
\end{align}

The maps \eqref{eq::matrix_units_embedding} descend to homomorphisms between corresponding reductions of Section~\ref{sect:parabolic_reduction} and Section~\ref{sect:kw_reduction}. Denote by
\begin{equation}
\label{eq::q_vars}
	Z_{iN} := {}^R\iota_N(\bar{v}_i), \qquad Z^*_{1i} := {}^R\iota^*_1\bigl(\bar{\phi}_i\bigr), \qquad Z_{iN},Z^*_{1j} \in \lhamp \dD^{\gen}(\GL_N) \sslash N^R_-,
\end{equation}
in the notations of \eqref{eq::parabolic_res_dual_difference_gens}, \eqref{eq::par_res_dual_vec}, and \eqref{eq::par_res_dual_vec}, where $\dD^{\gen}(\GL_N) := \U\bigl(\gl_N^R\bigr)^{\gen} \otimes_{\U(\gl_N^R)} \dD(\GL_N)$. Observe that the double reduction is an algebra.

\begin{Lemma}
	\label{lm:diff_ops_commute}
	For $i\neq j$, we have $Z^*_{1i} Z_{jN} = Z_{jN} Z^*_{1i}$.
	\begin{proof}
		Using \eqref{eq::matrix_units_embedding} and \eqref{eq::q_vars}, it is enough to prove that $\bar{\phi}_i \bar{v}_j = \bar{v}_j \bar{\phi}_i$. This is the content of Proposition~\ref{prop:diff_ops_commute}.
	\end{proof}
\end{Lemma}

\begin{Proposition}
\label{prop:diff_ops_inverses}
	For any $i$, we have $Z^*_{1i} Z_{iN} = 1 = Z_{iN} Z^*_{1i}$.
	\begin{proof}
		We prove only the first equality. Recall the quantum comatrix $\hat{T}^L(u)$ from Section~\ref{subsect:diff_ops} associated to vector fields generating left translations. Since they commute with the ones generating right translations, its matrix entries descend to elements of $\dD(\GL_N) \sslash N_-^R$, in other words, we have equality
		\[
			P^R \hat{T}^L_{1N}(u) = P^R\hat{T}^L_{1N}(u) P^R \in \dD(\GL_N) \sslash N_-^R.
		\]
		Using Propositions~\ref{prop:comatrix_left_right} and~\ref{prop:comatrix_interpolation}, we get
		\[
			P^R \hat{T}^L_{1N}(u) P_{-\psi}^L= \sum_{a=1}^N \prod_{b\neq a} \frac{u-w^R_b}{w^R_a-w^R_b} \cdot Z^*_{1a} Z_{aN} \in \lhamp \dD(\GL_N) \sslash N^R_-.
		\]
		In particular, we have $Z_{1a}^* Z_{aN} = P^R \hat{T}^L_{1N}\bigl(w_a^R\bigr) P_{-\psi}^L$. At the same time, one can show analogously to \cite[Proposition 4.8]{KalmykovYangians} that
		\[
			\hat{T}_{1N}^L(u) P_{-\psi}^L = (-1)^{N+1} T_{2\dots N}^{1\dots N-1}(u) P_{-\psi}^L = \psi^{N+1} P_{-\psi}^L \in \U\bigl(\gl_N^L\bigr) /\n_+^{L,-\psi}.
		\]
		Since $\psi=1$, we conclude that $Z^*_{1a} Z_{aN} = 1$ in $\lhamp \dD^{\gen}(\GL_N) \sslash N_-^R$, as required.
	\end{proof}
\end{Proposition}

The center $\Z_{\gl_N}$ is embedded in $\dD(\GL_N)$ as a subalgebra of bi-invariant differential operators, in particular, descends to the quantum Hamiltonian reduction $\lhamp \dD(\GL_N) \hamp N_+^R$. For any matrix $M$, we denote by $[M]_{ij}$ its $(i,j)$-matrix entry.
\begin{Lemma}\quad
\label{lm:generators}
	\begin{enumerate}\itemsep=0pt
		\item[$1.$] The algebra $\lhamp \dD(\GL_N) \hamp N_+^R$ is generated over $\Z_{\gl_N}$ by \smash{$\bigl[\bigl(E^R\bigr)^k X\bigr]_{1N} P_{\psi}^R$} and \linebreak \smash{$\bigl[X^{-1} \bigl(E^R\bigr)^k \bigr]_{1N} P_{\psi}^R$} for every $0 \leq k \leq N-1$.
		
		\item[$2.$] The algebra $\lhamp \dD^{\gen} \sslash N_+^R$ is generated by $Z_{iN}$ and $Z^*_{1i}$ from \eqref{eq::q_vars} over \smash{$\U\bigl(\h^R\bigr)^{\gen}$} for $1\leq i \leq N$.
	\end{enumerate}
	\begin{proof}
		By \eqref{eq::kw_trivialization}, we have an isomorphism
		\begin{align*}
		\Z_{\gl_N} \otimes \cO(\GL_N) &\xrightarrow{\sim} \lhamp \U\bigl(\gl_N^L\bigr) \ltimes \cO(\GL_N) \cong \lhamp \dD(\GL_N), \\
 f &\mapsto P_{-\psi}^L f P_{-\psi}^L.
		\end{align*}
		In particular, we see that $\lhamp \dD(\GL_N)$ is generated over $\Z_{\gl_N}$ by the classes \smash{$P_{-\psi}^L X_{ij} P_{-\psi}^L$} and \smash{$P_{-\psi}^L \bigl(X^{-1}\bigr)_{ij} P_{-\psi}^L$} for all $i$, $j$. By Proposition~\ref{prop:kw_tensor_algebra} and its dual version of Proposition~\ref{prop:kw_tensor_algebra_dual}, it is equivalently generated by the classes of matrix entries \smash{$\bigl[X \bigl(E^L\bigr)^k\bigr]_{\alpha N} P_{-\psi}^L$} and \smash{$\bigl[\bigl(E^L\bigr)^k X^{-1}\bigr]_{1\alpha} P_{-\psi}^L$} in view of Proposition~\ref{prop:canonical_matrix_desciption}, for every $\alpha$ and $k$. Therefore, by \eqref{eq::right_vs_left_matrix}, $\lhamp \dD(\GL_N)$ is generated by matrix entries \smash{$\bigl[\bigl(E^R\bigr)^kX\bigr]_{\alpha N} $} and \smash{$\bigl[X^{-1} \bigl(E^R\bigr)^k\bigr]_{1\alpha}$}; for readability, we drop the sign $P_{-\psi}^{L}$.
		
		Let us now consider the right action. For the first part of the lemma, the last statement implies that the double reduction \smash{$\lhamp \dD(\GL_N) \hamp N_+^R$} is generated by
		\[
			\pP^R \bigl[\bigl(E^R\bigr)^kX\bigr]_{\alpha N} \pP^R, \qquad \pP^R \bigl[X^{-1} \bigl(E^R\bigr)^k\bigr]_{1\alpha} \pP^R.
		\]
		By definition, we have
		\[
			\pP^R \bigl[\bigl(E^R\bigr)^kX\bigr]_{\alpha N} \pP^R = {}^R \iota_{N} \bigl(\pP \Omega^{k} (v_i) \pP\bigr), \qquad \pP^R \bigl[X^{-1} \bigl(E^R\bigr)^k\bigr]_{1\alpha} \pP^R ={}^R \iota_1^* \bigl(\pP (\Omega^*)^k (\phi_i) \bigr).
		\]
		Therefore, by Lemma~\ref{lm:kirillov_projector_canonical_homomorphism}, its dual version of Proposition~\ref{prop:dual_kw_presentation}, and Proposition~\ref{prop:cayley_hamilton}, the double Hamiltonian reduction is generated by the $N_+$-invariant classes
		\[
			{}^R \iota_N \bigl(\Omega^k (v_1) \pP\bigr) = \bigl[\bigl(E^R\bigr)^k X\bigr]_{1N} \pP^R, \qquad {}^R \iota_1^* \bigl((\Omega^*)^k (\phi_N) \pP\bigr) = \bigl[X^{-1} \bigl(E^R\bigr)^k\bigr]_{1N} \pP^R
		\]
		for $0\leq k \leq N-1$, as required.
	
		For the second part, it follows that the double reduction \smash{$\lhamp \dD^{\gen} (\GL_N) \sslash N_-^R$} is generated over \smash{$\U\bigl(\h^R\bigr)^{\gen}$} by
		\[
			P^R \bigl[\bigl(E^R\bigr)^kX\bigr]_{\alpha N} P^R, \qquad P^R \bigl[X^{-1} \bigl(E^R\bigr)^k\bigr]_{1\alpha} P^R.
		\]
		By Lemma~\ref{lm:parabolic_res_canonical_homomorphism} and its dual version \eqref{eq::omega_par_res_dual}, it is equivalently generated by \eqref{eq::q_vars}, as required.
	\end{proof}
\end{Lemma}

Recall the algebra of difference operators $\cA$ from Section~\ref{subsect:mickelsson}. We will need the following result which is an analog of simplicity of Weyl algebras.
\begin{Lemma}
\label{lm:difference_ops_simple}
	Any two-sided ideal in $\cA$ is generated by an ideal in $\C[\hbar]$.
	\begin{proof}
		It is enough to prove the lemma for the polynomial version $\U(\h)\bigl[\bu_1^{\pm},\dots,\bu_N^{\pm}\bigr]$. Assume there is an ideal $I$ inside the latter and $a\in I$ is some element. We can write
		\[
			a = \sum_{J=(j_1,\dots,j_N)} c_J(w) \bu_{1}^{j_1} \cdots \bu_{N}^{j_N}, \qquad c_J(w) \in \U(\h).
		\]
		We can assume that all the powers $\{\bu_i\}$ are positive. By repeatedly using
		\[
			\bu^{j_k}_k w_l = (w_l + \delta_{kl} \hbar j_k) \bu^{j_k},
		\]
		we can assume that there is only one monomial in $a$. Multiplying by its inverse, we obtain that~$I$ contains some polynomial $c(w)\in \U(\h)$. We have
		\[
			[\bu_i, c(w_1,\dots,w_N)] \bu_i^{-1} = c(w_1,\dots,w_i + \hbar,\dots,w_N) - c(w_1,\dots,w_i,\dots,w_N).
		\]
		Observe that the right-hand side is zero only when it is constant in $w_i$, at the same time, it has strictly lower degree in $w_i$ than $c(w)$. By repeatedly using this modification, we can assume that~$c(w)$ is $\C[\hbar]$-constant, as required.
	\end{proof}
\end{Lemma}

\begin{Theorem}
\label{thm:diff_ops_vs_difference}
	The $\U\bigl(\h^R\bigr)^{\gen}$-linear map
	\[
		\cA \xrightarrow{\sim} \lhamp \dD^{\gen}(\GL_N) \sslash N^R_-, \qquad \bu_i \mapsto Z^*_{1i}, \qquad \bu^{-1}_i \mapsto Z_{iN},
	\]
	is an isomorphism.
	\begin{proof}
		It is an algebra homomorphism by Lemma~\ref{lm:diff_ops_commute} and Proposition~\ref{prop:diff_ops_inverses}. Surjectivity follows from the second part of Lemma~\ref{lm:generators}. If it had a kernel, then it would be generated by some ideal in~$\C[\hbar]$ by Lemma~\ref{lm:difference_ops_simple}, in particular, the action of $\C[\hbar]$ on the class of \smash{$[1]\!\in \!\lhamp \dD(\GL_N)\sslash N_-^R$} would have a kernel as well. Consider the Bruhat cell $U = B_+ N_- \subset \GL_N$, where $B_+$ is the Borel subgroup. Since $\bigl(N_+^L,N^R_-\bigr)$-action is free, we have
		\[
			\lhamp \dD(\GL_N) \sslash N_-^R \hookrightarrow \lhamp \dD(U) \sslash N_-^R \cong \dD(H_N),
		\]
		where $\dD(H_N)$ is the differential operators on the torus. The action of $\C[\hbar]$ on $1\in \dD(H_N)$ is free, therefore, the map \smash{$\cA \rightarrow \lhamp \dD(\GL_N) \sslash N_-^R$} has no kernel and is an isomorphism.
	\end{proof}
\end{Theorem}

Finally, recall the finite Miura transform of Definition~\ref{def:miura_transform}:
\begin{equation}
\label{eq::miura_diff_ops}
	\miura_N \colon\ \lhamp \dD(\GL_N) \hamp N^R_+ \rightarrow \lhamp \dD^{\gen}(\GL_N) \sslash N^R_-.
\end{equation}

\begin{Lemma}
\label{lm:diff_ops_miura_embedding}
	The map $\miura_N$ is an embedding.
	\begin{proof}
		Consider the Miura transform
		\[
			 \dD(\GL_N) \hamp N^R_+ \rightarrow \dD^{\gen}(\GL_N) \sslash N^R_-.
		\]
		Since $\dD(\GL_N) \cong \U\bigl(\gl_N^R\bigr) \otimes \cO(\GL_N)$ is a free module, it follows from trivializations \eqref{eq::res_trivialization} and~\eqref{eq::kw_trivialization} that this map is injective. By exactness of the Kostant--Whittaker reduction \cite[Lemma~4\,(a)]{BezrukavnikovFinkelberg}, the induced map
		\[
			\lhamp \dD(\GL_N) \hamp N^R_+ \rightarrow \lhamp \dD^{\gen}(\GL_N) \sslash N^R_-
		\]
		is injective as well.
	\end{proof}
\end{Lemma}

The following theorem is an explicit presentation of the isomorphism from \cite[Appendix B]{BravermanFinkelbergNakajima}.

\begin{Theorem}
\label{thm:toda_lattice_yangian}
	There is an algebra isomorphism
	\[
		\rY^0_{-2N}(\sl_2) \xrightarrow{\sim} \lhamp \dD(\GL_N) \hamp N^R_+,
	\]
	defined by
	\begin{alignat*}{3}
		&d_1(u)\mapsto A^R(u), \qquad && d_2(u) \mapsto A^R(u-\hbar)^{-1},& \\
		&x_+(u)\mapsto -\bigl[\bigl(u-E^R\bigr)^{-1} X\bigr]_{1N}, \qquad && x_-(u)\mapsto \bigl[X^{-1}\bigl(u-E^R\bigr)^{-1}\bigr]_{1N},&
	\end{alignat*}
	such that the diagram
	\[
		\xymatrix{
			\rY^0_{-2N}(\sl_2) \ar[r]^-{\sim} \ar[d]_{\gklo_N} & \lhamp \dD(\GL_N) \hamp N^R_+ \ar[d]^{\miura_N} \\
			\cA \ar[r]^-{\sim} & \lhamp \dD^{\gen}(\GL_N) \sslash N^R_-
		}
	\]
	is commutative.
	\begin{proof}
		Consider the composition
		\[
			\rY_{-N,N}(\gl_2)^{\geq 0} \rightarrow \U\bigl(\gl^R_N\bigr) \otimes S^{\bullet} \C^N \hamp N^R_+ \xrightarrow{{}^R \iota_{N}} \lhamp \dD(\GL_N) \hamp N^R_+,
		\]
		where the first map is \eqref{eq::miura_yangian_positive}. By naturality of the Miura transform and Theorem~\ref{thm:kw_to_yangian}, the square
		\[
		\xymatrix{
			\rY_{-N,N}(\gl_2)^{\geq 0} \ar[r] \ar[d]_{\gklo_N} & \lhamp \dD(\GL_N) \hamp N^R_+ \ar[d]^{\miura_N} \\
			\cA^- \ar[r] & \lhamp \dD^{\gen}(\GL_N) \sslash N^R_-
		}
		\]
		is commutative. Likewise, the map ${}^R \iota^*_1$ gives a commutative square
		\[
		\xymatrix{
			\rY_{-N,N}(\gl_2)^{\leq 0} \ar[r] \ar[d]_{\gklo_N} & \lhamp \dD(\GL_N) \hamp N^R_+ \ar[d]^{\miura_N} \\
			\cA^+ \ar[r] & \lhamp \dD^{\gen}(\GL_N) \sslash N^R_-
		}
		\]
		by Theorem~\ref{thm:kw_to_yangian_dual}. To check the remaining relation \eqref{eq::yangian_e_f_comm}, we can use Lemma~\ref{lm:diff_ops_miura_embedding} and check it inside
		\[
			\cA \cong \lhamp \dD^{\gen}(\GL_N) \sslash N^R_-,
		\]
		which is a straightforward calculation, for instance, see \cite[Appendix B\,(v)]{BravermanFinkelbergNakajima}. Therefore, we have a~map
		\[
			\rY_{-N,N}(\gl_2) \rightarrow \lhamp \dD^{\gen}(\GL_N) \hamp N_+^R,
		\]
		whose composition with the Miura transform \eqref{eq::miura_diff_ops} coincides with the GKLO map. It is surjective by the first part of Lemma~\ref{lm:generators}. By Lemma~\ref{lm:diff_ops_miura_embedding}, the target is identified with its image in $\cA$. Then Definition~\ref{def:truncated_shifted_yangian} gives an isomorphism
		\[
			\rY^0_{-2N}(\sl_2) \xrightarrow{\sim} \lhamp \dD^{\gen}(\GL_N) \hamp N_+^R,
		\]
		as required.
	\end{proof}
\end{Theorem}

\begin{Remark}
\label{remark:w_action_difference_operators}
	In view of Remark~\ref{remark:miura_w_action}, the $W$-action on $\lhamp \dD(\GL_N) \sslash N_- \cong \cA$ is
	\[
	s \cdot w_i = w_{s(i)}, \qquad s\cdot \bu_i^{\pm} = \bu_{s(i)}^{\pm},
	\]
	which follows directly from $W$-invariance of the images of \smash{$x_{\pm}^{(1)}$} in \smash{$\rY^0_{-2N}(\sl_2)$}.
\end{Remark}

\subsection{Comparison with geometry}
\label{subsect:comparison}

Recall that the isomorphism \smash{$\rY^0_{-2N}(\sl_2) \cong \lhamp \dD(\GL_N) \hamp N_+^R$} is constructed in \cite{BravermanFinkelbergNakajima} by identifying the source with the homology of the affine Grassmannian and by using the derived Satake equivalence of \cite{BezrukavnikovFinkelberg}. In this subsection, we compare the isomorphism of loc.\ cit.\ with the one constructed in this paper.

On the topological side, consider the affine Grassmannian $\Gr = \GL_N((t))/ \GL_N[[t]]$. For any dominant coweight $\lambda \in \Lambda^+$, denote by $\Gr^{\lambda} =G[[t]] t^{\lambda} G[[t]]$ the corresponding orbit and by~$\IC^{\lambda}$ the IC sheaf on the closure \smash{$\overline{\Gr^{\lambda}}$}. The cohomology \smash{\raisebox{-0.5pt}{$H^{\bullet}_{\GL_N[[t]] \rtimes \C^*} \bigl(\overline{\Gr^{\lambda}},\IC^{\lambda}\bigr)$}} is a $H^{\bullet}_{\GL_N \times \C^*} (\mathrm{pt})$-bimodule via convolution.

On the algebraic side, consider the character
\[
	\psi^* \colon \n_- \rightarrow \C, \qquad \psi^*(E_{ij}) = -\delta_{j+1,i}.
\]
Analogously to Section~\ref{sect:kw_reduction}, one can define the negative version of the Kostant--Whittaker reduction
\[
	X \mapsto X \sslash_{\psi^*} N_- := \bigl(X/\n_-^{\psi^*}\bigr)^{N_-}, \qquad \n_-^{\psi^*} = \{\xi - \psi^*(\xi) \mid \xi \in \n_-\},
\]
for any Harish-Chandra bimodule $X$. Recall that the coefficients of the quantum minor \smash{$T^{1\dots N}_{1\dots N}(u)$} from \eqref{eq::quantum_determinant} generate the center $\Z_{\gl_N}$. Let us identify it with the generating function of the equivariant Chern classes of the \emph{dual} representation, see Remark~\ref{remark:center_identification},
\begin{equation}
\label{eq::center_topological}
	\Z_{\gl_N} \cong H^{\bullet}_{G}(\mathrm{pt}), \qquad T^{1\dots N}_{1\dots N}(u) \mapsto \sum_{i=0}^N c_i^{\GL_N} \bigl[\bigl(\C^N\bigr)^*\bigr] u^{N-i}.
\end{equation}
For every dominant weight $\lambda \in \Lambda$, the derived Satake equivalence of \cite{BezrukavnikovFinkelberg} gives an isomorphism of $\Z_{\gl_N}$-bimodules
\begin{equation*}
%\label{eq::derived_satake_iso}
	\U(\gl_N) \otimes V_{\lambda} \sslash_{\psi^*} N_- \xrightarrow{\sim} H^{\bullet}_{\GL_N[[t]] \rtimes \C^*} \bigl(\overline{\Gr^{\lambda}},\IC^{\lambda}\bigr) ,
\end{equation*}
as well as
\begin{equation}
\label{eq::toda_grassmannian}
	\lhamn \dD(\GL_N) \sslash_{\psi^*} N_-^R \xrightarrow{\sim} H^{\GL_N[[t]] \rtimes \C^*}_{\bullet} (\Gr),
\end{equation}
such that the diagram
\begin{equation*}
	\xymatrix{
		\U\bigl(\gl_N^R\bigr) \otimes V_{\lambda} \hamn N^R_- \ar[r] \ar[d] & H^{\bullet}_{\GL_N[[t]] \rtimes \C^*} \bigl(\overline{\Gr^{\lambda}},\IC^{\lambda}\bigr) \ar[d] \\
		\lhamn \dD(\GL_N) \sslash_{\psi^*} N_-^R \ar[r] & H^{\GL_N[[t]] \rtimes \C^*}_{\bullet} (\Gr)
	}
\end{equation*}
is commutative. Here the left vertical map is induced from
\[
	V_{\lambda} \rightarrow \cO(\GL_N), \qquad v \mapsto \langle g\cdot v_{\lambda}, v_{\lambda}^*\rangle,
\]
where $v_{\lambda}^*$ is the lowest-weight vector in $V_{\lambda}^*$.

\begin{Remark}
	Here and in what follows, we do \emph{not} keep track of the cohomological grading. Otherwise there should be a shift by the dimension of orbit.
\end{Remark}

To relate it to the constructions of this paper, consider the group automorphism
\[
	\GL_N \rightarrow \GL_N, \qquad g \mapsto \bigl(g^{-1}\bigr)^{\mathsf T}.
\]
One can easily check that it induces an algebra isomorphism
\begin{equation}
\label{eq::diff_ops_conj_inverse}
	\lhamp \dD(\GL_N) \hamp N_+^R \xrightarrow{\sim} \lhamn \dD(\GL_N) \sslash_{\psi^*} N_-^R.
\end{equation}
In particular, by Lemma~\ref{lm:quantum_determinant_transposition} and \eqref{eq::center_topological}, we have
\begin{equation}
\label{eq::center_topological_true}
	\Z_{\gl_N}[u] \ni A(u) \mapsto \sum_{i=0}^N c_i^{\GL_N} \bigl[\bigl(\C^N\bigr)^*\bigr] u^{N-i},
\end{equation}
where $A(u)$ is the series \eqref{eq::quantum_determinant}. For brevity, we will thus denote $A_i = c_i^{\GL_N}\bigl[\bigl(\C^N\bigr)^*\bigr]$.
Recall description of $\Gr^{\lambda}$ in terms of lattices, i.e., free $\C[[t]]$-modules of rank $N$ inside $\C((t))^N$~\cite{Zhu}. We will only consider a simple case $\lambda = \omega_1^*$ (and $\lambda = \omega_1$ in Appendix~\ref{subsect:dual_rep}), where $\omega_1^* = (0,\dots,0,-1)$. Then
\[
	\Gr^{\omega_1^*} = \overline{\Gr^{\omega_1^*}} = \bigl\{L_0 \subset L_1 \subset z^{-1} L_0 \mid \dim L_{1}/L_{0} = 1\bigr\} \cong \bP^{N-1} = \bP\bigl(\C^N\bigr),
\]
where \smash{$L_0 = \C[[t]]^N$} and identification on the right is $G$-equivariant. In particular, we have \smash{$\IC^{\omega_1^*} = \underline{\C}_{\Gr^{\omega_1^*}}$}. Using an isomorphism
\begin{equation}
\label{eq::varpi_hc}
	\tilde{\varpi}_N \colon\ \U\bigl(\gl_N^R\bigr) \ltimes S^{\bullet} \C^N \xrightarrow{\sim} \U\bigl(\gl^R_N\bigr) \ltimes S^{\bullet} \bigl(\C^N\bigr)^*, \qquad E_{ij}^R \mapsto -E_{ji}^R, \qquad v_i \mapsto \phi_i,
\end{equation}
one can check that there is a factorization
\begin{gather}
\label{eq::factorization}
\begin{split}
&	\xymatrix @C=1pc {
		\U\bigl(\gl_N^R\bigr) \otimes \C^N \hamp N_+^R \ar[r] \ar[d]_{\tilde{\varpi}_N} & \lhamp \dD(\GL_N) \hamp N_+^R \ar[r] & \lhamn \dD(\GL_N) \sslash_{\psi^*} N_-^R \ar[d] \\
		\U\bigl(\gl_N^R\bigr) \otimes \bigl(\C^N\bigr)^* \hamn N_-^R \ar[r]^{\sim} & H^{\bullet}_{\GL_N[[t]] \rtimes \C^*} \bigl(\Gr^{\omega_1^*},\underline{\C}_{\Gr^{\omega_1^*}}\bigr) \ar[r] & H^{\GL_N[[t]] \rtimes \C^*}_{\bullet} (\Gr).
	}\end{split}\hspace*{-10mm}
\end{gather}

Consider the bundle
\begin{equation}
\label{eq::tautological_bundles}
	\cV_i \rightarrow \Gr^{\omega_1^*}, \qquad \cV_i |_{L_0 \subset L_1} = z^{-1} L_i / L_i
\end{equation}
and the corresponding total Chern classes
\begin{equation}
\label{eq::total_chern_class}
	A^{(i)}(u) = \sum_{j=0}^N u^{N-j} A^{(i)}_j, \qquad A^{(i)}_j = (-1)^j c_j^{\GL_N} (\cV_{i-1}),
\end{equation}
for $i=1,2$. For instance, under identification $\Gr^{\omega_1^*} \cong \bP^{N-1}$, we have $A^{(1)}(u) = A(u)$ from \eqref{eq::center_topological_true}. The left (resp.\ right) action of $H^{\bullet}_{G}(\mathrm{pt})$ on \smash{$H^{\bullet}_{\GL_N[[t]] \rtimes \C^*} \bigl(\Gr^{\omega_1^*}\bigr)$} is via the coefficients of $A^{(1)}(u)$ \big(resp.\ $A^{(2)}(u)$\big). Also, consider the line bundle
\[
	\cS \rightarrow \Gr^{\omega_1^*}, \qquad \cS|_{L_0 \subset L_1} = L_1/L_0.
\]
Via identification \smash{$\Gr^{\omega_1^*} \cong \bP^{N-1}$}, we have $\cS \cong \cO_{\bP^{N-1}}(-1)$, and so
\begin{align*}
	H^{\bullet}_{\GL_N[[t]] \rtimes \C^*} \bigl(\Gr^{\omega_1^*},\underline{\C}_{\Gr^{\omega_1^*}}\bigr) &{}= H^{\bullet}_{\GL_N \times \C^*} \bigl(\bP^{N-1}\bigr)\\
 &{}\cong \C[\hbar][A_1,\dots,A_N]\bigl[c^{\GL_N}_1(\cS)\bigr]/\bigl(A\bigl(c^{\GL_N}_1(\cS)\bigr)\bigr).
\end{align*}
Recall that $G = \GL_N \times \C^*$, where $\C^*$ acts by loop rotations. In particular, we have
\begin{equation}
\label{eq::chern_class_c_star}
	c_1^{G}(\cS) = c_1^{\GL_N}(\cS) + \hbar.
\end{equation}
For the following proposition, we need an auxiliary result.
\begin{Lemma}
\label{lm:hc_unique_auto}
	Let $\phi \colon \U(\gl_N) \otimes \C^N \hamp N_+ \rightarrow \U(\gl_N) \otimes \C^N \hamp N_+$ be an automorphism of $\Z_{\gl_N}$-bimodules. Then $\phi = c\cdot \Id$, where $c\in \C$.
	\begin{proof}
		By \cite[Lemma~4]{BezrukavnikovFinkelberg}, the Kostant--Whittaker reduction functor is a full embedding when restricted to free Harish-Chandra bimodules. Therefore, we can lift $\phi$ to an automorphism of $\U(\gl_N) \otimes \C^N$ which induces an automorphism of the parabolic restriction $\U(\gl_N)^{\gen} \otimes \C^N \sslash N_-$ as of $\U(\h)^{\gen}$-bimodules. According to the decomposition \eqref{eq::hc_torus_decomposition}, it is isomorphic to a direct sum of $\U(\h)^{\gen}$-bimodules
		\[
			\U(\gl_N)^{\gen} \otimes \C^N \sslash N_- \cong \bigoplus_{i=1}^N \U(\h)^{\gen} \cdot \bar{v}_i,
		\]
		where $\bar{v}_i$ are vectors \eqref{eq::parabolic_res_dual_difference_gens}. Since $\phi$ was defined on the \emph{non-generic} bimodule $\U(\gl_N) \otimes \C^N$, the automorphism acts as $\bar{v}_i \mapsto f_i(w) \bar{v}_i$ for some \emph{polynomial} function $f_i(w) \in \U(\h)$ with a~\emph{polynomial} inverse, therefore, it should be constant $f_i(w) = c_i$. Moreover, since $\phi$ commutes with the diagonal action of the Weyl group $W = S_N \subset \GL_N$, we have $c_i = c_j$ for any $i$, $j$. Hence $\phi = c\cdot \Id$ by Lemma~\ref{lm:parabolic_res_injective}.
	\end{proof}
\end{Lemma}

\begin{Proposition}
\label{prop:top_vs_alg_vector}
	The $\Z_{\gl_N}$-bimodule isomorphism of \eqref{eq::factorization}
	\[
		\U\bigl(\gl_N^R\bigr) \otimes \C^N \hamp N_+^R \cong \Z_{\gl_N}[\Omega]/(A(\Omega)) \xrightarrow{\sim} \C[\hbar][A_1,\dots,A_N] \bigl[c^G_1(\cS)\bigr]/\bigl(A\bigl(c^G_1(\cS)\bigr)\bigr)
	\]
	 composed with Proposition~$\ref{prop:kw_vector}$ is given by \smash{$\Omega^k \mapsto (-1)^{N-1} c^{\GL_N}_1(\cS)^k$} for all $k$.
	\begin{proof}
		The left $\Z_{\gl_N}$-module structures are compatible by \eqref{eq::center_topological_true}. Let us show that this map intertwines the right actions. Following the proof of \cite[Lemma~8.13]{CautisKamnitzer}, consider two short exact sequences on $\Gr^{\omega_1^*}$:
		\begin{align*}
			& 0 \rightarrow L_1 / L_{0} \rightarrow z^{-1} L_{0} /L_{0} \rightarrow z^{-1} L_{0} / L_1 \rightarrow 0, \\
			&0 \rightarrow z^{-1} L_{0} / L_1 \rightarrow z^{-1} L_1 / L_1 \rightarrow z^{-1} L_1 / z^{-1} L_{0} \rightarrow 0,
		\end{align*}
		where each $L_i$ is considered as a vector bundle with the corresponding fiber. Recall that $\C^*$ acts by loop rotations on variable $z$. By multiplicativity of the total Chern class, we obtain
		\begin{align*}
			c^G(\cV_0) &{}= \bigl(u+c^G_1(\cS)\bigr) c^G\bigl(z^{-1} L_0/L_1\bigr) = \bigl(u+c^G_1(\cS)\bigr) \frac{c^G(\cV_1)}{c^G\bigl(z^{-1} L_1 / z^{-1} L_0\bigr)}\\
&{} = \bigl(u+c^G_1(\cS)\bigr) \frac{c^G(\cV_1)}{u+c^G_1(\cS) - \hbar},
		\end{align*}
		where \smash{$c^G(V) = \sum_i u^{\mathrm{rk}(V) - i} c_i(V)$} is the total Chern class of a vector bundle $V$ with respect to the action of $G= \GL_N \times \C^*$. Due to the additional $\C^*$-equivariance of $\cS$ and $\cV_i$, we have
		\[
			c_1^G(\cS) = c_1^{\GL_N}(\cS) - \hbar, \qquad c^G(\cV_i) = (-1)^N A^{(i+1)}(-u+\hbar)
		\]
		by \eqref{eq::total_chern_class}. Therefore, it implies
		\[
			A^{(2)}(u) = A^{(1)}(u) \frac{u - c_1(\cS) + \hbar}{u- c_1(\cS)},
		\]
		which agrees with the formula of Proposition~\ref{prop:conjugation_center}.
Therefore, this map is an isomorphism of $\Z_{\gl_N}$-bimodules such that \smash{$v_1 \mapsto (-1)^{N-1} \in H^{\bullet}_G\bigl(\bP^{N-1}\bigr)$}. The equivalence of Bezrukavnikov--Finkelberg gives an a priori different isomorphism with the same normalization condition by \cite[Lemma~2.13]{FKPRW}. However, it should coincide with the former by Lemma~\ref{lm:hc_unique_auto}.
	\end{proof}
\end{Proposition}

Finally, recall the truncated shifted Yangian of Definition~\ref{def:truncated_shifted_yangian}. By \cite[Theorem~B.18]{BravermanFinkelbergNakajima}, there is a surjective homomorphism
\begin{equation*}
	\rY^{0}_{-2N}(\sl_2) \rightarrow H_{\bullet}^{\GL_N[[t]] \rtimes \C^*} (\Gr),
\end{equation*}
such that
\begin{align}
		&d_1(u)\mapsto \sum_{p=0}^N u^{N-p} c_p^{\GL_N}\bigl[\bigl(\C^N\bigr)^*\bigr], \nonumber\\
		&x_-(u) \mapsto \frac{1}{u - c^G_1(\cQ) -\hbar} \cap [\Gr^{\omega_1}], \nonumber\\
		&x_+(u) \mapsto (-1)^N \frac{1}{u - c^G_1(\cS) -\hbar} \cap \bigl[\Gr^{\omega_1^*}\bigr],\label{eq::shifted_yangian_topological}
\end{align}
Here \smash{$[\Gr^{\omega_1}] \in H_{\bullet}^{\GL_N[[t]] \rtimes \C^*} (\Gr)$} is the fundamental class of the embedding $\Gr^{\omega_1} \hookrightarrow \Gr$, similarly for~\smash{$\bigl[\Gr^{\omega_1^*}\bigr]$}.
\begin{Theorem}
\label{thm:toda_comparison}
	The diagram
	\[
		\xymatrix{
			\rY^0_{-2N}(\sl_2) \ar[dr]\ar[d] & \\ \lhamp \dD(\GL_N) \hamp N_+^R \ar[r] & H_{\bullet}^{\GL_N[[t]] \rtimes \C^*} (\Gr)
		}
	\]
	is commutative. Here, the left vertical arrow is Theorem~$\ref{thm:toda_lattice_yangian}$ and the lower horizontal one is the composition of \eqref{eq::diff_ops_conj_inverse} and the derived Satake equivalence \eqref{eq::toda_grassmannian}.
	\begin{proof}
		Observe that the composition of Theorem~\ref{thm:kw_to_yangian} and Proposition~\ref{prop:top_vs_alg_vector} gives an algebra map
		\begin{align*}
			&\rY_{-N,N}(\gl_2)^{\geq 0} \rightarrow\U\bigl(\gl^R_N\bigr) \ltimes S^{\bullet} \C^N \hamp N_+^R \\
 &\hphantom{\rY_{-N,N}(\gl_2)^{\geq 0}}{}
 \rightarrow \lhamp \dD(\GL_N) \hamp N_+^R\rightarrow H_{\bullet}^{\GL_N[[t]] \rtimes \C^*} (\Gr), \\
			&d_1(u) \mapsto A(u), \qquad x_+(u) \mapsto (-1)^{N} \frac{1}{u - c_1^{\GL_N}(\cS)} \cap [\Gr^{\omega_1}].
		\end{align*}
		Likewise, the dual versions of Theorem~\ref{thm:kw_to_yangian_dual} and Proposition~\ref{prop:top_vs_alg_dual} give
		\begin{align*}
			&\rY_{-N,N}(\gl_2)^{\leq 0} \rightarrow\U\bigl(\gl^R_N\bigr) \ltimes S^{\bullet} \bigl(\C^N\bigr)^* \hamp N_+^R \\
 &\hphantom{\rY_{-N,N}(\gl_2)^{\leq 0}}{}
 \rightarrow \lhamp \dD(\GL_N) \hamp N_+^R \rightarrow H_{\bullet}^{\GL_N[[t]] \rtimes \C^*} (\Gr), \\
			&d_1(u) \mapsto A(u), \qquad x_-(u) \mapsto \frac{1}{u - c_1^{\GL_N}(\cQ) - \hbar} \cap \bigl[\Gr^{\omega_1^*}\bigr].
		\end{align*}
		By Theorem~\ref{thm:toda_lattice_yangian}, the combination of algebraic maps gives an isomorphism
		\[
			\rY^0_{-2N}(\sl_2) \xrightarrow{\sim} \lhamp \dD(\GL_N) \hamp N_+^R.
		\]
		Using explicit formulas \eqref{eq::shifted_yangian_topological} as well as relations \eqref{eq::chern_class_c_star} and \eqref{eq::chern_class_c_star_dual}, we conclude.
	\end{proof}
\end{Theorem}

\begin{Remark}
\label{remark:center_identification}
	The explicit identification \eqref{eq::center_topological} of \cite{BezrukavnikovFinkelberg} comes from two facts. On the algebraic side, in loc.\ cit., the authors identify the center $\Z_{\gl_N}$ with $W$-invariant functions on $\h^*$ shifted by $-\hbar \rho_N$ via the projection
	\[
		\Z_{\gl_N} \hookrightarrow \U(\gl_N) \rightarrow \n_+ \backslash \U(\gl_N),
	\]
	for instance, see \cite[Lemma~5]{BezrukavnikovFinkelberg}. By definition of the quantum minor \eqref{eq::quantum_minor}, we have
	\[
		T_{1\dots N}^{1\dots N}(u) \mapsto (u + E_{11} - N\hbar + \hbar) \cdots (u+E_{N-1,N-1} - \hbar) (u+E_{NN}),
	\]
	and its coefficients are exactly elementary symmetric functions in variables $\{E_{ii} - (N-i)\hbar\}$. On the topological side, this is formula \eqref{eq::shifted_yangian_topological}.
\end{Remark}

\subsection{Monopole operators}\label{subsect:monopole}

Recall the notations from Section~\ref{subsect:notations}. Let $\lambda$ be a dominant weight, and denote by ${v^*_{\lambda} \in V_{\lambda}^*}$ the~vector dual to the lowest-weight vector in the representation $V_{\lambda}$. Observe that the function~${\langle g\cdot v_{\lambda},v^*_{\lambda} \rangle}$ is $N_+^L \times N_+^R$-invariant, and hence defines an element in the Toda lattice \linebreak $\lhamp \dD(\GL_N) \hamp N_+^R$. The goal of this subsection is to compute these elements and to identify their images under the GKLO map for fundamental weights. This is an algebraic counterpart of the simplest case of the monopole operators formula \cite[formula~(A.5)]{BravermanFinkelbergNakajima}.

Under the embedding \smash{$V_{\omega_k} = \Lambda^i \C^N \hookrightarrow \bigl(\C^N\bigr)^{\otimes k}$}, let us choose a highest-weight vector and the dual of a lowest-weight vector as
\[
	v_{\omega_i} = \sum_{\sigma\in S_k} (-1)^{\sigma} v_{\sigma(k)} \otimes \dots \otimes v_{\sigma(1)}, \qquad v^*_{\omega_i} = \phi_{N} \wedge \dots \wedge \phi_{N-k+1} .
\]

\begin{Proposition}
\label{prop:image_of_minor_function}
	Under the embedding ${}^R \iota_N$ of Example~$\ref{ex:embedding_diff_ops}$, we have
	\[
		\bigl\langle g\cdot v_{\omega_k}, v^*_{\omega_k} \bigr\rangle = {}^R \iota_N \biggl[ \frac{1}{k!} \biggl(\prod_{1 \leq i \neq j \leq k} (\Omega_i - \Omega_j - \hbar)\biggr) (v_1)^{\otimes k} \biggr].
	\]
	\begin{proof}
		Observe that
$
			\bigl\langle g\cdot v_{\omega_k}, v^*_{\omega_k} \bigr\rangle = k! \langle g\cdot v_k \otimes \dots \otimes v_1, \phi_{N} \wedge \dots \wedge \phi_{N-k+1} \rangle$.
		Consider the map
		\[
			{}^L \iota^*_{k} \otimes \dots \otimes {}^L \iota^*_{1}\colon\ \U\bigl(\gl^L_N\bigr) \otimes \bigl(\bigl(\C^N\bigr)^*\bigr)^{\otimes k} \rightarrow \dD(\GL_N),
		\]
		which is the composition of
		\begin{align*}
			&\bigl(\U\bigl(\gl_N^L\bigr) \otimes \bigl(\C^N\bigr)^*\bigr) \otimes_{\U(\gl_N)} \dots \otimes_{\U(\gl_N)}\bigl(\U\bigl(\gl_N^L\bigr) \otimes \bigl(\C^N\bigr)^*\bigr) \\
			&\qquad \xrightarrow{{}^L \iota^*_{k} \otimes \dots \otimes {}^L \iota^*_{1}}\dD(\GL_N) \otimes_{\U(\gl_N^L)} \dots \otimes_{\U(\gl_N^L)} \dD(\GL_N) \xrightarrow{m} \dD(\GL_N)
		\end{align*}
		(the last arrow is the multiplication map) with the isomorphism
		\[
		\U\bigl(\gl_N^L\bigr) \otimes \bigl(\bigl(\C^N\bigr)^*\bigr)^{\otimes k} \cong \bigl(\U\bigl(\gl_N^L\bigr) \otimes \bigl(\C^N\bigr)^*\bigr) \otimes_{\U(\gl_N)} \dots \otimes_{\U(\gl_N)}\bigl(\U\bigl(\gl_N^L\bigr) \otimes \bigl(\C^N\bigr)^*\bigr).
		\]
		In what follows, we denote by the same letter the induced map of the Kostant--Whittaker reductions
		\[
		{}^L \iota^*_{k} \otimes \dots \otimes {}^L \iota^*_{1}\colon\ \lhamp \U\bigl(\gl^L_N\bigr) \otimes \bigl(\bigl(\C^N\bigr)^*\bigr)^{\otimes k} \rightarrow \lhamp \dD(\GL_N),
		\]
		we hope it will not lead to any confusion. By definition,
		\[
			\langle g\cdot v_k \otimes \dots \otimes v_1, \phi_{N-k+1} \wedge \dots \wedge \phi_N \rangle = \bigl({}^L \iota^*_{k} \otimes \dots \otimes {}^L \iota^*_{1}\bigr)(\phi_{N}\wedge \dots \wedge \phi_{N-k+1}),
		\]
		By Proposition~\ref{prop:skew_rep_vandermonde_dual},
		\begin{align*}
			&\bigl({}^L \iota^*_{k} \otimes \dots \otimes {}^L \iota^*_{1}\bigr)(\phi_{N}\wedge \dots \wedge \phi_{N-k+1}) \\
			&\qquad{}=\bigl({}^L \iota^*_{k} \otimes \dots \otimes {}^L \iota^*_{1}\bigr) \biggl[\frac{1}{k!} \prod_{i<j} \bigl(\Omega^*_i - \Omega^*_j + \hbar\bigr)\cdot \bigl(\phi_N^{\otimes k}\bigr) \biggr] \in \lhamp \dD(\GL_N)
		\end{align*}
		(recall that $\psi=1$). At the same time, it follows from Proposition~\ref{prop:canonical_matrix_desciption} and \eqref{eq::right_vs_left_matrix} that
		\[
			{}^L \iota_{\alpha}^* [\Omega^*(\phi_N)] = \bigl[XE^L\bigr]_{\alpha N} = \bigl[\bigl(-E^R + N\hbar\bigr) X\bigr]_{\alpha N} = {}^R \iota_{N} [(-\Omega+N\hbar)(v_{\alpha})]
		\]
		for any $\alpha$. Therefore, we have
		\begin{align*}
			&\bigl({}^L \iota^*_{k} \otimes \dots \otimes {}^L \iota^*_{1}\bigr) \biggl[\frac{1}{k!} \prod_{i<j} \bigl(\Omega^*_i - \Omega^*_j + \hbar\bigr)\cdot \bigl(\phi_N^{\otimes k}\bigr) \biggr] \\
			&\qquad{}= (-1)^{\frac{k(k-1)}{2}} \bigl({}^R \iota_{N} \otimes \dots \otimes {}^R \iota_{N}\bigr) \biggl[\frac{1}{k!}\prod_{i<j} (\Omega_i - \Omega_j - \hbar)\cdot (v_k \otimes \dots \otimes v_1) \biggr],
		\end{align*}
		where the map
		\[
			{}^R \iota_N \otimes \dots \otimes {}^R \iota_N \colon\ \U\bigl(\gl^R_N\bigr) \otimes \bigl(\C^N\bigr)^{\otimes k} \rightarrow \lhamp \dD(\GL_N)
		\]
		is defined similarly to the left version. Since
		\begin{gather*}
			\bigl({}^R \iota_{N} \otimes \dots \otimes {}^R \iota_{N}\bigr) \biggl[\prod_{i<j} (\Omega_i - \Omega_j - \hbar)\cdot (v_k \otimes \dots \otimes v_1) \biggr]\\
\qquad{} = \bigl\langle g\cdot v_{\omega_k}, v^*_{\omega_k} \bigr\rangle \in \lhamp \dD(\GL_N) \hamp N_+^R,
		\end{gather*}
		we have
		\begin{align*}
			&\bigl({}^R \iota_{N} \otimes \dots \otimes {}^R \iota_{N}\bigr) \biggl[\prod_{i<j} (\Omega_i - \Omega_j - \hbar)\cdot (v_k \otimes \dots \otimes v_1) \biggr] \\
			&\qquad{}=\bigl({}^R \iota_{N} \otimes \dots \otimes {}^R \iota_{N}\bigr) \biggl[\prod_{i<j} (\Omega_i - \Omega_j - \hbar)\cdot \pP (v_k \otimes \dots \otimes v_1) \pP \biggr]
		\end{align*}
(the argument of the right-hand side is considered as an element of $\U\bigl(\gl_N^R\bigr)\otimes \bigl(\C^N\bigr)^{\otimes k} \hamp N_+^R$). We have a commutative triangle
		\[
			\xymatrix{
				\U\bigl(\gl^R_N\bigr) \otimes \bigl(\C^N\bigr)^{\otimes k} \ar[rr]^{{}^R \iota_N \otimes \dots \otimes {}^R\iota_N} \ar[dr]_{e_k} & & \lhamp \dD(\GL_N) \\
				& \U\bigl(\gl^R_N\bigr) \otimes S^{\bullet} \C^N \ar[ur]^{{}^R \iota_N}&
			}
		\]
		where \smash{$e_k \colon \bigl(\C^N\bigr)^{\otimes k} \rightarrow S^{k} \C^N$} is the symmetrizer map. Also, by Proposition~\ref{prop:kw_triv_tensor_prod}, we have
		\[
			\pP^R v_k \otimes \dots \otimes v_1 \pP^R = \Omega_1^{k-1} \cdots \Omega_{k-1}^1 \Omega_k^0 \bigl(v_1^{\otimes k}\bigr) \pP^R \in \U\bigl(\gl_N^R\bigr) \otimes \bigl(\C^N\bigr)^{\otimes k} \hamp N_+^R.
		\]
		Therefore, we need to compute
		\[
			(-1)^{\frac{k(k-1)}{2}} {}^R\iota_N \biggl[ e_k\cdot \prod_{i<j} (\Omega_i - \Omega_j - \hbar) \Omega_1^{k-1} \cdots \Omega_{k-1}^1 \Omega_k^0 \bigl(v_1^{\otimes k}\bigr) \biggr] \in \U\bigl(\gl_N^R\bigr) \otimes \bigl(\C^N\bigr)^{\otimes k} \hamp N^R_+,
		\]
		which, by Proposition~\ref{prop:kw_tensor_algebra}, is given by the action of the element
		\[
			(-1)^{\frac{k(k-1)}{2}} \biggl({\rm e}^{\hbar}_k\cdot \prod_{i<j} (Y_i - Y_j -\hbar) Y_1^{k-1} \cdots Y_{k-1}^1 Y_k^0 \cdot {\rm e}^{\hbar}_k \biggr) \in \C[\hbar][Y_1,\dots,Y_k]^{S_k}
		\]
		in the spherical subalgebra of $\fH^{\hbar}_k$ as in \eqref{eq::hecke_spherical}. By Proposition~\ref{prop:symmetrizer_hecke_vs_shuffle}, it is equal to
		\begin{align*}
			& \frac{(-1)^{\frac{k(k-1)}{2}}}{k!} \Sym_k \biggl[ \prod_{i<j} \frac{(Y_i - Y_j+\hbar)(Y_i - Y_j - \hbar)}{Y_i - Y_j} Y_1^{k-1} \cdots Y_{k-1}^1 Y_k^0 \biggr] \\
			& \qquad{}=\frac{(-1)^{\frac{k(k-1)}{2}}}{k!} \prod_{i<j} (Y_i - Y_j +\hbar) (Y_i - Y_j - \hbar) \cdot \Sym_k \biggl[ \prod_{i<j} (Y_i - Y_j) Y_1^{k-1} \cdots Y_{k-1}^1 Y_k^0 \biggr].
		\end{align*}
		Since \smash{$\prod_{i<j} (Y_i - Y_j)$} is the Vandermonde determinant (in particular, skew-symmetric), we have
		\[
			\Sym_k \biggl[ \prod_{i<j} (Y_i - Y_j) Y_1^{k-1} \cdots Y_{k-1}^1 Y_k^0 \biggr] = \prod_{i<j} (Y_i - Y_j) \sum_{\sigma\in S_k} (-1)^{\sigma} Y_{\sigma(1)}^{k-1} \cdots Y_{\sigma(k-1)}^1 Y_{\sigma(k)}^0 = 1.
		\]
		Therefore, this element is equal to $\frac{1}{k!}\prod_{i\neq j} (Y_i - Y_j - \hbar)$, and so
		\[
			\bigl\langle g\cdot v_{\omega_i}, v^*_{\omega_i} \bigr\rangle = {}^R \iota_N \biggl[ \frac{1}{k!} \biggl(\prod_{1 \leq i \neq j \leq k} (\Omega_i - \Omega_j - \hbar)\biggr) (v_1)^{\otimes k} \biggr],
		\]
		as required.
	\end{proof}
\end{Proposition}

The next theorem computes the images of these elements under the GKLO map of Theorem~\ref{thm:toda_lattice_yangian}.

\begin{Theorem}
\label{thm:monopole_operators_gklo}
	Under the GKLO homomorphism
	\[
		\lhamp \dD(\GL_N) \hamp N_+^R \rightarrow \cA
	\]
	of Theorem~$\ref{thm:toda_lattice_yangian}$, we have
	\[
		\bigl\langle g\cdot v_{\omega_k}, v^*_{\omega_k} \bigr\rangle \mapsto \sum_{\substack{I \subset \{1,\dots,N\}, \\ |I| = k}} \prod_{\substack{i \in I, \\ j\in \bar{I}}} \frac{1}{w_i - w_j} \prod_{i\in I} \bu_i^{-1},
	\]
	where $\bar{I}$ is the complement of $I$ in $\{1,\dots, N\}$.
	\begin{proof}
		Since \smash{$\bigl\langle g\cdot v_{\omega_k}, v^*_{\omega_k} \bigr\rangle$} is in the image of ${}^R\iota_N$, it is enough to prove the statement for the Miura transform
		\[
			\U(\gl_N) \ltimes S^{\bullet} \C^N \hamp N_+ \rightarrow \cA^-
		\]
		of Theorem~\ref{thm:miura_vs_gklo}. Recall that
		\[
			v_1 \mapsto \sum_{i=1}^N \prod_{j\neq i} \frac{1}{w_i - w_j} \bu_i^{-1}.
		\]
		Since the Miura transform is an algebra map, we have
		\[
			S^k \C^N \ni v_1^{\otimes k} \mapsto \Biggl( \sum_{i=1}^N \prod_{j\neq i} \frac{1}{w_i - w_j} \bu_i^{-1} \Biggr)^k =: \sum_J c_J(w) \bu_{j_1}^{-1} \cdots \bu_{j_k}^{-1},
		\]
		where $J$ runs over all multisubsets of $\{1,\dots,N\}$ and $c_J(w)\in \U(\h)^{\gen}$ are some functions. In~what follows, we will use Proposition~\ref{prop:spherical_action_par_res} and notations thereof. By Proposition~\ref{prop:image_of_minor_function}, we need to compute
		\[
			\sum_J c_J(w) f\bigl(\vec{Y}\bigr) \bu^{-1}_{j_1} \cdots \bu^{-1}_{j_k}, \qquad f\bigl(\vec{Y}\bigr) = \frac{1}{k!} \prod_{i\neq j} (Y_i - Y_j - \hbar).
		\]
		If $J$ contains repeating indices that we can assume to be $j_1 = j_2$ without loss of generality, then
		\[
			(Y_1 - Y_2 - \hbar) \bu_{j_1}^{-1} \bu_{j_1}^{-1} = w_{j_1}\bu^{-2}_{j_1} - \bu_{j_1}^{-1} w_{j_1} \bu_{j_1}^{-1} - \hbar \bu^{-2}_{j_1} = 0.
		\]
		Therefore, the only non-zero indices in the sum are permutations of $\{1,\dots,k\}$. Moreover, by $W$-invariance of Remark~\ref{remark:w_action_difference_operators}, it is enough to compute the coefficient $c_{\{1,\dots,k\}}(w)$. By an easy induction on $k$, we have
		\[
			c_{\{1,\dots,k\}}(w) = k! \prod_{1 \leq i \neq j \leq k} \frac{1}{w_i - w_j - \hbar} \cdot \prod_{\substack{i=1,\dots, k, \\ j = k+1, \dots, N}} \frac{1}{w_i - w_j} \bu^{-1}_1 \cdots \bu^{-1}_k.
		\]
		 By Proposition~\ref{prop:spherical_action_par_res},
		 \[
		 	c_{\{1,\dots,k\}}(w) f\bigl(\vec{Y}\bigr) \bu_1^{-1} \cdots \bu_k^{-1} = \prod_{\substack{i =1,\dots, k, \\ j = k+1, \dots, N}} \frac{1}{w_i - w_j} \bu_1^{-1} \cdots \bu_k^{-1}.
		 \]
		 By $W$-invariance, we obtain
		 \[
			 \bigl\langle g\cdot v_{\omega_k}, v^*_{\omega_k} \bigr\rangle \mapsto \sum_{\substack{I \subset \{1,\dots,N\}, \\ |I| = k\}}} \prod_{\substack{i \in I, \\ j\in \bar{I}}} \frac{1}{w_i - w_j} \prod_{i\in I} \bu_i^{-1},
		 \]
		 as required.
	\end{proof}
\end{Theorem}

\begin{Remark}
	This is a rational analog of the simplest operator in the \emph{M-system} in its shuffle realization, see \cite[formula~(104)]{TsymbaliukGKLO} and \cite[Proposition 6.13]{TsymbaliukGKLO}. By \cite[Remark~9.19]{FinkelbergTsymbaliuk}, its operators correspond to the $K$-theory classes of line bundles $\cO(-k)$ restricted to the orbits~$\Gr^{\omega_i}$ as in Section~\ref{subsect:comparison}, where $\cO(1)$ is the theta bundle on the affine Grassmannian. Therefore, Proposition~\ref{prop:image_of_minor_function} and Theorem~\ref{thm:monopole_operators_gklo} are in accordance with the derived Satake isomorphism of Theorem~\ref{thm:toda_comparison}, since by \cite[Lemma~2.6]{FKPRW}, we have \smash{$\bigl\langle g\cdot v_{\omega_i}, v_{\omega_i}^* \bigr\rangle \mapsto [\Gr^{\omega_i}]$}, up to a scalar factor. It would be interesting to find a rational analog of the full $M$-system.
\end{Remark}

\appendix

\renewcommand{\theequation}{\thesection.\arabic{equation}}

\section{Technical lemmas}
\label{appendix:technical_lemmas}

In this appendix, we will prove some technical lemmas from the main section.

The following result is used in the proof of Proposition~\ref{prop:su_symmetrizer}. Recall that
\[
	G^{\theta}(z,w) = \frac{\theta(z+w)\theta'(0)}{\theta(z)\theta(w)},
\]
where $\theta(z)$ is any skew-symmetric function satisfying the three-term relation \eqref{eq::three_term_relation}.
\begin{Lemma}
\label{lm:generalized_three_term}
	We have
	\[
		\prod_{i=1}^l G^{\theta}(u-z_i,\kappa) = \sum_{i=1}^l G^{\theta}(u-z_i,l\kappa) \cdot \prod_{\substack{j=1,\dots, l ,\\ j\neq i}} G^{\theta}(z_i-z_j,\kappa).
	\]
	\begin{proof}
		We can show it by induction on $l$. The base $l=1$ is tautological. For the step, let us multiply both sides by $G^{\theta}(u-z_{l+1},\kappa)$. The three-term relation and skew-symmetricity of $\theta(t)$ imply
		\begin{align*}
			&\frac{\theta(u-z_i+l\kappa) \theta(u-z_{l+1} + \hbar)}{\theta(u-z_i) \theta(l\kappa)\theta(u-z_{l+1}) \theta(\kappa)} \\
			&\qquad{}=\frac{\theta(u-z_{l+1} +l\kappa+\kappa) \theta(z_{l+1}-z_i+l\kappa)}{\theta(u-z_{l+1}) \theta(l\kappa+\kappa) \theta(z_{l+1}-z_i) \theta(l\kappa)} + \frac{\theta(z_i-z_{l+1} + \kappa) \theta(u-z_i + l\kappa+\kappa)}{\theta(z_i-z_{l+1}) \theta(\kappa) \theta(u-z_i) \theta(l\kappa+\kappa)}
		\end{align*}
		by multiplying by common denominator and taking
		\begin{alignat*}{3}
			&x = \frac{l\kappa+\kappa +2u - z_i - z_{l+1}}{2}, \qquad&& y = \frac{l\kappa-\kappa-z_i + z_{l+1}}{2}, &\\
			&z = \frac{l\kappa+\kappa - z_i + z_{l+1}}{2}, \qquad&& w = \frac{l\kappa+\kappa + z_i - z_{l+1}}{2}&
		\end{alignat*}
		in the notations of \eqref{eq::three_term_relation}. Then the statement for $l+1$ follows from the induction assumption for $u=z_{l+1}$.
	\end{proof}
\end{Lemma}

The following result about symmetric polynomials is used in Proposition~\ref{prop:kw_tensor_algebra}. Define symmetric polynomials \smash{$p^{\hbar}_a\bigl(\vec{\Omega}\bigr) \in \C[\hbar][\Omega_1,\dots,\Omega_k]^{S_k}$} by
\[
\prod_{i=1}^k \frac{u-\Omega_i +\hbar}{u-\Omega_i} = 1 + \hbar \sum_{a=0}^{\infty} p_a^{\hbar} \bigl(\vec{\Omega}\bigr) u^{-a-1} \in \C[\hbar][\Omega_1,\dots,\Omega_k]^{S_k} \bigl[\hspace{-0.2mm}\bigl[u^{-1}\bigr]\hspace{-0.2mm}\bigr],
\]
where we understand each term in the product as in Remark~\ref{remark:pwoer_series}.
\begin{Lemma}\qquad
	\label{lm:symmetric_func}
	\begin{enumerate}\itemsep=0pt
		\item[$1.$] 	We have
		\begin{equation}
			\label{eq_prop:p_hbar_function}
			p_a^{\hbar}\bigl(\vec{\Omega}\bigr) = \sum_{i=1}^k \Omega_i^a \prod_{j\neq i} \frac{\Omega_i-\Omega_j + \hbar}{\Omega_i - \Omega_j}.
		\end{equation}
		In particular, \smash{$p^{0}_a\bigl(\vec{\Omega}\bigr) = \sum_{i=1}^k \Omega_i^a$} is the standard power sum symmetric function.
		
		\item[$2.$] The functions \smash{$p^{\hbar}_a\bigl(\vec{\Omega}\bigr)$} for $1\leq a \leq k$ freely generate the algebra of symmetric functions, in other words, there is an isomorphism
		\begin{equation*}
			\C[\hbar][\Omega_1,\dots,\Omega_k]^{S_k} \cong \C[\hbar]\bigl[p^{\hbar}_1\bigl(\vec{\Omega}\bigr),\dots,p^{\hbar}_k\bigl(\vec{\Omega}\bigr)\bigr].
		\end{equation*}
	\end{enumerate}
	\begin{proof}
		The first part of the proposition follows by using the partial fraction decomposition
		\[
		\prod_{i=1}^k \frac{u-\Omega_i +\hbar}{u-\Omega_i} = 1 + \sum_{i=1}^k \frac{\hbar}{u-\Omega_i} \prod_{j\neq i} \frac{\Omega_i-\Omega_j+\hbar}{\Omega_i-\Omega_j}.
		\]
		For the second part, recall that there is an isomorphism
		\[
		\C[\hbar][\Omega_1,\dots,\Omega_k]^{S_k} \cong \C[\hbar]\bigl[p^{0}_1\bigl(\vec{\Omega}\bigr),\dots,p^{0}_k\bigl(\vec{\Omega}\bigr)\bigr],
		\]
		since \smash{$p^{0}_a\bigl(\vec{\Omega}\bigr) = \sum_{i=1}^k \Omega_i^a$}. Let us consider the right-hand side as the coordinate ring of the affine space $\mathbb{A}^k \times \mathbb{A}^1$, where $\mathbb{A}^1$ corresponds to $\hbar$. There is a $\C^*$-action on $\C[\hbar][\Omega_1,\dots,\Omega_k]$ such that all the variables $\Omega_i$ and $\hbar$ have weight 1. It induces a $\C^*$-action on $\C[\hbar][\Omega_1,\dots,\Omega_k]^{S_k}$ with positive weights; in particular, it contracts $\mathbb{A}^k \times \mathbb{A}^1$ to the origin. In fact, under the isomorphism above, this action is linear and diagonal with weights $(1,2,\dots,k,1)$.
		
		Consider the map \smash{$\varphi_{\hbar} \colon \mathbb{A}^k \times \mathbb{A}^1 \rightarrow \mathbb{A}^k \times \mathbb{A}^1$} defined on the ring of functions by \smash{$p^{0}_a\bigl(\vec{\Omega}\bigr) \mapsto p^{\hbar}_a\bigl(\vec{\Omega}\bigr)$} for every $a$ and $\hbar \mapsto \hbar$. It follows from \eqref{eq_prop:p_hbar_function} that this map is $\C^*$-equivariant. Moreover, it is clear from weight considerations and explicit calculations for $a=1$ that the linear part of $\varphi_{\hbar}$ is
		\[
		\mathrm{d}\varphi_{\hbar} \bigl(p^0_a\bigl(\vec{\Omega}\bigr)\bigr) = p_a^0\bigl(\vec{\Omega}\bigr) + \delta_{a1} \frac{k(k-1)}{2}\hbar, \qquad \mathrm{d}\varphi_{\hbar}(\hbar) = \hbar.
		\]
		
		Therefore, the differential of $\varphi_{\hbar}$ at zero is invertible. By the proof of \cite[Lemma~2.1]{GanGinzburg}, it is an isomorphism. In particular, the polynomials \smash{$\bigl\{p^{\hbar}_a\bigl(\vec{\Omega}\bigr)\bigr\}$} generate the ring of symmetric functions.
	\end{proof}
	
\end{Lemma}

In what follows, we use the notations of Section~\ref{sect:parabolic_reduction} for $\U(\gl_N) \ltimes \bigl(S^{\bullet} \C^N \otimes S^{\bullet}\bigl(\C^N\bigr)^*\bigr)$ (more precisely, for its Mickelsson algebra).

\begin{Proposition}
	\label{prop:diff_ops_commute}
	For any $i\neq j$, we have $\bar{\phi}_i \bar{v}_j = \bar{v}_j \bar{\phi}_i$.
	\begin{proof}
		Assume $i<j$. By \eqref{eq::parabolic_res_vector_generators} and \eqref{eq::par_res_dual_difference_vs_extremal_proj}, it is enough to prove
		\[
		P\phi_i P \cdot P v_j P = P v_j P \cdot P\phi_i P.
		\]
		Using the properties of the extremal projector \eqref{eq::extremal_proj} \big(in particular $P^2=P$\big) and explicit formula~\eqref{eq::par_res_dual_vec}, the left-hand side is given by
		\begin{align*}
			&P \sum_{1\leq i_1 < \dots < i_s < i} (w_i - w_{i_1})^{-1} \cdots (w_i - w_{i_s})^{-1} \phi_{i_1} E_{i_1 i_2} \cdots E_{i_s i} \cdot v_j P \\
			& \qquad{}=P v_j \sum_{1\leq i_1 < \dots < i_s < i} (w_i - w_{i_1})^{-1} \cdots (w_i - w_{i_s})^{-1} \phi_{i_1} E_{i_1 i_2} \cdots E_{i_s i} P = Pv_j \phi_i P,
		\end{align*}
		where we used $i<j$ (hence $v_j$ commutes with everything). Likewise, using explicit formula~\eqref{eq::parabolic_res_dual_difference_gens}, the right-hand side is given by
		\begin{align*}
			P \sum_{N \geq j_1 > \dots > j_s > j} (w_j - w_{j_1})^{-1} \cdots (w_i - w_{j_s})^{-1} v_{j_1} E_{j_2 j_1} \cdots E_{jj_s} E_{i_1 i_2} \cdot \phi_i P = P \phi_i v_j P.
		\end{align*}
		Using $Pv_j \phi_i P = P \phi_i v_j P$, we conclude for $i<j$. For $i>j$, define the antiautomorphism
		\begin{align*}
			&\varpi_N \colon\ \U(\gl_N) \ltimes \bigl(S^{\bullet} \C^N \otimes S^{\bullet}\bigl(\C^N\bigr)^*\bigr) \rightarrow \bigl[\U(\gl_N) \ltimes \bigl(S^{\bullet} \C^N \otimes S^{\bullet}\bigl(\C^N\bigr)^*\bigr)\bigr]^{\mathrm{opp}}, \\
			&\varpi_N(v_i) = \phi_i, \qquad \varpi_N(\phi_i) = v_i,
		\end{align*}
		extending the one \eqref{eq::transposition} of $\U(\gl_N)$. By Remark~\ref{remark:transposition}, we have $\varpi_N(P) = P$, therefore,
		\[
		P \phi_j P v_i P = \varpi_N (P \phi_i P v_j P) = \varpi_N (P v_j P \phi_i P) = P v_i P \phi_j P,
		\]
		that implies the corresponding statement for the case $i>j$.
	\end{proof}
\end{Proposition}

The proof of the following proposition is a word-for-word repetition of \cite[Theorem~3.1]{MolevTransvectorAlgebras} adapted to the setting of this paper.
\begin{Proposition}
[{\cite[Theorem~3.1]{MolevTransvectorAlgebras}}]
\label{prop:comatrix_interpolation}
	We have
	\[
		(-1)^{N+1}P \sum_{i,j} \phi_i \hat{T}_{ij}(-u+N\hbar-\hbar) v_j P = \sum_{a=1}^N \prod_{b\neq a} \frac{u-w_b}{w_a - w_b} \cdot \bar{\phi}_a \bar{v}_a.
	\]
	\begin{proof}
		Denote by $\alpha = N+1$. Recall the trick of Example~\ref{ex:trick}. Denote by \smash{$T_{1\dots N}^{1\dots j\mapsto \alpha \dots N} (u)$} the quantum minor with values in $\U(\gl_{N+1})$ corresponding to the indices
		\[
		(b_1,\dots,b_N) = (1,\dots, N), \qquad (a_1,\dots,a_N) = (1,\dots,\alpha,\dots,N),
		\]
		in the notations of \eqref{eq::quantum_minor}, where $\alpha$ is on the $j$-th place. By Proposition~\ref{prop::minor_basic} and \eqref{eq::quantum_comatrix}, we have
		\[
			T_{1\dots N}^{1\dots j\mapsto \alpha \dots N} (u) = \sum_{i=1}^N (-1)^{i+j} E_{\alpha i} T_{1\dots \hat{i} \dots N}^{1\dots \hat{j} \dots N}(u) = \sum_{i=1}^N \phi_i \hat{T}_{ij}(u).
		\]
		Therefore, we can study the image of the left-hand side in the quotient $\n_- \backslash \U(\gl_{N+1})$, where $\n_- \subset \gl_N$. Recall the definition of the quantum minor \eqref{eq::quantum_minor} and the extremal projector notation from Section~\ref{sect:parabolic_reduction}. Using the properties \eqref{eq::extremal_proj} and pushing the terms of $\n_-$ to the left, we obtain
		\begin{align*}
			P T_{1\dots N}^{1\dots j\mapsto \alpha \dots N} (-u+N\hbar-\hbar) = \sum_I (-1)^s \prod_{k \in \bar{I}} (-u+w_k) P \phi_{i} E_{i i_2} \cdots E_{i_s j} \cdot \prod_{l > j} (-u + w_l),
		\end{align*}
		where $I = \{ 1 \leq i=i_1 < i_2 <\dots < i_s < j\}$ and $\bar{I}$ is the complement of $I$ in the set $\{1,\dots,j-1\}$. Using properties \eqref{eq::extremal_proj} again and pushing terms of $\n_+$ to the right, we get
		\begin{align*}
			(-1)^{N+1} \sum_{j=1}^N P T_{1\dots N}^{1\dots j\mapsto \alpha \dots N} (-u+N\hbar-\hbar) v_j P = \sum_{i=1}^N \prod_{k<i} (u-w_k)\cdot \prod_{k>i} (u-w_k + \hbar) P\phi_i v_i P,
		\end{align*}
		in particular, for any $a$, we have
		\[
			(-1)^{N+1} \sum_{i,j} P \phi_i \hat{T}_{ij}(-w_a+N\hbar-\hbar) v_j P = \sum_{i=1}^a \prod_{k<i} (w_a-w_k)\cdot \prod_{k>i} (w_a-w_k + \hbar) P\phi_i v_i P.
		\]
		At the same time, by explicit formula \eqref{eq::par_res_dual_vec}, we have
		\[
			\bar{\phi}_a v_a P= \sum_{i=1}^a \prod_{k=1}^{i-1}(w_a - w_k) \cdot \prod_{k=i+1}^{a}(w_a - w_{k}+\hbar) \cdot P \phi_i v_i P.
		\]
		Using relation \eqref{eq::parabolic_res_vector_generators}, we conclude that
		\[
			(-1)^{N+1} \sum_{i,j} P \phi_i \hat{T}_{ij}(-w_a+N\hbar-\hbar) v_j P = \bar{\phi}_a \bar{v}_a,
		\]
		therefore, by interpolation,
		\[
			(-1)^{N+1} \sum_{i,j} P \phi_i \hat{T}_{ij}(-u+N\hbar-\hbar) v_j P = \sum_{a=1}^N \prod_{b\neq a} \frac{u-w_b}{w_a - w_b} \cdot \bar{\phi}_a \bar{v}_a,
		\]
		as required.
	\end{proof}
\end{Proposition}

\section{Dual representation}
\label{subsect:dual_rep}

In this appendix, we present the facts concerning the Harish-Chandra bimodule $\U(\gl_N) \ltimes S^{\bullet} \bigl(\C^N\bigr)^*$. Unless otherwise stated, their proofs mirror the ones for $\U(\gl_N) \ltimes S^{\bullet} \C^N$ and will be omitted.

Define the dual analog of the canonical homomorphism \eqref{eq::canonical_homomorphism} by
\begin{equation}
\label{eq::canonical_homomorphism_dual}
	\Omega^* \colon\ \U(\gl_N) \otimes \bigl(\C^N\bigr)^* \rightarrow \U(\gl_N) \otimes \bigl(\C^N\bigr)^*, \qquad x\phi_i \mapsto \sum_j \phi_j E_{ji}.
\end{equation}
In matrix form,
\begin{equation*}
	(\Omega^*(\phi_1),\dots,\Omega^*(\phi_N))^{\mathsf T} = \bigl(E^{\mathsf T}+N\hbar\cdot\mathrm{Id}\bigr) \cdot (\phi_1,\dots,\phi_N)^{\mathsf T}.
\end{equation*}
In particular, the analog of the Cayley--Hamilton theorem of Proposition~\ref{prop:cayley_hamilton} reads as
\begin{equation*}
	A(\Omega^*-\hbar) = 0,
\end{equation*}
which follows directly from \cite[Theorem~7.2.1]{MolevYangians}. Similarly to Proposition~\ref{prop:center_left_vs_right}, one can show that%
\begin{equation}
\label{eq::center_conjugation_dual}
	A(u)^{-1} \phi A(u) = \frac{u-\Omega^*}{u-\Omega^*+\hbar} \phi.
\end{equation}

Consider the tensor product
\[
	\U(\gl_N) \otimes \bigl(\bigl(\C^N\bigr)^*\bigr)^{\otimes k} \cong \U(\gl_{N}) \otimes \bigl(\C^N\bigr)^* \otimes_{\U(\gl_N)} \dots \otimes_{\U(\gl_N)} \U(\gl_{N}) \otimes \bigl(\C^N\bigr)^*.
\]
As in \eqref{eq::affine_hecke_action_vector}, we can define an action of the degenerate affine Hecke algebra $\fH^{-\hbar}_k$ on it (observe a~different quantization constant):
\begin{align*}
		&\Omega^*_i := 1 \otimes \dots \otimes \Omega^* \otimes \dots \otimes 1 \curvearrowright \U(\gl_N) \otimes \bigl(\C^N\bigr)^* \otimes_{\U(\gl_N)} \dots \otimes_{\U(\gl_N)} \U(\gl_N) \otimes \bigl(\C^N\bigr)^*, \\
		&\sigma_i(x \otimes \phi_{a_1} \otimes \dots \otimes \phi_{a_{i}} \otimes \phi_{a_{i+1}} \otimes \dots \otimes \phi_{a_k}) := x \otimes \phi_{a_1} \otimes \dots \otimes \phi_{a_{i+1}} \otimes \phi_{a_{i}} \otimes \dots \otimes \phi_{a_k}.
\end{align*}

Recall the notations for the parabolic restriction of Section~\ref{sect:parabolic_reduction}. As in \eqref{eq::parabolic_res_dual_difference_gens}, we can define
\begin{equation}
	\label{eq::par_res_dual_vec}
	\bar{\phi}_i = \sum\limits_{1 \leq i_1 < \dots < i_s < i} P (w_i - w_{j_1}) \cdots (w_i - w_{j_r}) \phi_{i_1} E_{i_{1} i_{2}} \cdots E_{i_s i} \in \n_- \backslash \U(\gl_N) \otimes \bigl(\C^N\bigr)^*,
\end{equation}
where $s=0,1,\dots$ and $\{j_1,\dots,j_r\}$ is the complementary subset to $\{i_1,\dots,i_s\}$ in $\{1\dots,i-1\}$. By the trick of Example~\ref{ex:trick} and \cite[formula~(3.3)]{MolevTransvectorAlgebras}, it is essentially given by the action of the extremal projector
\begin{equation}
\label{eq::par_res_dual_difference_vs_extremal_proj}
	\bar{\phi}_i = (w_i - w_1) \cdots (w_i - w_{i-1}) P\phi_i P.
\end{equation}

For each $i$, we have
\begin{equation}
\label{eq::omega_par_res_dual}
	P\Omega^*(\phi_i) P = (w_i+\hbar) P\phi_i P \in \U(\gl_N) \otimes \bigl(\C^N\bigr)^* \sslash N_-
\end{equation}
as in Lemma~\ref{lm:parabolic_res_canonical_homomorphism}.

Recall the positive part $\cA^+$ of the algebra of difference operators $\cA$ from Section~\ref{subsect:mickelsson}. We have an analog of Theorem~\ref{thm:par_res_to_difference_ops}
\begin{equation*}
	\U(\gl_N)^{\gen} \ltimes S^{\bullet}\bigl(\C^N\bigr)^* \sslash N_- \xrightarrow{\sim} \cA^+, \qquad \bar{\phi}_i \mapsto \bu_i.
\end{equation*}

Next, recall the Kostant--Whittaker reduction of Section~\ref{sect:kw_reduction}, in particular the Kirillov projector~$\pP$. The following equality can be proved analogously to \cite[Proposition 5.20]{KalmykovYangians}.
\begin{equation}
\label{eq::last_row_vs_center}
	u^N - \sum_{i=1}^N \psi^{N-i} u^{i-1} \pP E_{Ni} \pP = \pP A(u) \pP.
\end{equation}

The following is an analog of Proposition~\ref{prop:kw_vector}.
\begin{Proposition}
\label{prop:dual_kw_presentation}
	The action on $\phi_N$ gives a left $\Z_{\gl_N}$-module isomorphism
	\[
		\Z_{\gl_N}[\Omega^*]/(A(\Omega^*-\hbar)) \xrightarrow{\sim} \U(\gl_N) \otimes \bigl(\C^N\bigr)^* \hamp N_+, \qquad f \mapsto f(\phi_N).
	\]
	\begin{proof}
		Observe that by \eqref{eq::kw_trivialization}, the $\Z_{\gl_N}$-linear map
		\[
			\Z_{\gl_N} \otimes \bigl(\C^N\bigr)^* \xrightarrow{\sim} \U(\gl_N) \otimes \bigl(\C^N\bigr)^* \hamp N_+, \qquad \phi_i \mapsto \pP \phi_i \pP
		\]
		is an isomorphism. Therefore, it will be enough to relate the action of $\Omega^*$ to that of the Kirillov projector. By pushing the $\gl_N$-terms to the left or to the right and using of \eqref{eq::kirillov_proj_properties} as in \cite[Proposition 6.18]{KalmykovYangians}, we deduce
		\[
			\pP \Omega^*(\phi_i) \pP = \sum_{j=1}^N \pP \phi_j E_{ji} \pP = \pP E_{Ni} \phi_N \pP +\psi \pP \phi_{i-1} \pP + \hbar \pP \phi_i \pP.
		\]
		Since $\phi_N$ is a highest-weight vector, we have $\phi_N \pP = \pP \phi_N \pP$. Applying \eqref{eq::last_row_vs_center}, we deduce
		\[
			\pP E_{Ni} \phi_N \pP = \pP E_{Ni} \pP \phi_N \pP = -\psi^{i-N} A_{N-i+1} \phi_N\pP,
		\]
		and so
		\[
			\pP \Omega^*(\phi_i) \pP = -\psi^{i-N} A_{N-i+1} \pP \phi_N \pP +\psi \pP \phi_{i-1} \pP +\hbar \pP \phi_i \pP.
		\]
		In particular, the matrix of $\zeta= \Omega^*-\hbar$ in the generators $\{ \pP \phi_i \pP \mid i=1,\dots, N\} $ has the form
		\[
			\begin{pmatrix}
				0 & \psi & \dots & 0 & 0 \\
				0 & 0 & \dots & 0 & 0 \\
				\vdots & \vdots & \ddots & \vdots & \vdots & \\
				0 & 0 & \dots & 0 & \psi \\
				-\psi^{-N+1} A_N & -\psi^{-N+2} A_{N-1} & \dots & -\psi^{-1} A_2 & -A_{1}
			\end{pmatrix}
		\]
		Define polynomials $Q_i(\zeta)\in \Z_{\gl_N}[\zeta]$ for $i=1,\dots, N$ inductively by the rule
		\[
			Q_N(\zeta) = 1, \qquad Q_{i-1} (\zeta) = \psi^{-1}\bigl[\psi^{N-i} A_{N-i} + \zeta Q_i (\zeta)\bigr].
		\]
		The transition matrix from $\bigl(1,\dots,\zeta^{N-1}\bigr)$ to $(\psi_N(\zeta),\dots,\psi_1(\zeta))$ is upper triangular with invertible diagonal, hence is invertible itself. However, observe that the action of $\zeta$ on the latter poly\-nomials is given by the matrix above. In particular, since the elements $\{ \pP \phi_i \pP \mid i=1,\dots, N\}$ \mbox{generate} \smash{$\U(\gl_N) \otimes \bigl(\C^N\bigr)^* \hamp N_+$} as a left $\Z_{\gl_N}$-module, the elements \smash{$\bigl\{(\Omega^*-\hbar)^{i-1}(\phi_N) \mid$} \linebreak $i=1,\dots, N\bigr\}$ do as well. However, the latter generate $\Z_{\gl_N}[\Omega^*]/(A(\Omega^*-\hbar))$.
	\end{proof}
\end{Proposition}

For $1\leq j \leq k+1$, define polynomials
\[
	A^{(-j)}(u) = \prod_{i=1}^{j-1}\frac{u-\Omega^*_i}{u-\Omega^*_i + \hbar} A(u)
\]
with values in $\Z_{\gl_N}\bigl[\Omega^*_1,\dots,\Omega^*_k\bigr]$ as in \eqref{eq::modified_center}. Likewise, consider intersection
\[
	I^*_k := \bigl(A^{(-j)}(\Omega_j^* - \hbar)\mid j=1\dots k\bigr) \cap \Z_{\gl_N} \bigl[\Omega^*_1,\dots, \Omega^*_k\bigr].
\]
The following is an analog of Proposition~\ref{prop:kw_tensor_algebra} and Corollary~\ref{cor:kw_symmetric_power}.
\begin{Proposition}\quad
\label{prop:kw_tensor_algebra_dual}
	\begin{enumerate}\itemsep=0pt
		\item[$1.$] The action on \smash{$\phi_N^{\otimes k}$} gives an isomorphism of left $\Z_{\gl_N}$-modules
		\[
			\Z_{\gl_N}\bigl[\Omega^*_1,\dots, \Omega^*_k\bigr]/ \bigl(A^{(-j)}(\Omega^*_j-\hbar) \mid j=1,\dots, k \bigr) \xrightarrow{\sim} \U(\gl_N) \otimes \bigl(\bigl(\C^N\bigr)^*\bigr)^{\otimes k} \hamp N_+,
		\]
		compatible with the action of $\fH^{-\hbar}_k$.
		\item[$2.$] There is an isomorphism
		\[
			\Z_{\gl_N}\bigl[\Omega^*_1,\dots, \Omega^*_k\bigr]^{S_k} / I^*_k \xrightarrow{\sim} \U(\gl_N) \otimes S^k \bigl(\C^N\bigr)^* \hamp N_+.
		\]
	\end{enumerate}
\end{Proposition}

Recall the definition of the non-positive part $\rY_{-N,N}(\gl_2)^{\leq 0}$ of the shifted Yangian from Definition~\ref{def:shifted_yangian}. The following is an analog of Theorem~\ref{thm:kw_to_yangian}.
\begin{Theorem}
\label{thm:kw_to_yangian_dual}
	There is a surjective homomorphism of $\C[\hbar]$-algebras
	\[
	\rY_{-N,N}(\gl_2)^{\leq 0} \rightarrow \U(\gl_N) \ltimes S^{\bullet} \bigl(\C^N\bigr)^* \hamp N_+,
	\]
	defined by
	\begin{align*}
		d_1(u) &\mapsto A(u), \qquad d_2(u) \mapsto A(u-\hbar)^{-1}, \\
		x_-(u) &\mapsto \frac{1}{u-\Omega^*_1} (\phi_N) \in \bigl(\U(\gl_N) \otimes \bigl(\C^N\bigr)^* \hamp N_+ \bigr)\bigl[\hspace{-0.2mm}\bigl[u^{-1}\bigr]\hspace{-0.2mm}\bigr].
	\end{align*}
\end{Theorem}

Finally, recall the Miura transform of Definition~\ref{def:miura_transform}:
\begin{equation}
\label{eq::miura_to_gklo_dual}
	\rY_{-N,N}(\gl_2)^{\leq 0} \rightarrow \U(\gl_N) \ltimes S^{\bullet}\bigl(\C^N\bigr)^* \hamp N_+ \rightarrow \U(\gl_N) \ltimes S^{\bullet} \bigl(\C^N\bigr)^* \sslash N_- \cong \cA^-.
\end{equation}
Using \eqref{eq::omega_par_res_dual}, the following theorem can be proved analogously to Theorem~\ref{thm:miura_vs_gklo}.
\begin{Theorem}
	The map \eqref{eq::miura_to_gklo_dual} coincides with the GKLO map from Theorem~\ref{thm:gklo}.
\end{Theorem}

The following result is used in Section~\ref{subsect:monopole}.
\begin{Proposition}
\label{prop:skew_rep_vandermonde_dual}
	Under the embedding $\Lambda^k \bigl(\C^N\bigr)^* \hookrightarrow \bigl(\bigl(\C^N\bigr)^*\bigr)^{\otimes k}$, we have
	\[
		\phi_{N-k+1} \wedge \dots \wedge \phi_N \pP = \frac{\psi^{-\frac{k(k-1)}{2}}}{k!} \biggl[\prod_{i<j} \big(\Omega_i^* - \Omega_j^* +\hbar\big) \biggr]\bigl(\phi_N^{\otimes k}\bigr) \pP,
	\]
	where both sides are in $\U(\gl_N) \otimes \Lambda^k \bigl(\C^N\bigr)^* \hamp N_+ $.
	\begin{proof}
		Since $\phi_{N-k+1} \wedge \dots \wedge \phi_N$ is a highest-weight vector, we have
		\[
			\phi_{N-k+1} \wedge \dots \wedge \phi_N \pP = \pP \phi_{N-k+1} \wedge \dots \wedge \phi_N \pP.
		\]
		Let us show by induction on $k$ that
		\begin{align*}
			\pP \phi_{N-k+1} \wedge \dots \wedge \phi_N \pP &{}= {\rm e}^-_k \cdot \pP \phi_{N-k+1} \otimes \dots \otimes \phi_N \pP\\
&{} = \psi^{-\frac{k(k-1)}{2}} ({\rm e}^-_k)^{-\hbar} \bigl[(\Omega^*_1)^{k-1} \cdots \Omega^*_{k-1}\bigr]\bigl(\phi_N^{\otimes k}\bigr),
		\end{align*}
		where ${\rm e}^-_k$ is the antisymmetrizer \eqref{eq::symmetrizer} acting on $\U(\gl_N) \otimes \bigl(\bigl(\C^N\bigr)^*\bigr)^{\otimes k}$, hence on its Kostant--Whittaker reduction. The case $k=1$ is obvious, and for the step, it follows from the proof of Proposition~\ref{prop:dual_kw_presentation} that
		\[
			\pP (\Omega^*)^{k-1} (\phi_N) \pP = \psi^{k-1} \pP \phi_{N-k+1} + \sum_{i<k} c_i \pP \phi_{N-i+1} \pP
		\]
		for some central elements $c_i \in \Z_{\gl_N}$. Therefore, using antisymmetrization and induction assumption, we have
		\[
			{\rm e}^-_k \pP \phi_{N-k+1} \otimes \dots \otimes \phi_N \pP = \psi^{-\frac{k(k-1)}{2}} {\rm e}^-_k \pP (\Omega^*_1)^{k-1} \cdots \Omega^*_{k-1} \bigl(\phi_N^{\otimes k}\bigr) \pP.
		\]
		By $\fH^{-\hbar}_k$-equivariance of Propositions~\ref{prop:dual_kw_presentation} and~\ref{prop:symmetrizer_hecke_vs_shuffle}, the right-hand side is equal to the action of the element
		\[
			\frac{\psi^{-\frac{k(k-1)}{2}}}{k!} \prod_{i<j} \frac{Y_i - Y_j + \hbar}{Y_i - Y_j}\cdot \Sym_k^- \bigl(Y_1^{k-1} \cdots Y_{k-1}^1 Y_{k}^0\bigr) = \frac{\psi^{-\frac{k(k-1)}{2}}}{k!} \prod_{i<j} (Y_i - Y_j + \hbar)
		\]
		by the Vandermonde determinant formula. Therefore,
		\[
			{\rm e}^-_k \cdot (\Omega^*_1)^{k-1} \cdots \Omega^*_{k-1} \bigl(\phi_{N}^{\otimes k}\bigr) \pP = \frac{\psi^{-\frac{k(k-1)}{2}}}{k!} \biggl[\prod_{i<j} \big(\Omega_i^* - \Omega_j^* + \hbar\big) \biggr]\bigl(\phi_N^{\otimes k}\bigr) \pP,
		\]
		as required.
	\end{proof}
\end{Proposition}

Finally, consider the topological setting of Section~\ref{subsect:comparison}. Let
\[
	\Gr^{\omega_1} = \overline{\Gr^{\omega_1}} = \{zL_0 \subset L_{-1} \subset L_0 \mid \dim L_0 / L_{-1} = 1\}.
\]
There is a factorization similar to \eqref{eq::factorization},
\begin{gather}
\label{eq::factorization_dual}
\begin{split}
	& \xymatrix @C=1pc {
		\U\bigl(\gl_N^R\bigr) \otimes \bigl(\C^N\bigr)^* \hamp N_+^R \ar[r] \ar[d]_{\tilde{\varpi}_N} & \lhamp \dD(\GL_N) \hamp N_+^R \ar[r] & \lhamn \dD(\GL_N) \sslash_{\psi^*} N_-^R \ar[d] \\
		\U\bigl(\gl_N^R\bigr) \otimes \C^N \hamn N_-^R \ar[r]^{\sim} & H^{\bullet}_{\GL_N[[t]] \rtimes \C^*} (\Gr^{\omega_1},\underline{\C}_{\Gr^{\omega_1}}) \ar[r] & H^{\GL_N[[t]] \rtimes \C^*}_{\bullet} (\Gr),
	}
\end{split}\hspace{-10mm}
\end{gather}
where $\tilde{\varpi}_N$ is an obvious analog of \eqref{eq::varpi_hc}. As in \eqref{eq::tautological_bundles}, define the bundles
\[
	\cV_i \rightarrow \Gr^{\omega_1}, \qquad \cV_i|_{L_{-1} \subset L_0} = z^{-1} L_i / L_i
\]
as well as their Chern classes
\[
	A^{(i-1)}(u) = \sum_{j=0}^N u^{N-j} A_j^{(i-1)}, \qquad A_j^{(i-1)} = (-1)^i c_j^{\GL_N} (\cV_{i})
\]
for $i=-1,0$. The left (resp.\ the right) action of $H^{\bullet}_G(\mathrm{pt})$ on $H^{\bullet}_{\GL_N[[t]] \rtimes \C^*}(\Gr^{\omega_1})$ is via the Chern classes of $\cV_0$ (resp.\ $\cV_{-1}$).

Define the quotient bundle
\[
	\cQ \rightarrow \Gr^{\omega_1}, \qquad \cQ|_{L_{-1} \subset L_0} = L_0 / L_{-1}.
\]
The following is an analog of Proposition~\ref{prop:top_vs_alg_vector}. Recall the isomorphism of Proposition~\ref{prop:dual_kw_presentation}. Also, since the loop rotation is trivial on $\cQ$, we have
\begin{equation}
\label{eq::chern_class_c_star_dual}
	c_1^G(\cQ) = c_1^{\GL_N} (\cQ).
\end{equation}
\begin{Proposition}
\label{prop:top_vs_alg_dual}
	The assignment $\Omega^* \mapsto c_1^G(\cQ)+\hbar$ defines an isomorphism of $\Z_{\gl_N}$-bimodules
	\[
		\Z_{\gl_N}[\Omega^*]/(A(\Omega^* - \hbar)) \xrightarrow{\sim} \U\bigl(\gl^R_N\bigr) \otimes \bigl(\C^N\bigr)^* \hamp N_+^R \xrightarrow{\sim} H^{\bullet}_{\GL_N[[t]] \rtimes \C^*}(\Gr^{\omega_1})
	\]
	that coincides with the Bezrukavnikov--Finkelberg isomorphism \eqref{eq::factorization_dual}.
	\begin{proof}
		Observe that we can identify
		\[
			\Gr^{\omega_1} \cong \bP\bigl[\bigl(\C^N\bigr)^*\bigr], \qquad [L_{-1} \subset L_0] \mapsto [(L_0/L_{-1})^* \subset (L_0/zL_0)^*].
		\]
		In view of \eqref{eq::center_topological_true}, we have
		\[
			H^{\bullet}_{\GL_N[[t]] \rtimes \C^*}(\Gr^{\omega_1}) \cong \C[\hbar][A_1,\dots,A_N] /(A(c_1(\cQ))).
		\]
		The rest of the proof goes as in Proposition~\ref{prop:top_vs_alg_vector} (recall that the right module structure on the target is given by \eqref{eq::center_conjugation_dual}).
	\end{proof}
\end{Proposition}

\subsection*{Acknowledgements} The author would like to thank Leonid Rybnikov, Joel Kamnitzer, Michael Finkelberg, Alexis Leroux-Lapierre, and Th\'{e}o Pinet for helpful discussions and explanations. We would also like to thank the anonymous referees whose comments greatly improved the paper. This work was supported by Fondation Courtois.

\pdfbookmark[1]{References}{ref}
\LastPageEnding

\end{document}